\documentclass[11pt,twoside,english]{elsarticle}
\usepackage{amssymb,amsmath,amsthm,enumitem}
\usepackage{amsfonts}
\usepackage[T1]{fontenc}
\usepackage[latin9]{inputenc}
\usepackage{geometry}
\usepackage{graphicx}
\usepackage{subfig}
\usepackage{babel}
\usepackage{float}
\usepackage{titlesec}

\titlespacing*{\section}{0pt}{0.3\baselineskip}{0.3\baselineskip}
\titlespacing*{\subsection}{0pt}{0.3\baselineskip}{0.3\baselineskip}
\setlist[itemize]{noitemsep, topsep=0pt}
\makeatletter
\providecommand{\tabularnewline}{\\}
\newtheorem{thm}{Theorem}

\theoremstyle{definition}
\newtheorem{rem}{Remark}
\ifx\proof\undefined
\newenvironment{proof}[1][\protect\proofname]{\par
\normalfont\topsep6\p@\@plus6\p@\relax
\trivlist
\itemindent\parindent
\item[\hskip\labelsep\scshape #1]\ignorespaces
}{
\endtrivlist\@endpefalse
}
\providecommand{\proofname}{Proof}
\fi
\@ifundefined{showcaptionsetup}{}{ \PassOptionsToPackage{caption=false}{subfig}}
\makeatother
\input{tcilatex}
\let\oldbibliography\thebibliography
\renewcommand{\thebibliography}[1]{
\oldbibliography{#1}
\setlength{\itemsep}{0pt}}

\begin{document}

\begin{frontmatter}{}

\title{An inverse problem of a simultaneous reconstruction of the dielectric constant and conductivity from experimental backscattering data}

\tnotetext[t1]{The work of Khoa, Bidney, Klibanov, L. H. Nguyen and Astratov was supported by US Army Research Laboratory and US Army
Research Office grant W911NF-19-1-0044. The work of Khoa was also
partly supported by the Research Foundation-Flanders (FWO) under the
project named ``Approximations for forward and inverse reaction-diffusion
problem related to cancer models''.}

\author[rvt]{Vo Anh Khoa}
\ead{vakhoa.hcmus@gmail.com, avo5@uncc.edu}

\author[rvt1]{Grant W. Bidney}
\ead{gbidney@uncc.edu}

\author[rvt]{Michael V. Klibanov\corref{cor1}}
\ead{mklibanv@uncc.edu}
\cortext[cor1]{Corresponding author.}

\author[rvt]{Loc H. Nguyen}
\ead{loc.nguyen@uncc.edu}

\author[rvt2]{Lam H. Nguyen}
\ead{lam.h.nguyen2.civ@mail.mil}

\author[rvt2]{Anders J. Sullivan}
\ead{anders.j.sullivan.civ@mail.mil}

\author[rvt1]{Vasily N. Astratov}
\ead{astratov@uncc.edu}

\address[rvt]{Department of Mathematics and Statistics, University of North Carolina
at Charlotte, Charlotte, North Carolina 28223, USA.}

\address[rvt1]{Department of Physics and Optical Science, University of North Carolina at Charlotte, Charlotte, North Carolina 28223, USA.}

\address[rvt2]{U.S. Army Research Laboratory, Adelphi, Maryland 20783-1197, USA.}

\begin{abstract}
This report extends our recent progress in tackling a challenging 3D inverse scattering problem governed by the Helmholtz equation. Our target application is to reconstruct 
dielectric constants, electric conductivities and shapes of front surfaces of objects buried very closely under the ground. These objects mimic explosives, like, e.g., antipersonnel land mines and improvised explosive devices. We solve a coefficient inverse problem with the backscattering data generated by a moving source at a fixed frequency. This scenario has been studied so far by our newly developed convexification method that consists in a new derivation of a boundary value problem for a coupled quasilinear elliptic system. However, in our previous work only the unknown dielectric constants of objects and shapes of 
their front surfaces were calculated. Unlike this, in the current work performance of our numerical convexification algorithm is verified for the case when the dielectric constants, the electric conductivities and those shapes
of objects are unknown. By running several tests with experimentally collected backscattering data, we find that we can accurately image both the dielectric constants and shapes of targets of interests including a challenging case of targets with voids. The computed electrical conductivity serves for reliably distinguishing conductive and non-conductive objects. The global convergence of our numerical procedure is shortly revisited.
\end{abstract}

\end{frontmatter}{}

\section{Introduction}

In this paper, we present an extension of our previous numerical studies of
the performance on both computationally simulated \cite{Khoa2019} and
experimental data \cite{Khoa2020} of our newly developed \emph{globally
convergent} \emph{convexification} numerical method for a Coefficient
Inverse Problem (CIP) for the Helmholtz equation in the 3D case. This new
approach has been established so far to solve a 3D CIP with a fixed
frequency and the point source moving along an interval of a straight line.
Our target application is in the reconstruction of physical properties of
explosive-like objects buried closely under the ground, such as
antipersonnel land mines and improvised explosive devices (IEDs). There are
four (4) important properties of such targets that we are currently
interested in:

\begin{enumerate}
\item Dielectric constants.

\item Locations.

\item Shapes of front surfaces.

\item Electrical conductivities.
\end{enumerate}

The previous work \cite{Klibanov2019b} was concerned with the first two
properties for a different set of experimental data and for a different
version of the convexification method. In \cite{Klibanov2019b} the point
source was fixed and the frequency was varied. However, shapes of targets of
interests or, at the very least, their front surfaces were not accurately
imaged in that scenario. This reveals an advantage of the context considered
both in \cite{Khoa2019,Khoa2020} and in this paper, compared with the
previous studies of the research group of the third author.

As we have seen, the dielectric constants and shapes of front surfaces of
the targets with different types of materials and geometries can be
recovered via the convexification from our experimentally collected data 
\cite{Khoa2020}. Thus, this paper is about the recovery of the fourth
property, the electrical conductivity, simultaneously with the three above
ones. We refer to publications of the research group of Novikov for a
variety of CIPs with the data at a single frequency \cite{Ag,Al,N1,N2,N3}.
Both statements of those CIPs and methods of their treatments are different
from ours. Also, we refer to \cite{Bakushinsky2019} for a different
numerical method for a similar CIP.

A conventional way to solve a CIP is numerically is via the minimization of
a conventional least squares cost functional; see, e.g., \cite%
{Chavent,Gonch1,Gonch2}. The major problem with this approach, however, is
that such a functional is non convex and typically suffers from the
phenomenon of multiple local minima and ravines. Since any gradient-like
method stops at any point of a local minimum, then any convergence result
for this method would be valid only if the starting point of iterations
would be located in a small neighborhood of the correct solution. However,
it is unclear how to a priori find such a point. In other words,
conventional numerical methods for CIPs are locally convergent ones.

\textbf{Definition. }\emph{We call a numerical method for a CIP globally
convergent if there is a theorem, which claims that this method delivers at
least one point in a sufficiently small neighborhood of the correct solution
without any advanced knowledge of this neighborhood. }

The convexification is a globally convergent numerical method. It
works with the most challenging case of data collection: our data are both
backscattering and non-overdetermined ones. Given a CIP, we call the data
for it non-overdetermined if the number $m$ of free variables in
the data equals the number $n$ of free variables in the unknown
coefficient, i.e. $m=n$. If, $m>n$, then the data are
overdetermined. In this paper, $m=n=3$. The authors are unaware about such numerical methods for the CIPs  with the non-overdetermined data at $m = n \ge 2$, which would be based on the minimization of a conventional least squares cost functional and, at the same time, would be globally convergent in terms of the above definition.

The latter is the reason why the third author and his collaborators have
been working on the convexification for a number of years; see, e.g., \cite%
{BKconv,Klib95,Klib97,KT} for some first publications in this direction.
After a certain change of variables the convexification constructs a
weighted cost functional $J_{\lambda }$, where $\lambda \geq 1$ is a
parameter. The main element of $J_{\lambda }$ is the presence of the
Carleman Weight Function (CWF) in it. This is the function, which
participates as a weight in the Carleman estimate for the corresponding
partial differential operator. The main theorem then is that, given a convex
bounded set $B\left( d\right) \subset H$ of an arbitrary diameter $d$ in an
appropriate Hilbert space $H,$ one can choose such a set of parameters $%
\lambda $ that the functional $J_{\lambda }$ is strictly convex on $B\left(
d\right) .$ This guarantees, of course, the absence of local minima and
ravines on $B\left( d\right) .$ Using that theorem, the global convergence
to the correct solution of the gradient projection method on $B\left(
d\right) $ is established.

The above cited first works on the convexification lacked some theorems
justifying the latter global convergence property and, therefore, numerical
results in them were not present, except of some first ones in \cite{KT}.
Fortunately, however, the paper \cite{Bakushinsky2017} has cleared that
global convergence issue, which has resulted in a number of recent
publications, where the theory of the convexification is complemented by
numerical studies, see, e.g., \cite%
{Khoa2019,Klibanov2019,Klibanov2019a,Klibhyp}. In particular, experimental
data are treated by the convexification in \cite%
{Khoa2020,convIPnew,Klibanov2019b}. 
In \cite{Baud1,Baud2}, the second version of the convexification is
developed for CIPs for hyperbolic PDEs. Ideas of \cite{Baud1,Baud2} are
explored further in \cite{BBS,LeNg}.

In both versions of the convexification, the ideas of the
Bukhgeim--Klibanov method play the foundational role. This method is based
on Carleman estimates. Initially it was proposed in 1981 only for proofs of
uniqueness theorems for multidimensional CIPs; see \cite{BukhKlib} for the
first publication, and a survey of this method can be found in \cite%
{Klibanov2013}.

This paper is organized as follows. In section 2 we state our CIP. In
section 3 we construct our functional $J_{\lambda }$ and formulate main
theoretical results about it, which we took from our recent work \cite%
{Khoa2020}. Finally, in section 4 we present our numerical results.

\section{Statement of the inverse problem}

Even though the propagation of electromagnetic waves is governed by the
Maxwell's equations, we use only the Helmholtz equation in this paper.
Indeed, it was demonstrated numerically in the appendix of the paper \cite%
{KNN} that this equation describes the propagation of one component of the
electric field equally well with the Maxwell's equations. Another
confirmation comes from our successful work with experimental data, both in 
\cite{Khoa2020,Klibanov2019b} and in this paper.

Let $\mathbf{x}=\left( x,y,z\right) \in \mathbb{R}^{3}$. Prior to the
statement of the inverse problem, we first consider the following
time-harmonic Helmholtz wave equation with conductivity: 
\begin{equation}\label{eq:helm2}
\Delta u+\left( \omega ^{2}\mu \varepsilon ^{\prime }\left( \mathbf{x}%
\right) -\text{i}\omega \mu \sigma \left( \mathbf{x}\right) \right)
u=-\delta \left( \mathbf{x}-\mathbf{x}_{\alpha }\right) \quad \text{in }%
\mathbb{R}^{3},\;\text{i}=\sqrt{-1}.
\end{equation}%
Cf. \cite[Section 3.3]{Balanis2012}, the function $u=u(\mathbf{x},\alpha )$
in (\ref{eq:helm2}) can be physically understood as a component of the
electric field $E=\left( E_{x},E_{y},E_{z}\right) $, which corresponds to a
single nonzero component of the incident field. In our case, this component
is $E_{y}$ (voltage) which is being incident upon the medium. The
backscattering signal of the same component was measured in our experiments.
In (\ref{eq:helm2}), $\omega $ is the angular frequency ($\text{rad}/\text{m}
$), $\mu ,\varepsilon ^{\prime }\left( \mathbf{x}\right) ,\sigma \left( 
\mathbf{x}\right) $ represent respectively the permeability ($\text{H}/\text{%
m}$), permittivity ($\text{F}/\text{m}$) and (effective) conductivity ($%
\text{S}/\text{m}$) of the medium. Furthermore, $\mathbf{x}_{\alpha }$ is
the point source which will be defined below.

Suppose that we only consider non-magnetic targets. Then their relative
permeability has to be unity. With $\varepsilon _{0},\mu _{0}$ being the
vacuum permittivity and vacuum permeability, respectively, and $k=\omega 
\sqrt{\mu _{0}\varepsilon _{0}}$, (\ref{eq:helm2}) becomes 
\begin{equation}\label{eq:helm3}
\Delta u+\left( k^{2}\frac{\mu }{\mu _{0}}\frac{\varepsilon ^{\prime }}{%
\varepsilon _{0}}-\text{i}k\sqrt{\frac{\mu }{\mu _{0}}}\sqrt{\frac{\mu }{%
\varepsilon _{0}}}\sigma \right) u=-\delta \left( \mathbf{x}-\mathbf{x}%
_{\alpha }\right) \quad \text{in }\mathbb{R}^{3}.
\end{equation}%
In (\ref{eq:helm3}), the fraction $\mu /\varepsilon _{0}$ is essentially
close to the so-called characteristic impedance of free space $\mu
_{0}/\varepsilon _{0}=:\eta _{0}$, which is approximately 377 ($\text{S}%
^{-1} $). Denote $c=\varepsilon ^{\prime }/\varepsilon $ as the dielectric
constant. We therefore rewrite the Helmholtz equation (\ref{eq:helm3}) and
impose the Sommerfeld radiation condition: 
\begin{align}
& \Delta u+\left( k^{2}c-\text{i}k\eta _{0}\sigma \right) u=-\delta \left( 
\mathbf{x}-\mathbf{x}_{\alpha }\right) \quad \text{in }\mathbb{R}^{3},
\label{eq:helm} \\
& \lim_{r\rightarrow \infty }r\left( \partial _{r}u-\text{i}ku\right)
=0\quad \text{for }r=\left\vert \mathbf{x}-\mathbf{x}_{\alpha }\right\vert ,%
\text{i}=\sqrt{-1}.  \label{eq:somm}
\end{align}

We now pose the inverse problem. Consider a rectangular prism $\Omega
=\left( -R,R\right) \times \left( -R,R\right) \times \left( -b,b\right) $ in 
$\mathbb{R}^{3}$ for numbers $R,b>0$. This prism our computational domain of
interest. Let $c=c\left( \mathbf{x}\right) $ be the spatially distributed
dielectric constant of the medium. Let $\sigma =\sigma \left( \mathbf{x}%
\right) $ be the spatially distributed electrical conductivity of the
medium. We assume that these functions are smooth and satisfy the following
conditions: 
\begin{equation}
\begin{array}{cc}
\begin{cases}
c\left( \mathbf{x}\right) \geq 1 & \text{in }\Omega , \\ 
c\left( \mathbf{x}\right) =1 & \text{in }\mathbb{R}^{3}\backslash \Omega ,%
\end{cases}
& \quad \text{and }\quad 
\begin{cases}
\sigma \left( \mathbf{x}\right) \geq 0 & \text{in }\Omega , \\ 
\sigma \left( \mathbf{x}\right) =0 & \text{in }\mathbb{R}^{3}\backslash
\Omega .%
\end{cases}%
\end{array}
\label{eq:2}
\end{equation}%
The second line of (\ref{eq:2}) means that we are assuming to have vacuum
outside of the domain of interest $\Omega $. Let the number $d>b$ and let $%
a_{1}<a_{2}$. We consider the line of sources, 
\begin{equation}
L_{\text{src}}:=\left\{ \left( \alpha ,0,-d\right) :a_{1}\leq \alpha \leq
a_{2}\right\} .  \label{Lsrc}
\end{equation}%
This line is parallel to the $x$-axis and is located outside of the closed
domain $\overline{\Omega }$. The distance between $L_{\text{src}}$ and the $%
xy$-plane is $d$, and the length of the line of sources is $\left(
a_{2}-a_{1}\right) $. Using this setting, we arrange for each $\alpha \in %
\left[ a_{1},a_{2}\right] $ the point source $\mathbf{x}_{\alpha }:=\left(
\alpha ,0,-d\right) $ located on the straight line $L_{\text{src}}$. We also
define the near-field measurement site as the lower side of the prism $%
\Omega ,$ 
\begin{equation*}
\Gamma :=\left\{ \mathbf{x}:\left\vert x\right\vert ,\left\vert y\right\vert
<R,z=-b\right\} .
\end{equation*}

To this end, we use either $\alpha $ or $\mathbf{x}_{\alpha }$ to indicate
the dependence of a function/parameter/number on those point sources. We
denote by $u$, $u_{i}$ and $u_{s}$ the total wave, incident wave and
scattered wave, respectively. Also, we note that $u=u_{i}+u_{s}$.

\subsection*{Forward problem}

Given the wavenumber $k>0$ and the functions $c\left( \mathbf{x}\right)
,\sigma \left( \mathbf{x}\right) ,$ the forward problem is to seek the
function $\left. u\left( \mathbf{x},\alpha \right) \right\vert _{\Gamma }$
such that the function $u=u\left( \mathbf{x},\alpha \right) $ is the
solution of problem (\ref{eq:helm})--(\ref{eq:somm}).

Here, the incident wave is 
\begin{equation}
u_{i}\left( \mathbf{x},\alpha \right) =\frac{\text{exp}\left( \text{i}%
k\left\vert \mathbf{x}-\mathbf{x}_{\alpha }\right\vert \right) }{4\pi
\left\vert \mathbf{x}-\mathbf{x}_{\alpha }\right\vert }.  \label{eq:3}
\end{equation}%
Moreover, we can deduce that the scattered wave is: 
\begin{equation}\label{eq:4}
u_{s}\left( \mathbf{x},\alpha \right) =\int_{\Omega }\frac{\text{exp}\left( 
\text{i}k\left\vert \mathbf{x}-\mathbf{x}^{\prime }\right\vert \right) }{%
4\pi \left\vert \mathbf{x}-\mathbf{\ x}^{\prime }\right\vert }\left[
k^{2}\left( c\left( \mathbf{x}^{\prime }\right) -1\right) -\text{i}k\eta
_{0}\sigma \left( \mathbf{x}^{\prime }\right) \right] u\left( \mathbf{x}%
^{\prime },\alpha \right) d\mathbf{x}^{\prime },\quad \text{$\mathbf{x}$}\in 
\mathbb{R}^{3},
\end{equation}%
since functions $c-1$ and $\sigma $ are compactly supported in $\Omega $;
see (\ref{eq:2}). We then combine \eqref{eq:3} and \eqref{eq:4} to obtain
the Lippmann--Schwinger equation \cite{Colton}: 
\begin{equation*}
u\left( \mathbf{x},\alpha \right) =u_{i}\left( \mathbf{x},\alpha \right)
+\int_{\Omega }\frac{\text{exp}\left( \text{i}k\left\vert \mathbf{x}-\mathbf{%
x}^{\prime }\right\vert \right) }{4\pi \left\vert \mathbf{x}-\mathbf{\ x}%
^{\prime }\right\vert }\left[ k^{2}\left( c\left( \mathbf{x}^{\prime
}\right) -1\right) -\text{i}k\eta _{0}\sigma \left( \mathbf{x}^{\prime
}\right) \right] u\left( \mathbf{x}^{\prime },\alpha \right) d\mathbf{x}%
^{\prime },\quad \text{$\mathbf{x}$}\in \mathbb{R}^{3}.
\end{equation*}%
In fact, we generate our data $\left. u\left( \mathbf{x},\alpha \right)
\right\vert _{\Gamma }$ via solving this equation.

\subsection*{Coefficient Inverse Problem (CIP)}

Given $k>0$, the CIP is to reconstruct the two smooth functions: the
dielectric constant $c\left( \mathbf{x}\right) $ and the conductivity $%
\sigma \left( \mathbf{x}\right) $ for $\mathbf{x}\in \Omega $ satisfying
conditions (\ref{eq:2}) from the boundary measurement $F_{0}\left( \mathbf{x}%
,\mathbf{x}_{\alpha }\right) $ of the near-field data, 
\begin{equation}
F_{0}\left( \mathbf{x},\mathbf{x}_{\alpha }\right) =u\left( \mathbf{x}%
,\alpha \right) \quad \text{for }\mathbf{x}\in \Gamma ,\mathbf{x}_{\alpha
}\in L_{\text{src}},  \label{eq:bdr}
\end{equation}%
where $u\left( \mathbf{x},\alpha \right) $ is the total wave associated with
the incident wave $u_{i}$ of \eqref{eq:3}.

\begin{rem}
\label{rem:11} Our CIP is posed for the case of the near-field data,
although the far-field data are collected in our experimental setup.
However, we have observed in our previous works on experimental data that
those far-field data do not look nice \cite{Khoa2020,Nguyen2018}. Thus, we
apply a data propagation procedure, which is described in detail in \cite%
{Khoa2020,Nguyen2018}. This procedure delivers a good approximation for both
the near-field data $F$ (i.e. data at the assigned boundary $\Gamma $) and
the $z-$derivative of the function $u\left( \mathbf{x},\alpha \right) $ at $%
\Gamma ,$%
\begin{equation}
F_{1}\left( \mathbf{x},\mathbf{x}_{\alpha }\right) =u_{z}\left( \mathbf{x}%
,\alpha \right) \quad \text{for }\mathbf{x}\in \Gamma ,\mathbf{x}_{\alpha
}\in L_{\text{src}}.  \label{1}
\end{equation}%
In our experience, the function $F_{0}$ looks much better than the raw
far-field backscattering data \cite{Khoa2020,Nguyen2018}.
\end{rem}

\section{A globally convergent numerical method\label%
{sec:A-globally-convergent}}

\subsection{A system of coupled quasilinear elliptic PDEs}

\label{subsec:3.1}

Since our line of sources $L_{\text{src}}$ is located outside of $\overline{%
\Omega }$, we deduce this system from the homogeneous version of equation %
\eqref{eq:helm} and for each $\alpha \in \left[ a_{1},a_{2}\right] $ 
\begin{equation}
\Delta u+\left( k^{2}c\left( \mathbf{x}\right) -\text{i}k\eta _{0}\sigma
\left( \mathbf{x}\right) \right) u=0\quad \text{in }\Omega .
\label{eq:helm1}
\end{equation}

We set 
\begin{equation*}
\log u_{i}\left( \mathbf{x},\alpha \right) =\text{i}k\left\vert \mathbf{x}-%
\mathbf{x}_{\alpha }\right\vert -\log \left( 4\pi \left\vert \mathbf{x}-%
\mathbf{x}_{\alpha }\right\vert \right) ,
\end{equation*}%
which then leads to 
\begin{equation}
\nabla \left( \log u_{i}\left( \mathbf{x},\alpha \right) \right) =\frac{%
\text{i}k\left( \mathbf{x}-\mathbf{x}_{\alpha }\right) }{\left\vert \mathbf{x%
}-\mathbf{x}_{\alpha }\right\vert }-\frac{\mathbf{x}-\mathbf{x}_{\alpha }}{%
\left\vert \mathbf{x}-\mathbf{x}_{\alpha }\right\vert ^{2}}.
\label{eq:gradui}
\end{equation}%
It was established in \cite{KR} that, under certain conditions, the function 
$u\left( \mathbf{x},\alpha \right) $ is nowhere nonzero at every point in $%
\Omega $ for sufficiently large values of $k$, at least in the case when $%
\sigma \left( \mathbf{x}\right) \equiv 0.$ This, in turn, has allowed in 
\cite{Khoa2019} to uniquely define the function $\log \left( u\left( \mathbf{%
x},\alpha \right) \right) $ for those values of $k$. Thus, we assume below
that we can uniquely define the function $\log \left( u\left( \mathbf{x}%
,\alpha \right) \right) $ as in \cite{Khoa2019}. We note, however, that the
condition $u\left( \mathbf{x},\alpha \right) \neq 0$ is more important for
our derivations below than just defining $\log \left( u\left( \mathbf{x}%
,\alpha \right) \right) $. This is because only derivatives $\nabla \log
\left( u\left( \mathbf{x},\alpha \right) \right) =\nabla u\left( \mathbf{x}%
,\alpha \right) /u\left( \mathbf{x},\alpha \right) $ and $\partial _{\alpha
}\log \left( u\left( \mathbf{x},\alpha \right) \right) =\partial _{\alpha
}u\left( \mathbf{x},\alpha \right) /u\left( \mathbf{x},\alpha \right) $ are
involved below and these derivatives are defined uniquely of course.

Denote $v_{0}\left( \mathbf{x},\alpha \right) =u\left( \mathbf{x},\alpha
\right) /u_{i}\left( \mathbf{x},\alpha \right) $. We define the function $%
v\left( \mathbf{x},\alpha \right) $ as 
\begin{equation*}
v\left( \mathbf{x},\alpha \right) :=\log \left( v_{0}\left( \mathbf{x}%
,\alpha \right) \right) =\log \left( u\left( \mathbf{x},\alpha \right)
\right) -\log \left( u_{i}\left( \mathbf{x},\alpha \right) \right) \quad 
\text{for }\mathbf{x}\in \Omega ,\alpha \in \left[ a_{1},a_{2}\right] .
\end{equation*}%
Hence, one computes that 
\begin{equation}
\nabla v\left( \mathbf{x},\alpha \right) =\frac{\nabla v_{0}\left( \mathbf{x}%
,\alpha \right) }{v_{0}\left( \mathbf{x},\alpha \right) },\quad \Delta
v\left( \mathbf{x},\alpha \right) =\frac{\Delta v_{0}\left( \mathbf{x}%
,\alpha \right) }{v_{0}\left( \mathbf{x},\alpha \right) }-\left( \frac{%
\nabla v_{0}\left( \mathbf{x},\alpha \right) }{v_{0}\left( \mathbf{x},\alpha
\right) }\right) ^{2}.  \label{eq:gradel}
\end{equation}%
Note that $u_{i}(\mathbf{x},\alpha )$ is the fundamental solution of
equation (\ref{eq:helm1}) when $c\equiv 1$ and $\sigma \equiv 0$. We have 
\begin{align}  \label{aaaa}
\Delta u\left( \mathbf{x},\alpha \right) =u_{i}\left( \mathbf{x},\alpha
\right) \Delta v_{0}\left( \mathbf{x},\alpha \right) +2\nabla u_{i}\left( 
\mathbf{x},\alpha \right) \cdot \nabla v_{0}\left( \mathbf{x},\alpha \right)
-v_{0}\left( \mathbf{x},\alpha \right) k^{2}u_{i}\left( \mathbf{x},\alpha
\right) .
\end{align}%
Thus, it follows from (\ref{eq:helm1}) and (\ref{aaaa}) that 
\begin{equation*}
u_{i}\left( \mathbf{x},\alpha \right) \Delta v_{0}\left( \mathbf{x},\alpha
\right) +2\nabla u_{i}\left( \mathbf{x},\alpha \right) \cdot \nabla
v_{0}\left( \mathbf{x},\alpha \right) =-\left[ k^{2}\left( c\left( \mathbf{x}%
\right) -1\right) -\text{i}k\sigma \left( \mathbf{x}\right) \right]
v_{0}\left( \mathbf{x},\alpha \right) u_{i}\left( \mathbf{x},\alpha \right) .
\end{equation*}%
Therefore, combining \eqref{eq:helm1}--\eqref{eq:gradel} we derive the
equation for $v$, 
\begin{equation}
\Delta v+\left( \nabla v\right) ^{2}+2\nabla v\cdot \nabla \left( \log
\left( u_{i}\left( \mathbf{x},\alpha \right) \right) \right) =-\left[
k^{2}\left( c\left( \mathbf{x}\right) -1\right) -\text{i}k\sigma \left( 
\mathbf{x}\right) \right] \quad \text{for }\mathbf{x}\in \Omega .
\label{eq:v}
\end{equation}

We now differentiate \eqref{eq:v} with respect to $\alpha $ and use %
\eqref{eq:gradui} to obtain the following third-order PDE: 
\begin{equation}
\Delta \partial _{\alpha }v+2\nabla v\cdot \nabla \partial _{\alpha
}v+2\nabla \partial _{\alpha }v\cdot \tilde{\mathbf{x}}_{\alpha }+2\hat{%
\mathbf{x}}_{\alpha }\cdot \nabla v=0\quad \text{for }\mathbf{x\in \Omega },
\label{eq:vv}
\end{equation}%
where, for $\mathbf{x}-\mathbf{x}_{\alpha }=\left( x-\alpha ,y,z+d\right)
\in \mathbb{R}^{3},$ 
\begin{align*}
\tilde{\mathbf{x}}_{\alpha }& =\frac{\text{i}k\left( \mathbf{x}-\mathbf{x}%
_{\alpha }\right) }{\left\vert \mathbf{x}-\mathbf{x}_{\alpha }\right\vert }-%
\frac{\mathbf{x}-\mathbf{x}_{\alpha }}{\left\vert \mathbf{x}-\mathbf{x}%
_{\alpha }\right\vert ^{2}}, \\
\hat{\mathbf{x}}_{\alpha }& =\frac{\text{i}k}{\left\vert \mathbf{x}-\mathbf{x%
}_{\alpha }\right\vert ^{3}}\left( -y^{2}-\left( z+d\right) ^{2},\left(
x-\alpha \right) y,\left( x-\alpha \right) z\right) \\
& -\frac{1}{\left\vert \mathbf{x}-\mathbf{x}_{\alpha }\right\vert ^{4}}%
\left( \left( x-\alpha \right) ^{2}-y^{2}-\left( z+d\right) ^{2},2\left(
x-\alpha \right) y,2\left( x-\alpha \right) z\right) .
\end{align*}

\begin{rem}
\label{rem:X} It is obvious that if the function $v\left( \mathbf{x},\alpha
\right) $ is known, then we can immediately find the target coefficients $%
c\left( \mathbf{x}\right) $ and $\sigma \left( \mathbf{x}\right) $ by taking
the real and imaginary parts of the left hand side of (\ref{eq:v}). In order
to get rid of the presence of those unknowns in (\ref{eq:v}), we take
advantage of the moving point source via the differentiation with respect to 
$\alpha $. This leads to solving the nonlinear third order PDE \eqref{eq:vv}%
, which is not an easy task. To overcome this, we apply the Fourier series
approach using a special orthonormal basis with respect to $\alpha $.
\end{rem}

For $\alpha \in \left( a_{1},a_{2}\right) $, let $\left\{ \Psi _{n}\left(
\alpha \right) \right\} _{n=0}^{\infty }$ be the special orthonormal basis
in $L^{2}\left( a_{1},a_{2}\right) $, which was first proposed in \cite%
{Klibanov2017}. Herewith, the construction of this basis is shortly
revisited. For each $n\in \mathbb{N}$, let $\varphi _{n}\left( \alpha
\right) =\alpha ^{n}e^{\alpha }$ for $\alpha \in \left[ a_{1},a_{2}\right] $%
. The set $\left\{ \varphi _{n}\left( \alpha \right) \right\} _{n=0}^{\infty
}$ is linearly independent and complete in $L^{2}\left( a_{1},a_{2}\right) $%
. Using the Gram--Schmidt orthonormalization procedure, we can obtain the
orthonormal basis $\left\{ \Psi _{n}\left( \alpha \right) \right\}
_{n=0}^{\infty }$ in $L^{2}\left( a_{1},a_{2}\right) $, which possesses the
following special properties:

\begin{itemize}
\item $\Psi_{n}\in C^{\infty}\left[a_{1},a_{2}\right]$ for all $n\in\mathbb{N%
}$;

\item Let $s_{mn}=\left\langle \Psi _{n}^{\prime },\Psi _{m}\right\rangle $,
where $\left\langle \cdot ,\cdot \right\rangle $ denotes the scalar product
in $L^{2}\left( a_{1},a_{2}\right) $. Then the square matrix $S_{N}=\left(
s_{mn}\right) _{m,n=0}^{N-1}$ for $N\in \mathbb{N}$ is invertible with $%
s_{mn}=1$ if $m=n$ and $s_{mn}=0$ if $n<m$.
\end{itemize}

We note that neither classical orthogonal polynomials nor the classical
basis of trigonometric functions do not hold the second property. The matrix 
$S_{N}$ is an upper diagonal matrix with $\det \left( S_{N}\right) =1$. On
the other hand, the special second property allows us to reduce the
third-order PDE \eqref{eq:vv} to a system of coupled elliptic PDEs.

Given $N\in \mathbb{N}$, we consider the following truncated Fourier series
for $v$: 
\begin{equation}
v\left( \mathbf{x},\alpha \right) =\sum_{n=0}^{N-1}\left\langle v\left( 
\mathbf{x},\cdot \right) ,\Psi _{n}\left( \cdot \right) \right\rangle \Psi
_{n}\left( \alpha \right) \quad \text{for }\mathbf{x}\in \Omega ,\alpha \in %
\left[ a_{1},a_{2}\right] .  \label{eq:fourier}
\end{equation}%
As in \cite{Khoa2019,Khoa2020}, we substitute \eqref{eq:fourier} into %
\eqref{eq:vv}, multiplying both side of the resulting equation by $\Psi
_{m}\left( \alpha \right) $ for $0\leq m\leq N-1$ and then taking the
integration with respect to $\alpha $, we obtain the following system of
coupled elliptic equations: 
\begin{align}
& \Delta V\left( \mathbf{x}\right) +K\left( \nabla V\left( \mathbf{x}\right)
\right) =0\quad \text{for }\mathbf{x}\in \Omega ,  \label{eq:pdeV} \\
& \partial _{\nu }V\left( \mathbf{x}\right) =0\quad \text{for }\mathbf{x}\in
\partial \Omega \backslash \Gamma ,  \label{eq:bdrV1} \\
& V\left( \mathbf{x}\right) =\psi _{0}\left( \mathbf{x}\right) ,V_{z}\left( 
\mathbf{x}\right) =\psi _{1}\left( \mathbf{x}\right) \quad \text{for }%
\mathbf{x}\in \Gamma ,  \label{eq:bdrV2}
\end{align}%
where $\nu =\nu \left( \mathbf{x}\right) $ is the outward looking normal
vector on $\partial \Omega \backslash \Gamma .$ Here, $V\left( \mathbf{x}%
\right) \in \mathbb{R}^{N}$ is the unknown vector function given by 
\begin{equation}
V\left( \mathbf{x}\right) =%
\begin{pmatrix}
v_{0}\left( \mathbf{x}\right) & v_{1}\left( \mathbf{x}\right) & \cdots & 
v_{N-1}\left( \mathbf{x}\right)%
\end{pmatrix}%
^{T}.  \label{3}
\end{equation}%
Boundary conditions (\ref{eq:bdrV2}) are Cauchy data and are overdetermined
ones. A Lipschitz stability estimate for problem (\ref{eq:pdeV})--(\ref%
{eq:bdrV2}) is obtained in \cite{Khoa2019} using a Carleman estimate with
the same CWF as we work with below.

In PDE \eqref{eq:pdeV}, we denote $K\left( \nabla V\left( \mathbf{x}\right)
\right) =S_{N}^{-1}f\left( \nabla V\left( \mathbf{x}\right) \right) $, where 
$f=\left( \left( f_{m}\right) _{m=0}^{N-1}\right) ^{T}\in \mathbb{R}^{N}$ is
quadratic with respect to the first derivative of components of $V\left( 
\mathbf{x}\right) $, 
\begin{align}
f_{m}\left( \nabla V\left( \mathbf{x}\right) \right) &
=2\sum_{n,l=0}^{N-1}\nabla v_{n}\left( \mathbf{x}\right) \cdot \nabla
v_{l}\left( \mathbf{x}\right) \int_{a_{1}}^{a_{2}}\Psi _{m}\left( \alpha
\right) \Psi _{n}\left( \alpha \right) \Psi _{l}^{\prime }\left( \alpha
\right) d\alpha  \label{eq:17} \\
& +2\sum_{n=0}^{N-1}\int_{a_{1}}^{a_{2}}\Psi _{m}\left( \alpha \right) \Psi
_{n}^{\prime }\left( \alpha \right) \nabla v_{n}\left( \mathbf{x}\right)
\cdot \tilde{\mathbf{x}}_{\alpha }d\alpha
+2\sum_{n=0}^{N-1}\int_{a_{1}}^{a_{2}}\Psi _{m}\left( \alpha \right) \Psi
_{n}\left( \alpha \right) \hat{\mathbf{x}}_{\alpha }\cdot \nabla v_{n}\left( 
\mathbf{x}\right) d\alpha .  \notag
\end{align}

\begin{rem}
In \eqref{eq:bdrV2}, the boundary information $\psi _{0}$ at $\Gamma $ is a
direct application of the expansion \eqref{eq:fourier} to the near field
data \eqref{eq:bdr}. As mentioned in Remark \ref{rem:11}, using the data
propagation technique will result in an approximation of the $z$--derivative 
$F_{1}$ of the function \ u at $\Gamma ,$ see (\ref{1}). The function $F_{1}$
generates the Neumann boundary condition $\psi _{1}$ in (\ref{eq:bdrV2}) at $%
\Gamma $. It is shown in \cite{Khoa2020} that the zero Neumann boundary
condition \eqref{eq:bdrV1} at $\partial \Omega \backslash \Gamma $ follows
from an approximation of the radiation condition (\ref{eq:somm}).
\end{rem}

\begin{rem}
From now on we work only within the framework of an approximate mathematical
model. This means that we fix the number $N$ of Fourier harmonics in the
truncated Fourier series (\ref{eq:fourier}) and ignore the residual in (\ref%
{eq:pdeV}). Furthermore, we cannot prove convergence of our method at $%
N\rightarrow \infty $ since such a proof is an extremely challenging
problem. We refer to \cite{Khoa2019} for a detailed discussion of this
issue. We also note that similar approximate mathematical models without
proofs of convergence at $N\rightarrow \infty $ are actually used quite
often in the field of inverse problems. In this regard we refer to \cite%
{Ag,Al,Kab1,Kab2,N1,N2,N3}.
\end{rem}

\subsection{Weighted cost functional in partial finite differences}

We set 
\begin{equation}
L\left( V\right) \left( \mathbf{x}\right) =\Delta V\left( \mathbf{x}\right)
+K\left( \nabla V\left( \mathbf{x}\right) \right) .  \label{eq:LL}
\end{equation}%
Let the numbers $\theta >b$ and $\lambda \geq 1$. We define our CWF as 
\begin{equation}
\mu _{\lambda }\left( z\right) =\exp \left( 2\lambda \left( z-\theta \right)
^{2}\right) \quad \text{for }z\in \left[ -b,b\right] .  \label{eq:CWF}
\end{equation}%
The choice $\theta >b$ is based on the fact that the gradient of the CWF
should not vanish in the closed domain $\overline{\Omega }$. Our CWF is
decreasing for $z\in \left( -b,b\right) $ and 
\begin{equation*}
\max_{z\in \left[ -b,b\right] }\mu _{\lambda }\left( z\right) =\mu _{\lambda
}\left( -b\right) =e^{2\lambda \left( b+\theta \right) ^{2}},\quad
\min_{z\in \left[ -b,b\right] }\mu _{\lambda }\left( z\right) =\mu _{\lambda
}\left( b\right) =e^{2\lambda \left( b-\theta \right) ^{2}}.
\end{equation*}%
Therefore, the CWF \eqref{eq:CWF} attains its maximal value on the site $%
\Gamma $, and it attains its minimal value on the opposite side. By this
way, it \textquotedblleft maximizes\textquotedblright\ the influence of the
measured boundary data at $z=-b$. Furthermore, the notion behind this use of
the CWF is to \textquotedblleft convexify\textquotedblright\ the cost
functional globally and, especially, to control the nonlinear term $K\left(
\nabla V\left( \mathbf{x}\right) \right) $ in \eqref{eq:pdeV}.

\subsubsection{Abstract setting and notation}

While the convergence analysis in \cite{Khoa2019} was done for the case when
derivatives in (\ref{eq:pdeV})--(\ref{eq:bdrV2}), (\ref{eq:LL}) are
considered in their regular continous seting, in \cite{Khoa2020}, so as in
the current paper, we consider the convergence analysis for the case when (%
\ref{eq:pdeV})--(\ref{eq:bdrV2}), (\ref{eq:LL}) are written via
\textquotedblleft partial finite differences\textquotedblright . This means
that we use finite differences with respect the variables $x,y$ and keep the
standard derivatives with respect to $z$. It was pointed out in \cite%
{Khoa2020} that an important reason of doing so is that partial finite
differences allow us not to use the penalty regularization term in our
weighted cost functional, which is convenient for computations. When doing
so, we use the same grid step size $h$ in $x$ and $y$ directions and do not
allow $h$ to tend to zero. Note that it is impractical to use too small grid
step sizes in computations. Consider two partitions of the interval $\left[
-R,R\right] ,$ 
\begin{align*}
& -R=x_{0}<x_{1}<\ldots <x_{Z_{h}-1}<x_{Z_{h}}=R,\quad x_{p}-x_{p-1}=h, \\
& -R=y_{0}<y_{1}<\ldots <y_{Z_{h}-1}<y_{Z_{h}}=R,\quad y_{q}-y_{q-1}=h.
\end{align*}%
Then, for any $N$-dimensional vector function $u\left( \mathbf{x}\right) $,
we denote by $u_{p,q}^{h}\left( z\right) =u\left( x_{p},y_{q},z\right) $ the
corresponding semi-discrete function defined at grid points $\left\{ \left(
x_{p},y_{q}\right) \right\} _{p,q=0}^{Z_{h}}$. Thus, the interior grid
points are $\left\{ \left( x_{p},y_{q}\right) \right\} _{p,q=1}^{Z_{h}-1}$.
Denote 
\begin{equation*}
\Omega _{h}=\left\{ \left( x_{p},y_{q},z\right) :\left\{ \left(
x_{p},y_{q}\right) \right\} _{p,q=0}^{Z_{h}-1}\subset \left[ -R,R\right]
\times \left[ -R,R\right] ,z\in \left( -b,b\right) \right\} ,
\end{equation*}%
\begin{equation*}
\Gamma _{h}=\left\{ \left( x_{p},y_{q},-b\right) :\left\{ \left(
x_{p},y_{q}\right) \right\} _{p,q=0}^{Z_{h}-1}\subset \left[ -R,R\right]
\times \left[ -R,R\right] \right\} .
\end{equation*}%
Henceforth, the corresponding Laplace operator in partial finite differences
is given by $\Delta ^{h}u^{h}=u_{zz}^{h}+u_{xx}^{h}+u_{yy}^{h}$, where, for
interior points of $\Omega _{h}$, we use 
\begin{equation*}
u_{xx}^{h}=h^{-2}\left( u_{p+1,q}^{h}\left( z\right) -2u_{p,q}^{h}\left(
z\right) +u_{p-1,q}^{h}\left( z\right) \right) ,\text{ }p,q\in \left[
1,Z_{h}-1\right]
\end{equation*}%
and similarly for $u_{yy}^{h}$. We define for interior points $\nabla
^{h}u_{p,q}\left( z\right) =\left( \partial _{x}^{h}u_{p,q}\left( z\right)
,\partial _{y}^{h}u_{p,q}\left( z\right) ,\partial _{z}u_{p,q}^{h}\left(
z\right) \right) $. Here, $\partial _{x}^{h}u_{p,q}^{h}\left( z\right)
=\left( 2h\right) ^{-1}\left( u_{p+1,q}^{h}\left( z\right)
-u_{p-1,q}^{h}\left( z\right) \right) $. Hence, the differential operator %
\eqref{eq:LL} has the following form in the partial finite differences: 
\begin{equation}
L^{h}\left( V^{h}\left( z\right) \right) =\Delta ^{h}V^{h}\left( z\right)
+K\left( \nabla ^{h}V^{h}\left( z\right) \right) .  \label{98}
\end{equation}

To simplify the presentation, we consider any $N$--D complex valued function 
$W=\func{Re}W+$i$\func{Im}W$ as the $2N$--D vector function with real valued
components $\left( \func{Re}W,\func{Im}W\right) :=\left( W_{1},W_{2}\right)
:=W\in \mathbb{R}^{2N}$.

Denote $w^{h}\left( z\right) $ the vector function $w^{h}\left( z\right)
=\left\{ w_{p,q}\left( z\right) \right\} _{p,q=0}^{Z_{h}},$ where $%
w_{p,q}^{h}\left( z\right) =w\left( x_{p},y_{q},z\right) .$ Next, we
consider $2N-$D vector functions 
\begin{equation}
W^{h}\left( z\right) =\left( w_{0,1}^{h}\left( z\right) ,w_{0,2}^{h}\left(
z\right) ,w_{1,1}^{h}\left( z\right) ,w_{1,2}^{h}\left( z\right)
,...,w_{N-1,1}^{h}\left( z\right) ,w_{N-1,2}^{h}\left( z\right) \right) ^{T},
\label{4}
\end{equation}%
where $w_{s,1}^{h}\left( z\right) =\func{Re}w_{s}^{h}\left( z\right) $ and $%
w_{s,2}^{h}\left( z\right) =\func{Im}w_{s}^{h}\left( z\right) .$ In
notations of spaces below the subscript \textquotedblleft $2N$'' means that
this space consists of such vector functions. As in \cite{Khoa2020}, we
consider the Hilbert spaces $H_{2N}^{2,h}=H_{2N}^{2,h}\left( \Omega
_{h}\right) $ and $L_{2N}^{2,h}=L_{2N}^{2,h}\left( \Omega _{h}\right) $ of
semi-discrete real valued functions: 
\begin{align*}
& H_{2N}^{2,h}=\left\{ W^{h}\left( z\right) :\left\Vert W^{h}\right\Vert
_{H_{2N}^{2,h}}^{2}:=\dsum\limits_{s=1}^{N-1}\dsum\limits_{j=1}^{2}%
\sum_{p,q=1}^{Z_{h}-1}\sum_{m=0}^{2}h^{2}\int_{-b}^{b}\left\vert \partial
_{z}^{m}w_{s,j,p,q}^{h}\left( z\right) \right\vert ^{2}dz<\infty \right\} ,
\\
& L_{2N}^{2,h}=\left\{ w^{h}\left( z\right) :\left\Vert w^{h}\right\Vert
_{L_{2N}^{2,h}}^{2}:=\dsum\limits_{s=1}^{N-1}\dsum\limits_{j=1}^{2}%
\sum_{p,q=1}^{Z_{h}-1}h^{2}\int_{-b}^{b}\left\vert w_{s,j,p,q}^{h}\left(
z\right) \right\vert ^{2}dz<\infty \right\} .
\end{align*}%
Denote $\left( \cdot ,\cdot \right) $ the scalar product in the space $%
H_{2N}^{2,h}\left( \Omega _{h}\right) $. The subspace $H_{2N,0}^{2,h}\subset
H_{2N}^{2,h}$ is defined as 
\begin{equation*}
H_{2N,0}^{2,h}=\left\{ w^{h}\left( z\right) \in H_{2N}^{2,h}:\left. \nabla
^{h}w_{p,q}^{h}\left( z\right) \right\vert _{\partial \Omega _{h}\backslash
\Gamma _{h}}\cdot \nu =0,\left. w_{p,q}^{h}\left( z\right) \right\vert
_{\Gamma _{h}}=\left. \partial _{z}w_{p,q}^{h}\left( z\right) \right\vert
_{\Gamma _{h}}=0\right\} .
\end{equation*}%
Let $h_{0}>0$ be a fixed positive number. We assume below that 
\begin{equation}
h\geq h_{0}>0.  \label{190}
\end{equation}%
For an arbitrary $M>0$ let the sets $B\left( M\right) \subset
H_{2N}^{2,h}\left( \Omega _{h}\right) $ and $B_{0}\left( M\right)
\subset H_{2N,0}^{2,h}$ be: 
\begin{equation}
B\left( M\right) :=\left\{ 
V^{h}\in H_{2N}^{2,h}:\left\Vert V^{h}\right\Vert _{H_{2N}^{2,h}}<M,\text{ }%
\left. \nabla ^{h}V^{h}\right\vert _{\partial \Omega _{h}\backslash \Gamma
_{h}}\cdot \nu =0,
\left. V^{h}\right\vert _{\Gamma _{h}}=\psi _{0}^{h},\text{ }\left. \partial
_{z}V^{h}\right\vert _{\Gamma _{h}}=\psi _{1}^{h}%
\right\} ,  \label{eq:BM}
\end{equation}%
\begin{equation}
B_{0}\left( M\right) =\left\{ V^{h}\in H_{2N,0}^{2,h}:\left\Vert
V^{h}\right\Vert _{H_{2N}^{2,h}}<M\right\} .  \label{100}
\end{equation}

\subsubsection{Minimization problem and convergence results}

Under the partial finite differences setting, we seek an approximate
solution of system \eqref{eq:pdeV}--\eqref{eq:bdrV2} using the minimization
of the following weighted Tikhonov-type functional. We define the cost
functional $J_{h,\lambda }:H_{2N}^{2,h}\left( \Omega _{h}\right) \rightarrow 
\mathbb{R}_{+}$ as follows: 
\begin{equation}
J_{h,\lambda }\left( V^{h}\right)
=\sum_{p,q=1}^{Z_{h}-1}h^{2}\int_{-b}^{b}\left\vert L^{h}\left( V^{h}\left(
z\right) \right) \right\vert ^{2}\mu _{\lambda }\left( z\right) dz,
\label{eq:J}
\end{equation}%
In (\ref{eq:J}), we have denoted%
\begin{equation*}
V^{h}\left( z\right) =\left( v_{0,1}^{h}\left( z\right) ,v_{0,2}^{h}\left(
z\right) ,v_{1,1}^{h}\left( z\right) ,v_{1,2}^{h}\left( z\right)
,...,v_{N-1,1}^{h}\left( z\right) ,v_{N-1,2}^{h}\left( z\right) \right) ^{T},
\end{equation*}
see (\ref{3}) and (\ref{4}). Also, the CWF $\mu _{\lambda }\left( z\right) $
is defined in \eqref{eq:CWF}, and the operator $L^{h}\left( V^{h}\left(
z\right) \right) $ is as in (\ref{98}). The minimization problem is
formulated as:

\emph{Minimize the cost functional }$J_{\lambda }\left( V^{h}\right) $\emph{%
\ on the set }$\overline{B\left( M\right) }$\emph{, where the set }$B\left(
M\right) $ \emph{is defined in (\ref{eq:BM}).}

We now formulate theorems of our convergence analysis. Those theoretical
results were proven in \cite{Khoa2020}. It is worth mentioning that our
results below are valid only under condition (\ref{190}). Restricting from
below our grid step size $h$ by a constant $h_{0}>0$ is well-suited to our
numerical context. In fact, it is impractical to use too fine
discretization. We begin with the Carleman estimate for the Laplace operator
in partial finite differences.

\begin{thm}[Carleman estimate in partial finite differences]
\label{thm:carle}There exist a sufficient large constant $\lambda
_{0}=\lambda _{0}\left( b,\theta ,h_{0}\right) \geq 1$ and a number $%
C=C\left( b,\theta ,h_{0}\right) >0$ depending only on numbers $b,\theta
,h_{0}$ such that for all $\lambda \geq \lambda _{0}$ and for all vector
functions $u^{h}\in H_{2N,0}^{2,h}\left( \Omega _{h}\right) $ the following
Carleman estimate holds: 
\begin{align}
& \sum_{p,q=1}^{Z_{h}-1}h^{2}\int_{-b}^{b}\left( \Delta
^{h}u_{p,q}^{h}\left( z\right) \right) ^{2}\mu _{\lambda }\left( z\right)
dz\geq C\sum_{p,q=1}^{Z_{h}-1}h^{2}\int_{-b}^{b}\left( \partial
_{z}^{2}u_{p,q}^{h}\left( z\right) \right) ^{2}\mu _{\lambda }\left(
z\right) dz  \label{eq:Carleman} \\
& +C\lambda \sum_{p,q=1}^{Z_{h}-1}h^{2}\int_{-b}^{b}\left( \partial
_{z}u_{p,q}^{h}\left( z\right) \right) ^{2}\mu _{\lambda }\left( z\right)
dz+C\lambda ^{3}\sum_{p,q=1}^{Z_{h}-1}h^{2}\int_{-b}^{b}\left[ \left( \nabla
^{h}u_{p,q}^{h}\left( z\right) \right) ^{2}+\left( u_{p,q}^{h}\left(
z\right) \right) ^{2}\right] \mu _{\lambda }\left( z\right) dz.  \notag
\end{align}
\end{thm}

The next theorem is about the global strict convexity of the cost functional 
$J_{h,\lambda }$.

\begin{thm}[Global strict convexity: the central theorem]
\label{thm:convex}For any $\lambda >0$ the functional $J_{h,\lambda }\left(
V^{h}\right) $ defined in \eqref{eq:J} has its Fréchet derivative $%
J_{h,\lambda }^{\prime }\left( V^{h}\right) \in H_{2N,0}^{2,h}$ at any point 
$V^{h}\in \overline{B\left( M\right) }$. Let $\lambda _{0}>1$ be the number
of Theorem \ref{thm:carle}. There exist numbers $\lambda _{1}=\lambda
_{1}\left( b,\theta ,h_{0},N,M\right) \geq \lambda _{0}>1$ and $%
C_{1}=C_{1}\left( b,\theta ,h_{0},N,M\right) >0$ depending only on listed
parameters such that for all $\lambda \geq \lambda _{1}$ the functional $%
J_{h,\lambda }\left( V^{h}\right) $ is strictly convex on the set $\overline{%
B\left( M\right) }$. More precisely, the following estimate holds: 
\begin{equation}
J_{h,\lambda }\left( V^{h}+r^{h}\right) -J_{h,\lambda }\left( V^{h}\right)
-J_{h,\lambda }^{\prime }\left( V^{h}\right) \left( r^{h}\right) \geq
C_{1}e^{2\lambda \left( b-\theta \right) ^{2}}\left\Vert r^{h}\right\Vert
_{H_{2N}^{2,h}}^{2}\quad \text{for all }V^{h},V^{h}+r^{h}\in \overline{%
B\left( M\right) }.  \label{eq:convex}
\end{equation}
\end{thm}

Below $C_{1}=C_{1}\left( b,\theta ,h_{0},N,M\right) >0$ denotes different
numbers depending only on listed parameters. We now formulate a theorem
about the Lipschitz continuity of the Fréchet derivative $J_{h,\lambda
}^{\prime }\left( V^{h}\right) $ on the set $\overline{B\left( M\right) }$.
We omit the proof of this result because it is similar to the proof of
Theorem 3.1 in \cite{Bakushinsky2017}.

\begin{thm}
\label{thm:lip}The Lipschitz continuity on the set $\overline{B\left(
M\right) }$ holds for the Fréchet derivative $J_{h,\lambda }^{\prime }\left(
V^{h}\right) .$ Namely, there exists a number $\widehat{C}=\widehat{C}\left(
b,\theta ,h_{0},N,M,\lambda \right) >0$ depending only on numbers $b,\theta
,h_{0},N,M,\lambda $ such that for any pair $V_{\left( 1\right)
}^{h},V_{\left( 2\right) }^{h}\in \overline{B\left( M\right) }$ the
following estimate is valid: 
\begin{equation*}
\left\Vert J_{h,\lambda }^{\prime }\left( V_{\left( 2\right) }^{h}\right)
-J_{h,\lambda }^{\prime }\left( V_{\left( 1\right) }^{h}\right) \right\Vert
_{H_{2N}^{2,h}}\leq \tilde{C}\left\Vert V_{\left( 2\right) }^{h}-V_{\left(
1\right) }^{h}\right\Vert _{H_{2N}^{2,h}}.
\end{equation*}
\end{thm}

The following theorem establishes the existence and uniqueness of the
minimizer of the functional $J_{h,\lambda }$ on the set $\overline{B\left(
M\right) }.$ This theorem follows from a combination of Theorems \ref%
{thm:convex} and \ref{thm:lip} with Lemma 2.1 and Theorem 2.1 of \cite%
{Bakushinsky2017}.

\begin{thm}
\label{thm:6} Let $\lambda _{1}>1$ be the number of Theorem \ref{thm:convex}%
. For any $\lambda \geq \lambda _{1}$ there exists a unique minimizer $%
Q_{\min ,\lambda }^{h}\in \overline{B\left( M\right) }$ of the functional $%
J_{h,\lambda }\left( V^{h}\right) $ on the set $\overline{B\left( M\right) }$%
. In addition, 
\begin{equation}
\left( J_{h,\lambda }^{\prime }\left( Q_{\min ,\lambda }^{h}\right) ,Q_{\min
,\lambda }^{h}-S\right) \leq 0\quad \text{for all }S\in \overline{B\left(
M\right) }.  \label{3.0}
\end{equation}
\end{thm}

In the regularization theory, minimizers of the functional $J_{h,\lambda
}\left( V^{h}\right) $ for different levels of the noise in the data are
called ``regularized solutions'' \cite{Beilina2012,Tikhonov1995}. We now
establish the accuracy estimate of our regularized solutions depending on
the noise level in the data. Using (\ref{98}), we obtain the following
analog of problem (\ref{eq:pdeV})--\eqref{eq:bdrV2} in partial finite
differences:%
\begin{align}
& L^{h}\left( V^{h}\left( z\right) \right) =\Delta ^{h}V^{h}\left( \mathbf{x}%
^{h}\right) +K\left( \nabla ^{h}V^{h}\left( \mathbf{x}^{h}\right) \right)
=0\quad \text{for }\mathbf{x}^{h}\in \Omega _{h},  \label{3.1} \\
& \nabla ^{h}V^{h}\left( \mathbf{x}^{h}\right) \cdot \nu \left( \mathbf{x}%
^{h}\right) =0\quad \text{ for }\mathbf{x}^{h}\in \partial \Omega
_{h}\backslash \Gamma _{h},  \label{3.2} \\
& V\left( \mathbf{x}^{h}\right) =\psi _{0}^{h}\left( \mathbf{x}^{h}\right)
,\partial _{z}V\left( \mathbf{x}^{h}\right) =\psi _{1}\left( \mathbf{x}%
^{h}\right) \quad \text{for }\mathbf{x}^{h}\in \Gamma _{h}.  \label{3.3}
\end{align}%
Following the Tikhonov regularization concept \cite{Beilina2012,Tikhonov1995}%
, we assume that there exists an exact solution $V_{\ast }^{h}\in $ $%
H_{2N}^{2,h}$ of problem (\ref{3.1})--(\ref{3.3}) with the noiseless data $%
\psi _{0\ast }^{h}\left( \mathbf{x}^{h}\right) $ and $\psi _{1\ast
}^{h}\left( \mathbf{x}^{h}\right) $. The subscript \textquotedblleft $\ast $%
\textquotedblright\ is used only for the exact solution. Note that the data $%
\psi _{0}^{h},\psi _{1}^{h}$ are always noisy and we denote the noise level
by $\delta \in \left( 0,1\right) $. Besides, we assume that 
\begin{equation}
\left\Vert V_{\ast }^{h}\right\Vert _{H_{2N}^{2,h}}<M-\delta .  \label{3.40}
\end{equation}%
We assume that there exist two vector functions $\Psi _{\ast }^{h},\Psi
^{h}\in H_{2N}^{2,h}$ such that 
\begin{align}
& \nabla ^{h}\Psi _{\ast }^{h}\left( \mathbf{x}^{h}\right) \cdot \nu \left( 
\mathbf{x}^{h}\right) =0,\quad \nabla ^{h}\Psi ^{h}\left( \mathbf{x}%
^{h}\right) \cdot \nu \left( \mathbf{x}^{h}\right) =0\quad \text{for }%
\mathbf{x}^{h}\in \partial \Omega _{h}\backslash \Gamma _{h},  \label{3.4} \\
& \Psi _{\ast }^{h}\left( \mathbf{x}^{h}\right) =\psi _{0\ast }^{h}\left( 
\mathbf{x}^{h}\right) ,\quad \partial _{z}\Psi _{\ast }^{h}\left( \mathbf{x}%
^{h}\right) =\psi _{1\ast }\left( \mathbf{x}^{h}\right) \quad \text{for }%
\mathbf{x}^{h}\in \Gamma _{h},  \label{3.5} \\
& \Psi ^{h}\left( \mathbf{x}^{h}\right) =\psi _{0}^{h}\left( \mathbf{x}%
^{h}\right) ,\quad \partial _{z}\Psi ^{h}\left( \mathbf{x}^{h}\right) =\psi
_{1}\left( \mathbf{x}^{h}\right) \quad \text{for }\mathbf{x}^{h}\in \Gamma
_{h},  \label{3.6} \\
& \left\Vert \Psi _{\ast }^{h}\right\Vert _{H_{2N}^{2,h}}<M,\quad \left\Vert
\Psi ^{h}\right\Vert _{H_{2N}^{2,h}}<M,  \label{3.7} \\
& \left\Vert \Psi ^{h}-\Psi _{\ast }^{h}\right\Vert _{H_{2N}^{2,h}}<\delta .
\label{3.8}
\end{align}%
The following theorem provides the accuracy estimate of the minimizer $%
V_{\min ,\lambda }^{h}$.

\begin{thm}[accuracy estimate of regularized solutions]
\label{thm:7} Assume that conditions (\ref{3.40})--(\ref{3.8}) are valid.
Let $\lambda _{1}=\lambda _{1}\left( b,\theta ,h_{0},N,M\right) >1$ be the
number of Theorem \ref{thm:lip}. Let $V_{\min ,\lambda }^{h}\in \overline{%
B\left( M\right) }$ be the minimizer of the functional (\ref{eq:J}), which
is found in Theorem \ref{thm:6}. Then the following accuracy estimate holds
for all $\lambda \geq \lambda _{1}$ 
\begin{equation}
\left\Vert V_{\min ,\lambda }^{h}-V_{\ast }^{h}\right\Vert
_{H_{2N}^{2,h}}\leq C_{1}e^{4\lambda b\theta }\delta .  \label{3.80}
\end{equation}
\end{thm}

For each $V^{h}\in B\left( M\right) $ consider the vector function $%
W^{h}=V^{h}-\Psi ^{h}.$ Then (\ref{3.7}) implies that 
\begin{align}
& W^{h}\in B_{0}\left( 2M\right) \subset H_{2N,0}^{2,h}\quad \text{for all }%
V^{h}\in B\left( M\right) ,  \label{3.17} \\
& W^{h}+\Psi ^{h}\in B\left( 3M\right) \quad \text{for all }W^{h}\in
B_{0}\left( 2M\right) ,  \label{3.18}
\end{align}%
see (\ref{100}) for $B_{0}\left( M\right) .$ Consider the functional $%
I_{h,\lambda }:B_{0}\left( 2M\right) \rightarrow \mathbb{R}$ defined as%
\begin{equation}
I_{h,\lambda }\left( W^{h}\right) =J_{h,\lambda }\left( W^{h}+\Psi
^{h}\right) \quad \text{for all }W^{h}\in B_{0}\left( 2M\right) .
\label{3.19}
\end{equation}%
Theorem \ref{thm:8} follows immediately from Theorems \ref{thm:lip}--\ref%
{thm:7} and (\ref{3.17})--(\ref{3.19}).

\begin{thm}
\label{thm:8} For any $\lambda >0$ the functional $I_{h,\lambda }\left(
W^{h}\right) $ has its Fréchet derivative $I_{h,\lambda }^{\prime }\left(
W^{h}\right) \in H_{2N,0}^{2,h}$ at any point $W^{h}\in \overline{%
B_{0}\left( 2M\right) }$ and this derivative is Lipschitz continuous on $%
\overline{B_{0}\left( 2M\right) }.$ Let $\lambda _{1}=\lambda _{1}\left(
b,\theta ,h_{0},N,M\right) >1$ and $C_{1}=C_{1}\left( b,\theta
,h_{0},N,M\right) >0$ be the numbers of Theorem 2. Denote $\widetilde{%
\lambda }=\lambda _{1}\left( b,\theta ,h_{0},N,3M\right) >1$ and $\widetilde{%
C}_{1}=C_{1}\left( b,\theta ,h_{0},N,3M\right) >0$. Then for any $\lambda
\geq \widetilde{\lambda }$ the functional $I_{h,\lambda }\left( W^{h}\right) 
$ is strictly convex on the ball $B_{0}\left( 2M\right) ,$ i.e. the
following analog of estimate (\ref{eq:convex}) holds for all $%
W^{h},W^{h}+r^{h}\in \overline{B_{0}\left( 2M\right) }:$ 
\begin{equation*}
I_{h,\lambda }\left( W^{h}+r^{h}\right) -I_{h,\lambda }\left( W^{h}\right)
-I_{h,\lambda }^{\prime }\left( W^{h}\right) \left( r^{h}\right) \geq 
\widetilde{C}_{1}e^{2\lambda \left( b-\theta \right) ^{2}}\left\Vert
r^{h}\right\Vert _{H_{2N}^{2,h}}^{2}.
\end{equation*}%
Furthermore, there exists a unique minimizer $W_{\min ,\lambda }^{h}$ of the
functional $I_{h,\lambda }\left( W^{h}\right) $ on the closed ball $%
\overline{B_{0}\left( 2M\right) }$ and the following inequality holds: 
\begin{equation*}
\left( I_{h,\lambda }^{\prime }\left( W_{\min ,\lambda }^{h}\right) ,W_{\min
,\lambda }^{h}-S\right) \leq 0\quad \text{for all }S\in \overline{%
B_{0}\left( 2M\right) }.
\end{equation*}%
Finally, let $Y_{\min ,\lambda }^{h}=W_{\min ,\lambda }^{h}+\Psi ^{h}.$ Then
the direct analog of (\ref{3.80}) holds where $V_{\min ,\lambda }^{h}$ is
replaced with $Y_{\min ,\lambda }^{h}$ and $V_{\ast }^{h}$ remains.
\end{thm}

We now construct the gradient projection method of the minimization of the
functional $I_{h,\lambda }\left( W^{h}\right) $ on the closed ball $%
\overline{B_{0}\left( 2M\right) }$. Let $\mathbf{P}:H_{2N,0}^{2,h}%
\rightarrow \overline{B_{0}\left( 2M\right) }$ be the orthogonal projection
operator of the space $H_{2N,0}^{2,h}$ on $\overline{B_{0}\left( 2M\right) }$%
. Let $W_{0}^{h}\in B_{0}\left( 2M\right) $ be an arbitrary point of this
ball. Let $\gamma \in \left( 0,1\right) $ be a number which we chose in
Theorem \ref{thm:9}. The sequence of the gradient projection method is:%
\begin{equation}
W_{n,\lambda ,\gamma }^{h}=\mathbf{P}\left( W_{n-1,\lambda ,\gamma
}^{h}-\gamma I_{h,\lambda }^{\prime }\left( W_{n-1,\lambda ,\gamma
}^{h}\right) \right) ,\text{ }n=1,2,...  \label{3.20}
\end{equation}%
By Theorem \ref{thm:8}, it holds that $I_{h,\lambda }^{\prime }\left(
W_{n-1,\lambda ,\gamma }^{h}\right) \in H_{2N,0}^{2,h}$. Since $%
W_{n-1,\lambda ,\gamma }^{h}\in \overline{B_{0}\left( 2M\right) }\subset
H_{2N,0}^{2,h}$, then all three terms in (\ref{3.20}) belong to $%
H_{2N,0}^{2,h}$.

\begin{thm}[Global convergence of the gradient projection method]
\label{thm:9} Assume that conditions of Theorem \ref{thm:8} hold and let $%
\lambda \geq \widetilde{\lambda }$. Then there exists a number $\gamma
_{0}=\gamma _{0}\left( b,\theta ,h_{0},N,3M\right) \in \left( 0,1\right) $
depending only on listed parameters such that for every $\gamma \in \left(
0,\gamma _{0}\right) $ there exists a number $\xi =\xi \left( \gamma \right)
\in \left( 0,\gamma _{0}\right) $ depending on $\gamma $ such that for these
values of $\gamma $ the sequence (\ref{3.20}) converges to $W_{\min ,\lambda
}^{h}$ and the following convergence rate holds: 
\begin{equation}
\left\Vert W_{n,\lambda ,\gamma }^{h}-W_{\min ,\lambda }^{h}\right\Vert
_{H_{2N}^{2,h}}\leq \xi ^{n}\left\Vert W_{\min ,\lambda
}^{h}-W_{0}^{h}\right\Vert _{H_{2N}^{2,h}}.  \label{3.21}
\end{equation}%
In addition, 
\begin{equation}
\left\Vert W_{\ast }^{h}-W_{n,\lambda ,\gamma }^{h}\right\Vert
_{H_{2N}^{2,h}}\leq \widetilde{C}_{2}\delta e^{4\lambda b\theta }+\xi
^{n}\left\Vert W_{\min ,\lambda }^{h}-W_{0}^{h}\right\Vert _{H_{2N}^{2,h}}.
\label{3.22}
\end{equation}
\end{thm}

\begin{rem}
Since the radius of the ball $\overline{B_{0}\left( 2M\right) }$ is $2M$,
where $M$ is an arbitrary number and since the starting point $W_{0}^{h}$ of
iterations of the sequence (\ref{3.20}) is an arbitrary point of $%
B_{0}\left( 2M\right) ,$ then Theorem \ref{thm:9} actually claims the \emph{%
global convergence }of the sequence (\ref{3.20}) to the exact solution, as
long as the level of the noise in the data $\delta $ tends to zero; see
Introduction for the definition of the global convergence.
\end{rem}

To close this section, let the function $c_{n,\lambda ,\gamma }^{h}\left( 
\mathbf{x}^{h}\right) $ be the one obtained after the substitution of the
components of the vector function $V_{n,\lambda ,\gamma }^{h}=W_{n,\lambda
,\gamma }^{h}+\Psi ^{h}$ in the real part of equation (\ref{eq:v}).
Similarly, let $c_{\ast }^{h}\left( \mathbf{x}^{h}\right) $ be obtained
after the substitution of the components of the vector function $V_{\ast
}^{h}=W_{\ast }^{h}+\Psi _{\ast }^{h}$. In the same vein, considering the
imaginary part of equation (\ref{eq:v}) we obtain the functions $\sigma
_{n,\lambda ,\gamma }^{h}\left( \mathbf{x}^{h}\right) $ and $\sigma _{\ast
}^{h}\left( \mathbf{x}^{h}\right) $, respectively. Then, the following
convergence estimates follow from Theorem \ref{thm:9}: 
\begin{equation*}
\left\Vert c_{\ast }^{h}-c_{n,\lambda ,\gamma }^{h}\right\Vert
_{L_{N}^{2,h}}+\left\Vert \sigma _{\ast }^{h}-\sigma _{n,\lambda ,\gamma
}^{h}\right\Vert _{L_{N}^{2,h}}\leq \widetilde{C}_{2}\delta e^{4\lambda
b\theta }+\xi ^{n}\left\Vert W_{\min ,\lambda }^{h}-W_{0}^{h}\right\Vert
_{H_{2N}^{2,h}}
\end{equation*}

\section{Experimental study\label{sec:Experimental-results}}

\subsection{Measured data}

Our raw data were experimentally collected at the microwave facility of The
University of North Carolina at Charlotte (UNCC). For brevity, we skip the
experimental setup and data acquisition because they were detailed in \cite%
{Khoa2020}. However, we recall that these are far field backscattering data
being collected for objects buried in a sandbox. As mentioned in Remark \ref%
{rem:11}, we need to apply the data propagation technique which approximates
the near field data for our CIP. This technique was derived in \cite%
{Nguyen2018} and was revisited in \cite{Khoa2020}.

We introduce dimensionless spatial variables as $\mathbf{x}^{\prime }=%
\mathbf{x}/(10\text{ cm})$. This means that the dimensions we use below are
10 times less than the real ones in centimeters. Cf. Figure \ref%
{fig:Experimental-setup}, we briefly describe our experiment with relatively
small targets. Typically the sizes of antipersonnel land mines and
improvised explosive devices are between 0.5 and 1.5 (i.e., between 5 cm and
15 cm), see, e.g. \cite{Nguyen2018}. We have used a wooden-like framed box
filled with the dry sand and covered by bending styrofoam layers from the
front and back.The experimental objects were buried in that sand. The
dielectric constant and the electrical conductivity of that sand, i.e. the
these parameters of the background were $c_{\text{bckgr}}=4$ and $\sigma _{%
\text{bckgr}}=0.$ Note that in the styrofoam $c=1$ and $\sigma =0.$%
Therefore, the styrofoam should not affect neither the incident nor the
scattered electric waves. As to the line of sources $L_{\text{src}}$ defined
in (\ref{Lsrc}), we have $d=9$, $a_{1}=0.1$ and $a_{2}=0.6$ with 0.1
mesh-width.

For each source position, we have collected two sets of complex valued data $%
u$. The first set is the reference data, i.e. the data for the case when
only the sand was present in that box, and a target of interest was not
present. The second set of the data was for the case when that unknown
target, which we want to image, was buried in that sandbox. We do not know
yet whether or not we can work with the scenario when the reference data are
not measured. We point out, however, that in previous works of this group on
experimental data for buried targets, the data for the reference medium were
also collected \cite{Klibanov2019b,Nguyen2018,Thanh2015}.

For each position of the source, the raw data consist of multi-frequency
backscatter data associated with 300 frequency points uniformly distributed
between 1 GHz to 10 GHz. However, for each selected target, only a single
frequency for each target was used in our computations, i.e. our
data are non-overdetermined ones with $m=n=3;$ see Introduction of
the definition of non overdetermined data. The backscattering data were
measured on a square, which was a part of a plane. The dimensions of this
square were $10\times 10$ in dimensionless units. These measurements were
done by a detector which was moving in both vertical and horizontal
direction with the moving step size 0.2. Thus, total $50\times 50=2500$
positions of the detector on that measurement plane were used for each
position of the source. It is worth noting that not all propagated data have
a good quality. Therefore, we need to choose proper data among those
frequency-dependent data sets for each target of interest. We tabulate in
Table \ref{table:1} the chosen wavenumber $k$ as well as the corresponding
frequency $\tilde{f}$ for each test, using the standard dimensionless
formulation $k=2\pi \tilde{f}/2997924580$. We now briefly summarize the
crucial steps of the data preprocessing to obtain fine near-field data for
our CIP from the raw backscattering ones. Comparison of the raw and
propagated data in Figures \ref{fig:RawProp1}--\ref{fig:imgProp1}, \ref%
{fig:RawProp2}--\ref{fig:imgProp2}, \ref{fig:RawProp3}--\ref{fig:imgProp3}, %
\ref{fig:RawProp4}--\ref{fig:imgProp4}, \ref{fig:RawProp5}--\ref%
{fig:imgProp5} justifies this data preprocessing procedure.

\begin{itemize}
\item \textbf{Step 1.} Subtract the reference data from the far-field
measured data for every frequency and for each source position. The
reference data mean the ones measured when the sandbox is without a target.
This subtraction helps to extract the pure signals, which always contain
unwanted noises, from buried objects from the whole signal. In our
computations, we treat the resulting data as the ones for the background
values as in vacuum, i.e. $\widetilde{c}_{\text{bckgr}}=1$ and $\sigma _{%
\text{bckgr}}=0.$ The latter is our heuristic assumption, which,
nevertheless, works quite well for our reconstructions.

\item \textbf{Step 2.} For each position of the source, apply the data
propagation procedure to approximate the near field data. We propagated the
data up to the the sand surface, i.e. after getting through the styrofoam
layer. As a result, good estimates of $x,y$ coordinates of buried objects
were obtained, see Figures \ref{fig:RawProp1}--\ref{fig:imgProp1}, \ref%
{fig:RawProp2}--\ref{fig:imgProp2}, \ref{fig:RawProp3}--\ref{fig:imgProp3}, %
\ref{fig:RawProp4}--\ref{fig:imgProp4}, \ref{fig:RawProp5}--\ref%
{fig:imgProp5}. Besides, this procedure reduces the size of the
computational domain in the $z$--direction.

\item \textbf{Step 3.} Truncate the obtained near field data to get rid of
random oscillations that appear randomly during the data propagation.
Suppose the function $\mathcal{K}\left( x,y,\alpha \right) $ represents
those near field data. We truncate this function in two steps A and B:

\textbf{A.} Replace $\mathcal{K}\left( x,y,\alpha \right) $ with the
function $\tilde{\mathcal{K}}\left( x,y,\alpha \right) $ defined as: 
\begin{equation}
\tilde{\mathcal{K}}\left( x,y,\alpha \right) =%
\begin{cases}
\mathcal{K}\left( x,y,\alpha \right) & \text{if }\left\vert \mathcal{K}%
\left( x,y,\alpha \right) \right\vert \geq \kappa _{1}\max_{\left\vert
x\right\vert ,\left\vert y\right\vert \leq R}\left\vert \mathcal{K}\left(
x,y,\alpha \right) \right\vert , \\ 
0 & \text{otherwise.}%
\end{cases}
\label{eq:trun}
\end{equation}%
Here, $\kappa _{1}>0$ represents the truncation number and we take $\kappa
_{1}=0.4$ based upon the trial and error procedure. This means that we only
preserve those propagated near field data whose absolute values are at least
40\% of their global maximum value.

\textbf{B.} Smooth the function $\tilde{\mathcal{K}}$ using the Gaussian
filter. However, we have observed that the maximal absolute value of the
smoothed function $\tilde{\mathcal{K}}_{\text{sm}}$ is less than that value
of the original function $\tilde{\mathcal{K}}.$ This, however, results in
lesser values of the computed dielectric constants of targets we image.
Thus, we preserve those maximal absolute values. More precisely, we replace
the function $\tilde{\mathcal{K}}_{\text{sm}}$ with the function $\tilde{%
\mathcal{K}}_{\text{new}}\left( x,y,\alpha \right) =\kappa _{2}\tilde{%
\mathcal{K}}_{\text{ sm}}\left( x,y,\alpha \right) $. We call $\kappa _{2}>0$
the \textquotedblleft retrieval number". This number is defined as $\kappa
_{2}=\max \left( \left\vert \tilde{\mathcal{K}}\right\vert \right) /\tilde{m}
$, where $\tilde{m}$ is the maximal absolute value of $\tilde{\mathcal{K}}_{%
\text{sm}}$.
\end{itemize}

The distance between the measurement plane and the sandbox with the
styrofoam layer is about 11.05. Meanwhile, the length in the $z$ direction
of the sandbox without the styrofoam is 4.4. Since the thickness of the
bending front styrofoam layer is 0.5, then the domain of interests $\Omega $
should be $\Omega =\left\{ \mathbf{x}\in \mathbb{R}^{3}:\left\vert
x\right\vert ,\left\vert y\right\vert <5,\left\vert z\right\vert <2\right\}
, $ which implies that $R=5$ and $b=2$. The near-field measurement site is
then $\Gamma =\left\{ \mathbf{x}\in \mathbb{R}^{3}:\left\vert x\right\vert
,\left\vert y\right\vert <5,z=-2\right\} .$ Also, we report that the
distance between the far field measurement site and the zero point of our
coordinate system was 14. Given that by Table \ref{table:1} our
average frequency was 4.28 GHz, our average wave length was 7 cm (cf. \cite{GHz}),
which is 0.7 in our dimensionless units. Therefore, the distance between the
measurement site and the zero point of our coordinate system was about 20
wave lengths, which is a far field zone.

\subsection{Reconstruction results}

Five (5) examples for our reconstructions of buried objects are presented
here. They are basically in-store items that one can purchase easily.
Nevertheless, they really mimic some well-known metallic and non-metallic
antipersonnel land mines met during the World War eras. More precisely, this is true for images presented on Figure \ref{fig:tube-1} that follows the NO-MZ 2B mine and on Figure \ref%
{fig:bottle-1} that imitate the Glassmine 43 of the Germans; see \cite[Section 4.2]{Khoa2020}. For brevity, photos and descriptions of these five
experimental objects are presented in Figures \ref{fig:tube-1}, \ref%
{fig:bottle-1}, \ref{fig:Walnut}, \ref{fig:A}, \ref{fig:O} and Table \ref%
{table:2}, where the reader can find their distinctive levels of geometry
and materials.

Cf. those cited in \cite{Khoa2020}, the sizes of such military antipersonnel
land mines are only between 0.5 and 1.5. We focus on finding them in a
sub-domain of $\Omega $ with only 2 in depth in the $z-$direction. Denote
this sub-domain by $\Omega _{1}=\left\{ -b\leq z\leq -b+2\right\} =\left\{
-2\leq z\leq 0\right\} $. This results in the following choice of our
starting point of iterations in the minimization of the functional $%
J_{\lambda ,h}\left( V^{h}\right) $. Denote that starting point by the
vector $V_{0}^{h}=V_{0}\left( x_{p},y_{q},z_{s}\right) $ with 
\begin{equation}
V_{0}^{h}=%
\begin{pmatrix}
v_{00}^{h} & v_{01}^{h} & \cdots & v_{0\left( N-1\right) }^{h}%
\end{pmatrix}%
^{T},\quad v_{0n}^{h}=\left( \psi _{0n}^{h}+\psi _{1n}^{h}\left( z+b\right)
\right) \chi \left( z\right) .  \label{30}
\end{equation}%
Recall that $\psi _{0n}^{h}$ and $\psi _{1n}^{h}$ are the Fourier
coefficients of the propagated data in (\ref{3.3}). Here, $\chi :\left[ -b,b%
\right] \rightarrow \mathbb{R}$ is the smooth function given by 
\begin{equation*}
\chi \left( z\right) =%
\begin{cases}
\exp \left( \frac{2\left( z+b\right) ^{2}}{\left( z+b\right) ^{2}-b^{2}}%
\right) & \text{if }z<0, \\ 
0 & \text{otherwise}.%
\end{cases}%
\end{equation*}%
This function attains the maximal value 1 at $z=-b$ exactly where the
near-field data are given. Then, it holds that $v_{0n}^{h}\mid _{z=-b}=\psi
_{0n}^{h}$, $\partial _{z}v_{0n}^{h}\mid _{z=-b}=\psi _{1n}^{h}$. Moreover, $%
\chi $ tends to 0 as $z\rightarrow 0^{+}$ leading to $v_{0n}^{h}\mid
_{z=b}=\partial _{z}v_{0n}^{h}\mid _{z=b}=0.$ Hence, our starting point (\ref%
{30}) satisfies the boundary conditions (\ref{3.3}).

\begin{figure*}[t]
\vspace{-7mm} 
\begin{centering}
\includegraphics[scale=0.5]{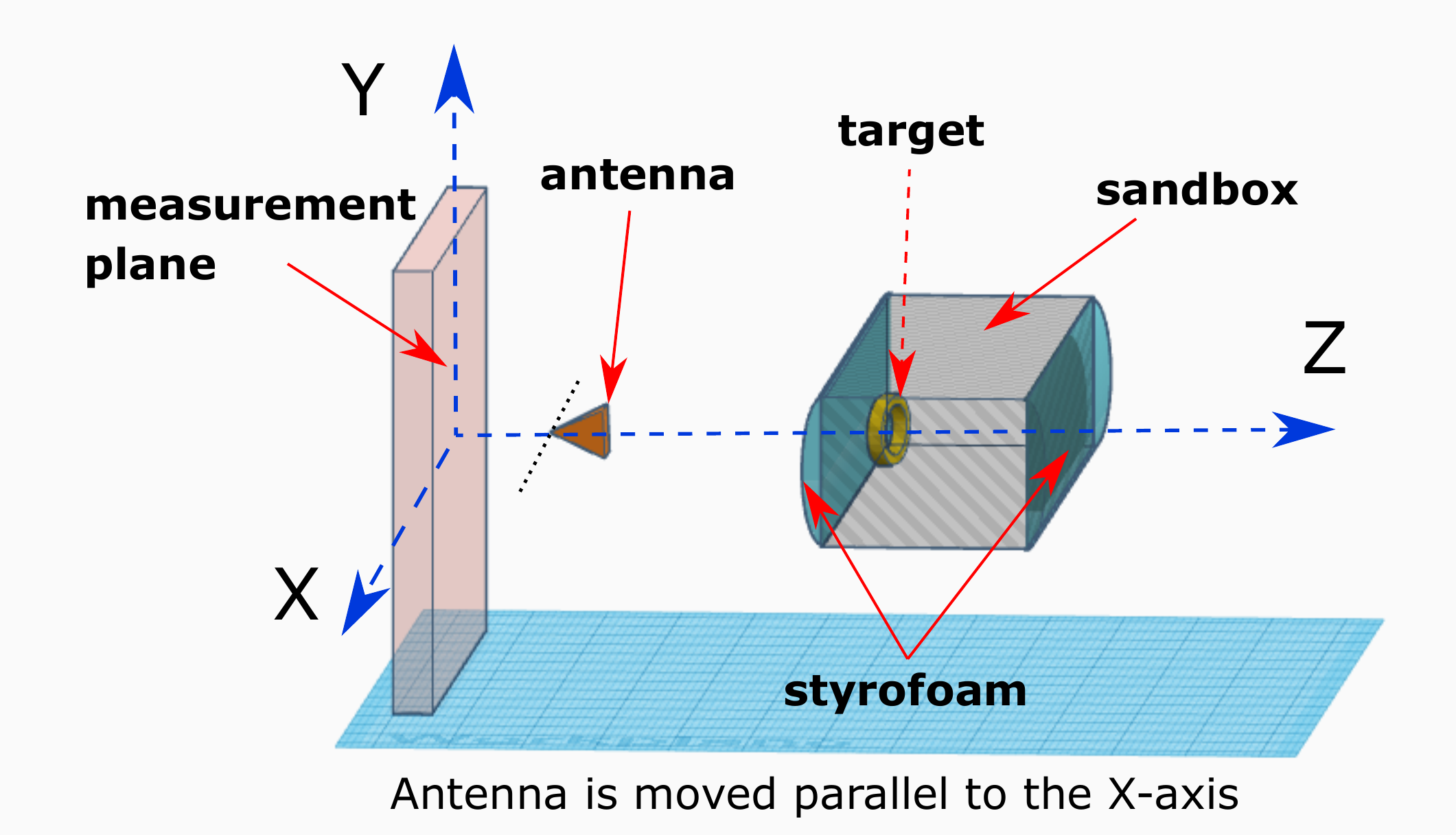}
\par\end{centering}
\caption{A schematic diagram illustrating our experimental configuration. The
transmitter is the antenna put in front of the detector placing the
far-field measurement site. The backscattering waves hit the antenna
before reaching the detectors. This causes a certain noise in the data.}
\label{fig:Experimental-setup}
\end{figure*}

\begin{table}
	\begin{center}
		\begin{tabular}{|c|c|c|c|c|c|}
			\hline Example & 1 & 2 & 3 & 4 & 5 \tabularnewline
			\hline 
			$k$ & 8.51 & 6.62 & 11.43 & 9.55 & 8.79 \tabularnewline
			\hline 
			Frequency (GHz) & 4.06 & 3.16 & 5.45 & 4.55 & 4.19\tabularnewline
			\hline 
			
		\end{tabular}
	\end{center}
	\caption{Chosen wavenumbers and frequencies for Examples 1--5.\label{table:1}}
\end{table}

\begin{table}
	\begin{centering}
		\begin{tabular}{|c|c|c|c|c|c|}
			\hline 
			Example & 1 & 2 & 3 & 4 & 5\tabularnewline
			\hline 
			Object & Metallic cylinder & Water & Wood & Metallic letter ``A'' & Metallic letter ``O''\tabularnewline
			\hline 
			$\max\left(c_{\text{comp}}\right)$ & 28.78 & 23.18 & 6.33 & 16.24 & 16.22\tabularnewline
			\hline 
			$c_{\text{true}}$ & $\left[10,30\right]$ & 23.8 & $\left[2,6\right]$ & $\left[10,30\right]$ & $\left[10,30\right]$\tabularnewline
			\hline 
			$\max\left(\sigma_{\text{comp}}\right)$ & 2.33 & 0.94 & 0.94 & 1.32 & 1.66\tabularnewline
			\hline 
		\end{tabular}
		\par\end{centering}
	\caption{True $c_{\text{true}}$, computed $\max \left( c_{\text{%
				comp}}\right) $ dielectric constants and computed conductivity $\max \left( \sigma_{\text{%
				comp}}\right) $ of Examples 1--5 of experimental data.
		True values of dielectric constants were taken from: (a) Examples 1, 4, 5: formula (7.2) of \cite{Kuzh}%
		, (b)  Example 2 (clear water) \cite{Thanh2015}, (c) Example 3 \cite{Table}. \label{table:2}}
\end{table}

Even though our analysis in section \ref{sec:A-globally-convergent} is
applicable to the semi-discrete form of the functional $J_{h,\lambda }\left(
V^{h}\right) $ defined in (\ref{eq:J}), to perform computations, we
naturally write it in the fully discrete form. In this form, we take the
uniform mesh in $x,y,z$ directions, $\left\{ \left( x_{p},y_{q},z_{s}\right)
\right\} _{p,q,s=0}^{Z_{h}}$, where the mesh sizes in $x,y,z$ directions are
the same. For brevity, we do not bring in here this fully discrete form of $%
J_{h,\lambda }\left( V^{h}\right) .$ After obtaining the global minimum $%
V_{p,q,s}$ of the functional $J_{h,\lambda }\left( V^{h}\right) $, we
compute the unknown coefficients $c_{p,q,s}$ and $\sigma _{p,q,s}$ as
follows: 
\begin{align*}
& c_{p,q,s}=\text{mean}_{\alpha }\left\vert \func{Re}\left\{ -\frac{\Delta
^{h}v_{p,q,s,\alpha _{l}}+\left( \nabla ^{h}v_{p,q,s,\alpha _{l}}\right)
^{2}+2\nabla ^{h}v_{p,q,s,\alpha _{l}}\cdot \tilde{\mathbf{x}}_{p,q,s,\alpha
_{l}}}{k^{2}}\right\} \right\vert +1, \\
& \sigma _{p,q,s}=\text{mean}_{\alpha }\left\vert \func{Im}\left\{ -\frac{%
\Delta ^{h}v_{p,q,s,\alpha _{l}}+\left( \nabla ^{h}v_{p,q,s,\alpha
_{l}}\right) ^{2}+2\nabla ^{h}v_{p,q,s,\alpha _{l}}\cdot \tilde{\mathbf{x}}%
_{p,q,s,\alpha _{l}}}{0.1k\eta _{0}}\right\} \right\vert ,
\end{align*}%
aided by \eqref{eq:v}; see also Remark \ref{rem:X}. Here mean$_{\alpha }$
denotes the average value with respect to the positions of the source $%
\alpha .$ The number 0.1 presented in the computed conductivity is due to
its physical unit $\text{S/m}$ in this dimensionless regime. Since the
number of point sources is very limited, we use the Gauss--Legendre
quadrature method to compute the measured near-field data in the Fourier
series.

Here, we use the gradient descent method for the minimization of the target
functional $J_{h,\lambda }\left( V^{h}\right) $ of (\ref{eq:J}) due to its
easy implementation. Even though Theorem \ref{thm:9} claims the global
convergence of the gradient projection method, our success in working with
the gradient descent method is similar with those in all previous
publications that study the convexification \cite%
{Khoa2019,convIPnew,Klibanov2019,Klibanov2019b,Klibanov2019a,Klibhyp}. As to
the value of the parameter $\lambda $ in $J_{h,\lambda }\left( V^{h}\right)
, $ even though the above theorems require large values of $\lambda ,$ our
numerical experience tells us that we can choose a moderate value $\lambda
=1.1$, which was in the range $\lambda \in \left[ 1,3\right] $ chosen in all
above cited publications on the convexification.

Concerning the step size $\gamma $ of the gradient descent method, we start
from $\gamma _{1}=10^{-1}$. On each step of iterations $m\geq 1$, the next
step size $\gamma _{m+1}=\gamma _{m}/2,$ if the value of the functional on
the step $m$ exceeds its value of the previous step. Otherwise, we set $%
\gamma _{m+1}=\gamma _{m}$. The minimization process is stopped when either $%
\gamma _{m}<10^{-10}$ or $\left\vert J_{h,\lambda }\left( V_{m}^{h}\right)
-J_{h,\lambda }\left( V_{m-1}^{h}\right) \right\vert <10^{-10}$. As to the
gradient $J_{h,\lambda }^{\prime }$ of the discrete functional $J_{h,\lambda
}$, we apply the technique of Kronecker deltas (cf. e.g. \cite{Kuzhuget2010}%
) to derive its explicit formula. For brevity, we do not provide this
formula here.

After the minimization procedure is stopped, we obtain the discrete
coefficient of $c_{p,q,s}$. We apply to $c_{p,q,s}$ the truncation procedure
described in the above \textbf{Step 3} and in (\ref{eq:trun}). But now we
use $\kappa _{1}=0.2$ instead of $\kappa _{1}=0.4$ which was used for the
data preprocessing. Denote the resulting function by $\tilde{c}$. Our
reconstructed solution, denoted by $c_{\text{comp}}$, is obtained after we
smooth $\tilde{c}$ by the standard filtering via the \texttt{smooth3}
built-in function in MATLAB. We find $c_{\text{comp}}$ by using $c_{\text{%
comp}}=\rho \text{smooth}\left( \left\vert \tilde{c}\right\vert \right) $
where the number $\rho >0$ is found the same way as the number $\kappa _{2}$
in the above \textbf{Step 3}. As to the coefficient $\sigma _{\text{comp}}$,
we follow the same vein with the same truncation parameter.


\begin{figure}[H]\vspace{-7mm}
	\begin{centering}\hspace*{\fill}
		\subfloat[Real part of raw and propagated data at $\alpha=0.4$\label{fig:RawProp1}]{\begin{centering}
				\includegraphics[scale=0.3]{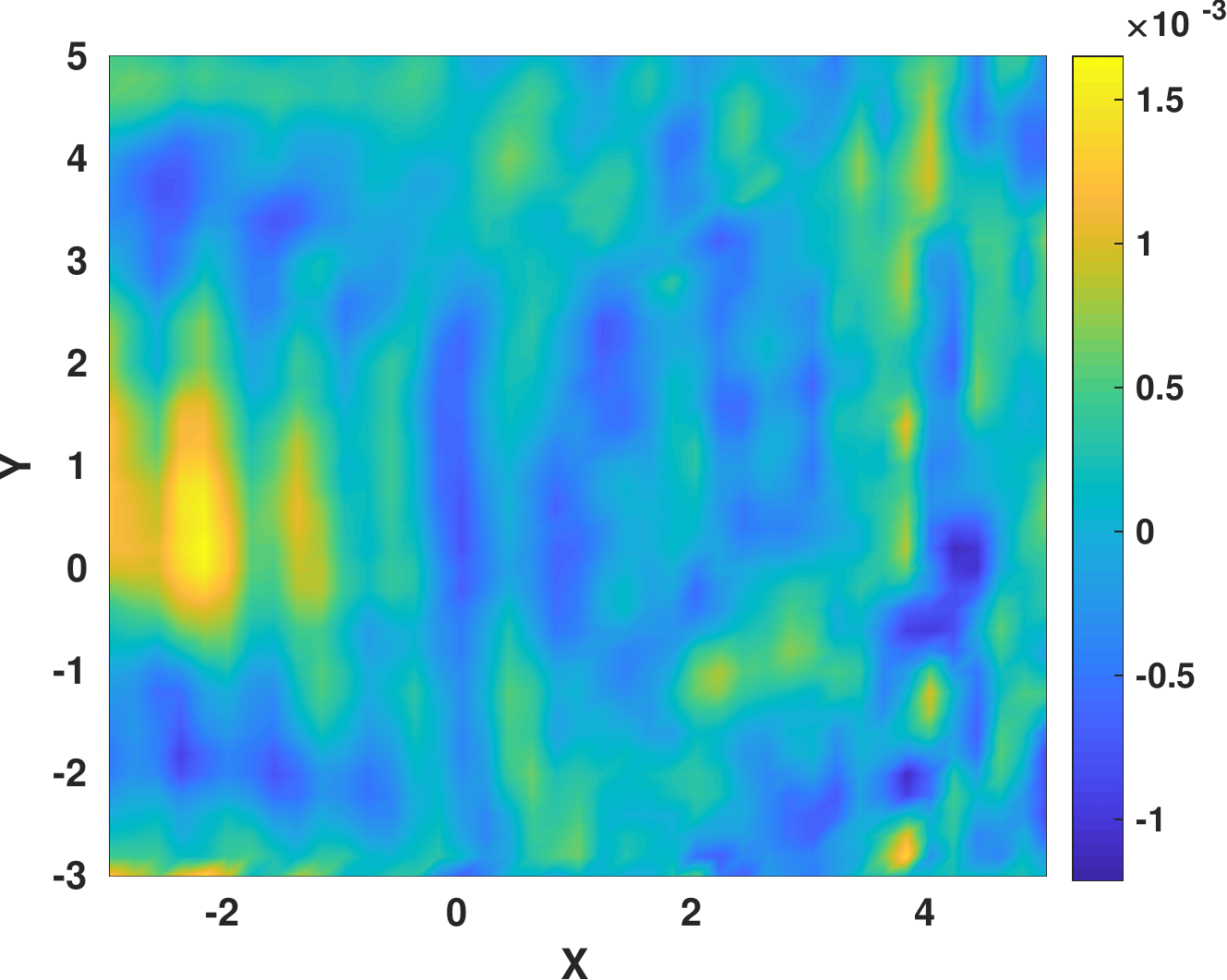}
				\includegraphics[scale=0.3]{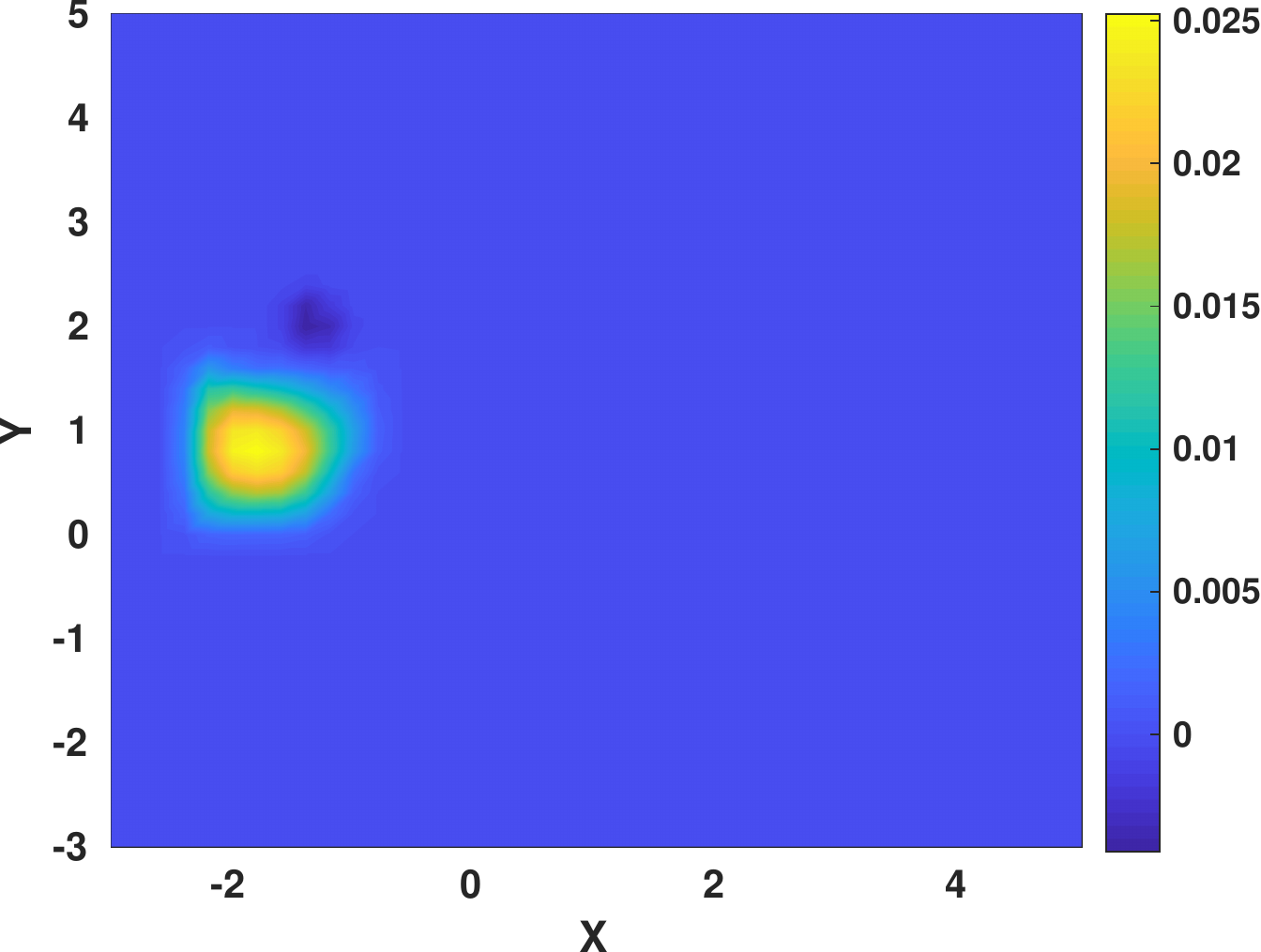}
				\par\end{centering}
		}\subfloat[Imaginary part of raw and propagated data at $\alpha=0.4$\label{fig:imgProp1}]{\begin{centering}
				\includegraphics[scale=0.3]{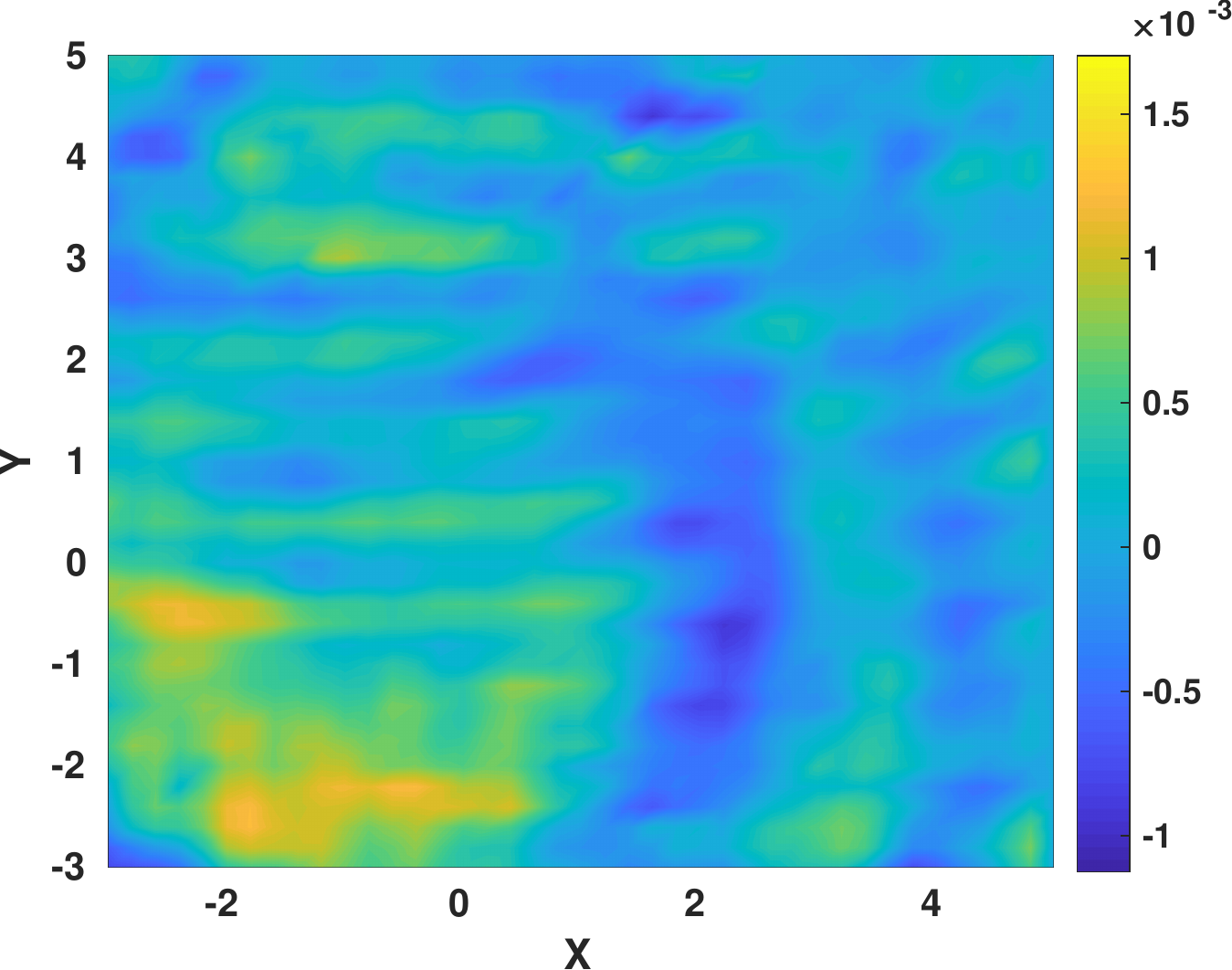}
				\includegraphics[scale=0.3]{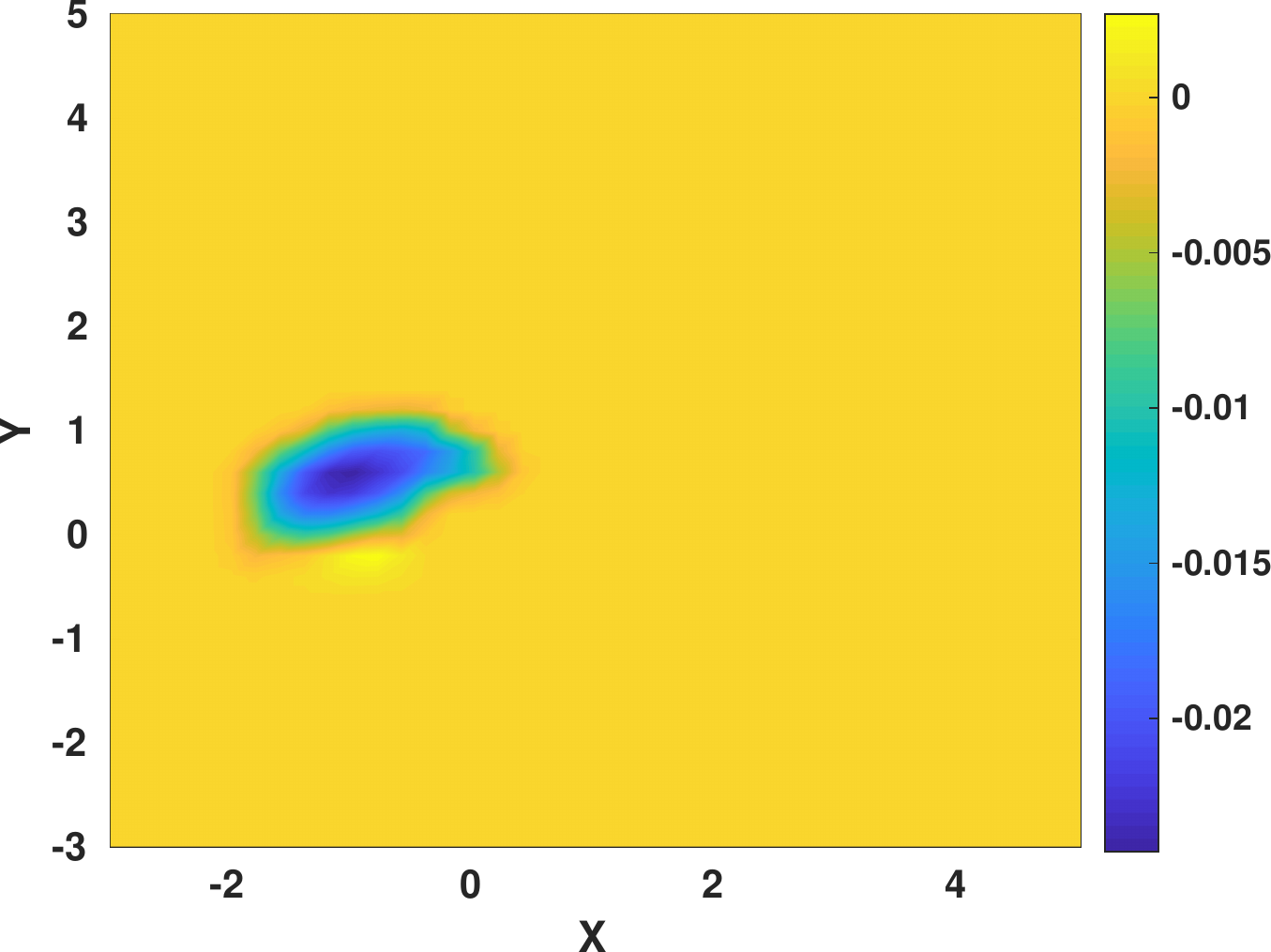}
				\par\end{centering}
		}
		\par\end{centering}
	\begin{centering}\vspace*{\fill}
		\subfloat[Left: Aluminum cylinder (cf. \cite{Khoa2020}). Middle: Image of computed dielectric constant. Right: Image of computed conductivity\label{fig:tube-1}]{ \begin{centering}
				\includegraphics[scale=0.2]{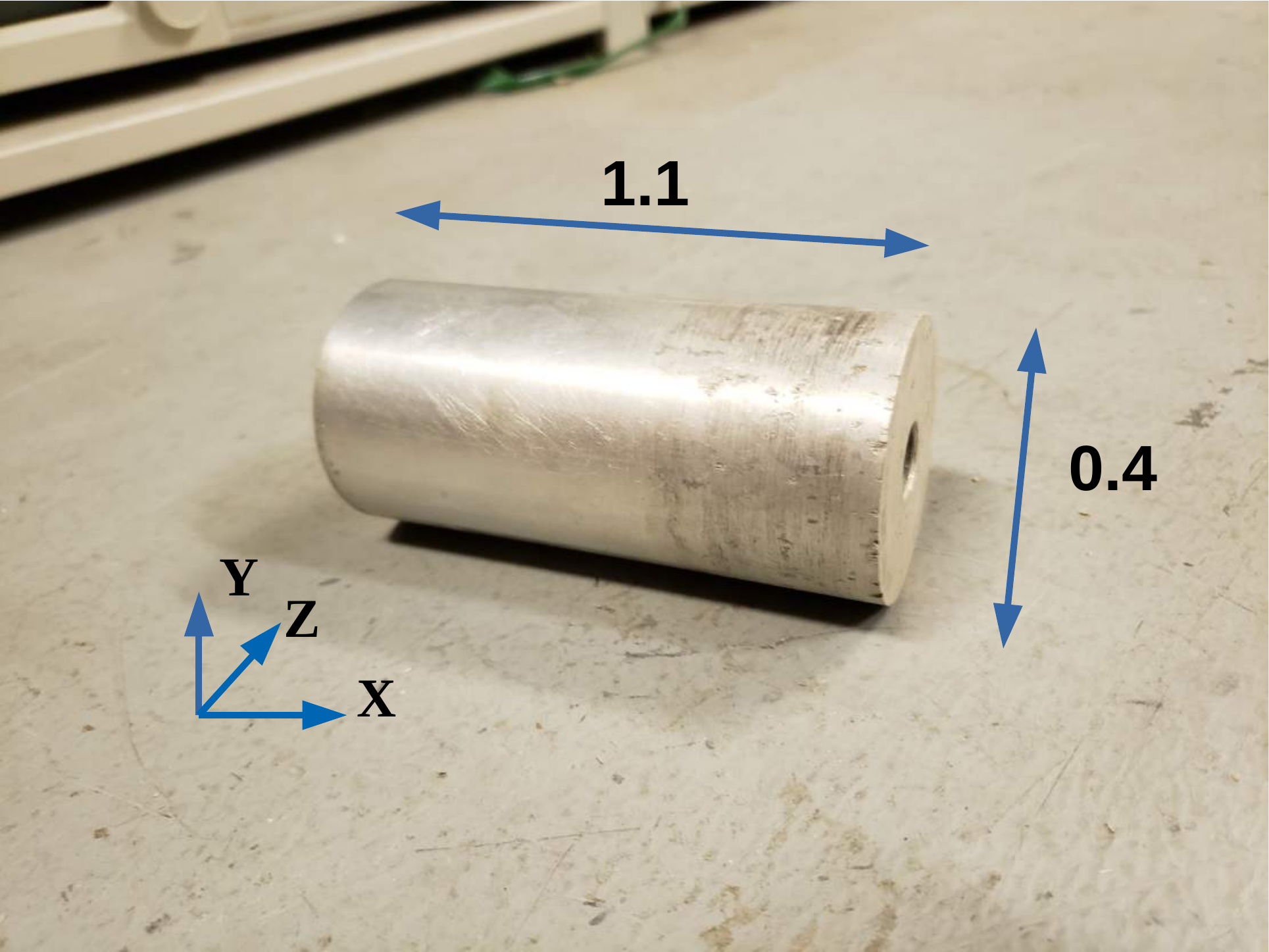}
				\includegraphics[scale=0.4]{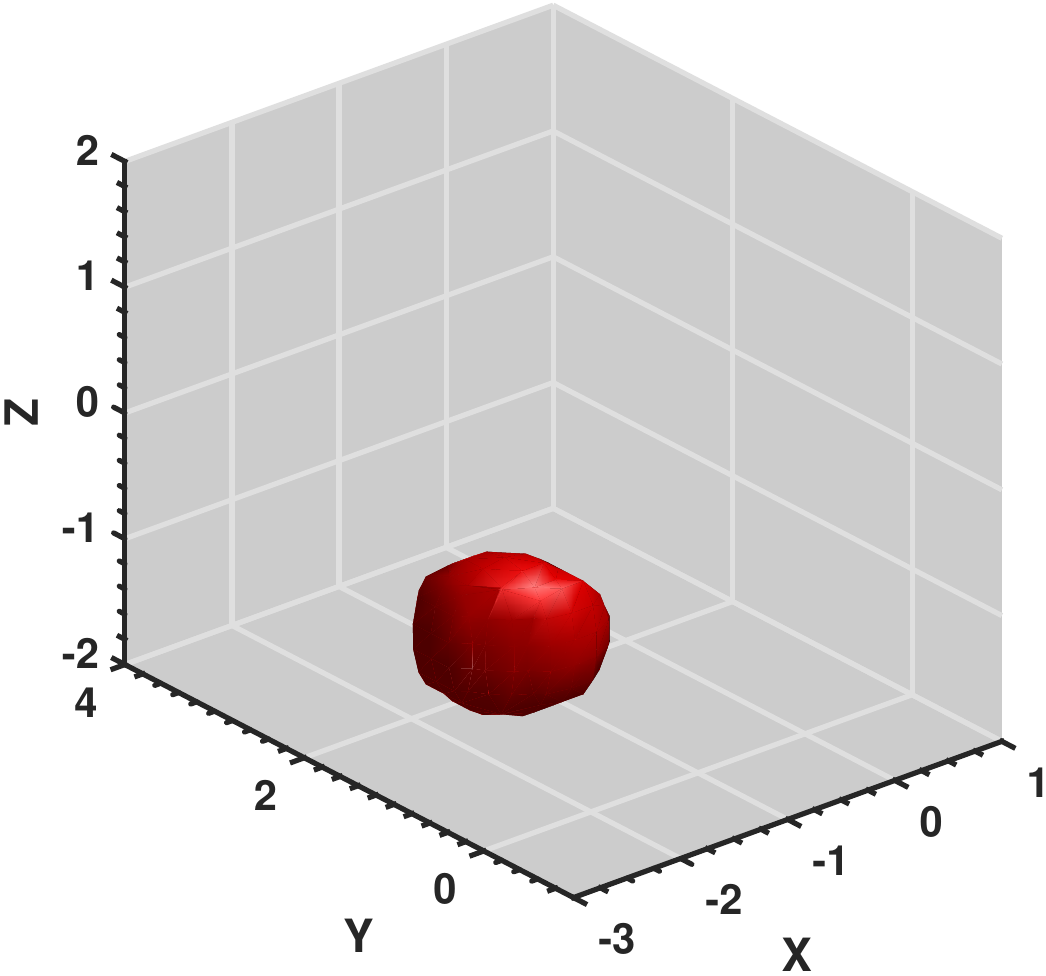}
				\includegraphics[scale=0.4]{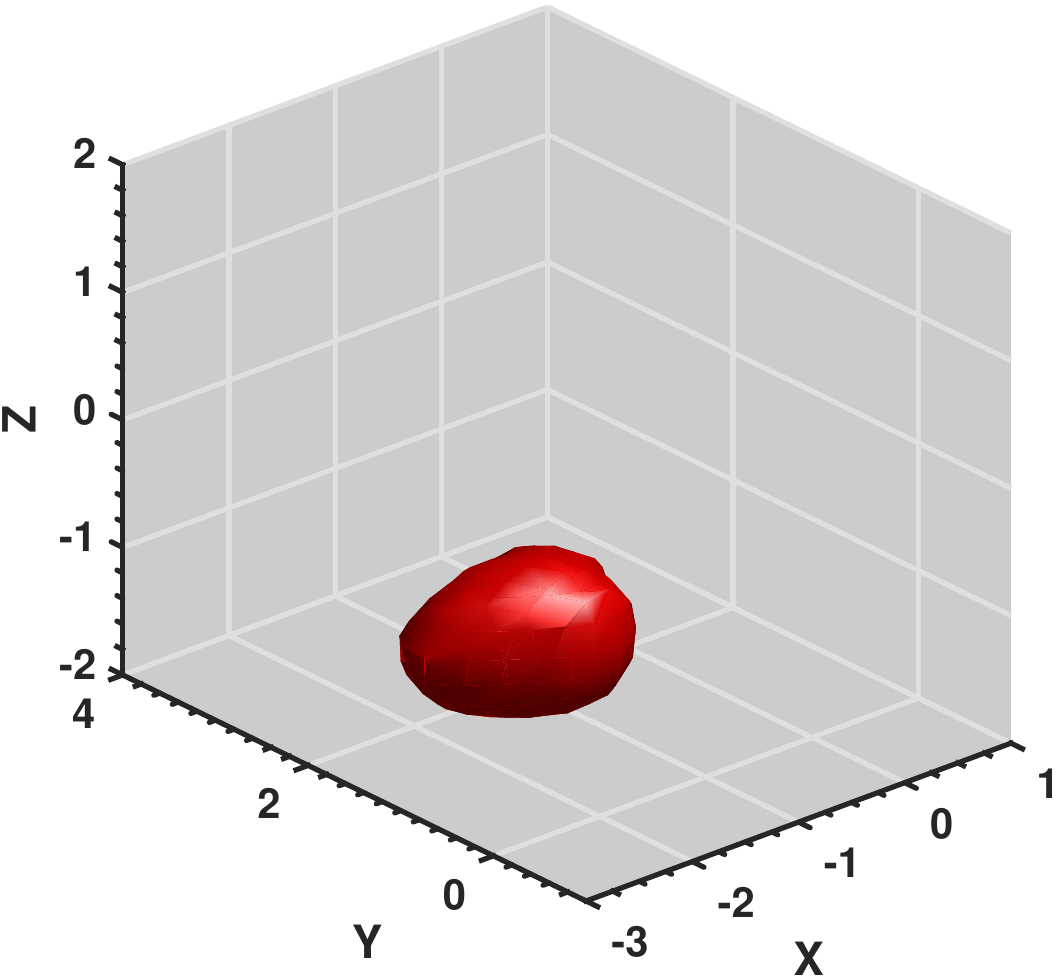}
				\par\end{centering}
		}
		\par\end{centering}\hspace*{\fill}
	\caption{Aluminum cylinder; see Tables \ref{table:1}--\ref{table:2} for further details. (a) and (b) demonstrate a clear advantage of the data propagation procedure in data preprocessing.}\vspace{-5mm}
\end{figure}

\begin{figure}[H]\vspace{-7mm}
	\begin{centering}\hspace*{\fill}
		\subfloat[Real part of raw and propagated data at $\alpha=0.4$\label{fig:RawProp2}]{\begin{centering}
				\includegraphics[scale=0.3]{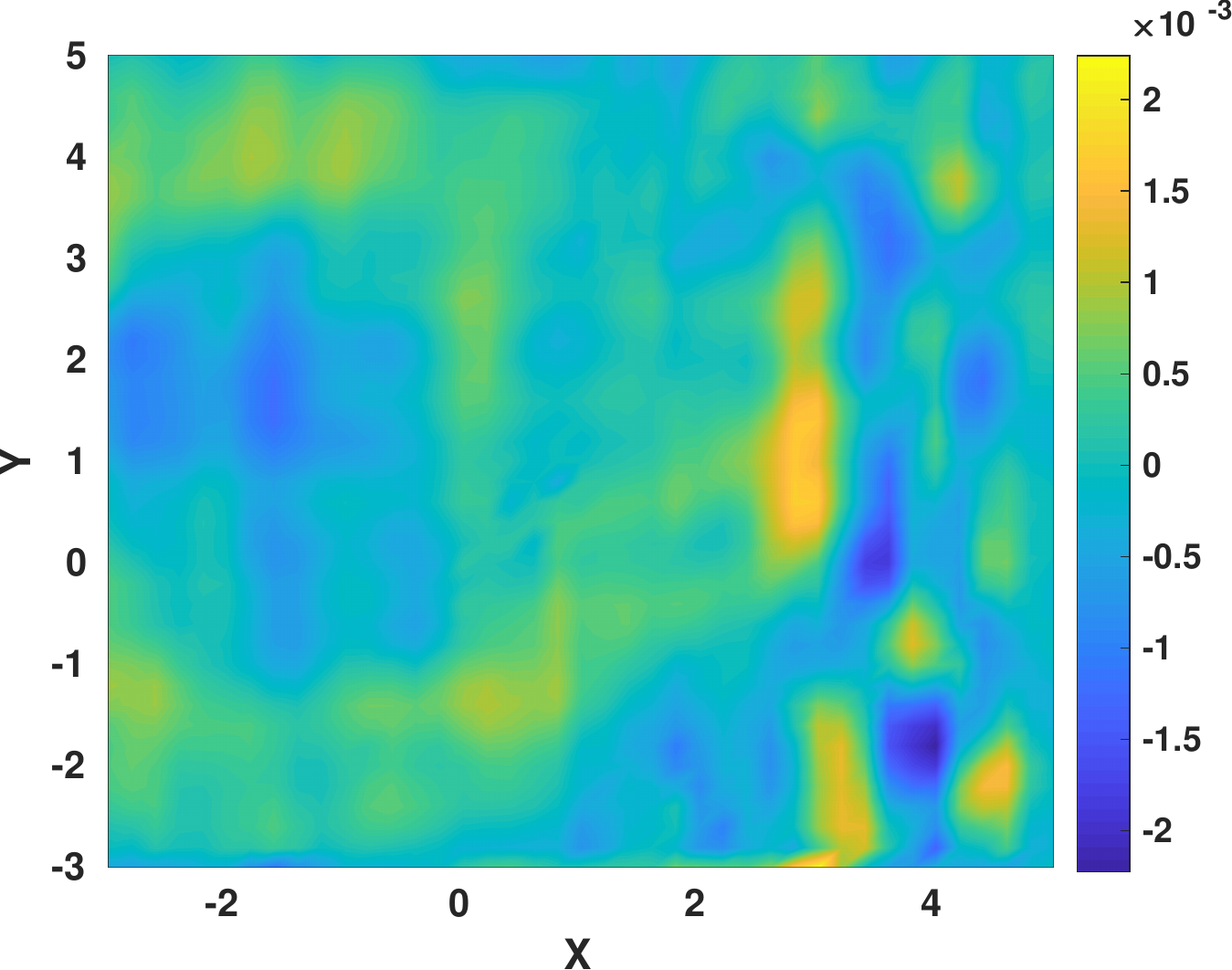}
				\includegraphics[scale=0.3]{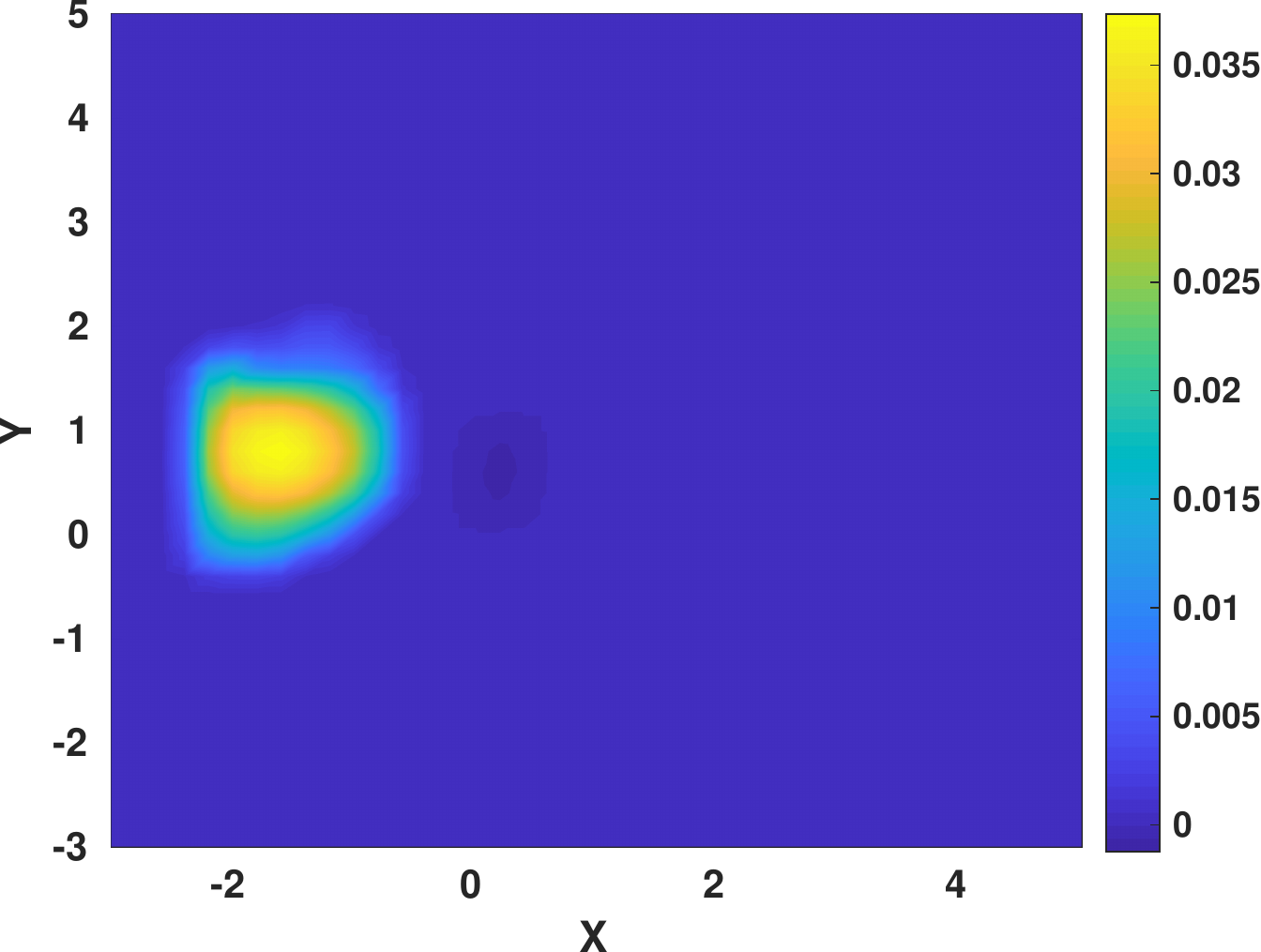}
				\par\end{centering}
		}\subfloat[Imaginary part of raw and propagated data at $\alpha=0.4$\label{fig:imgProp2}]{\begin{centering}
				\includegraphics[scale=0.3]{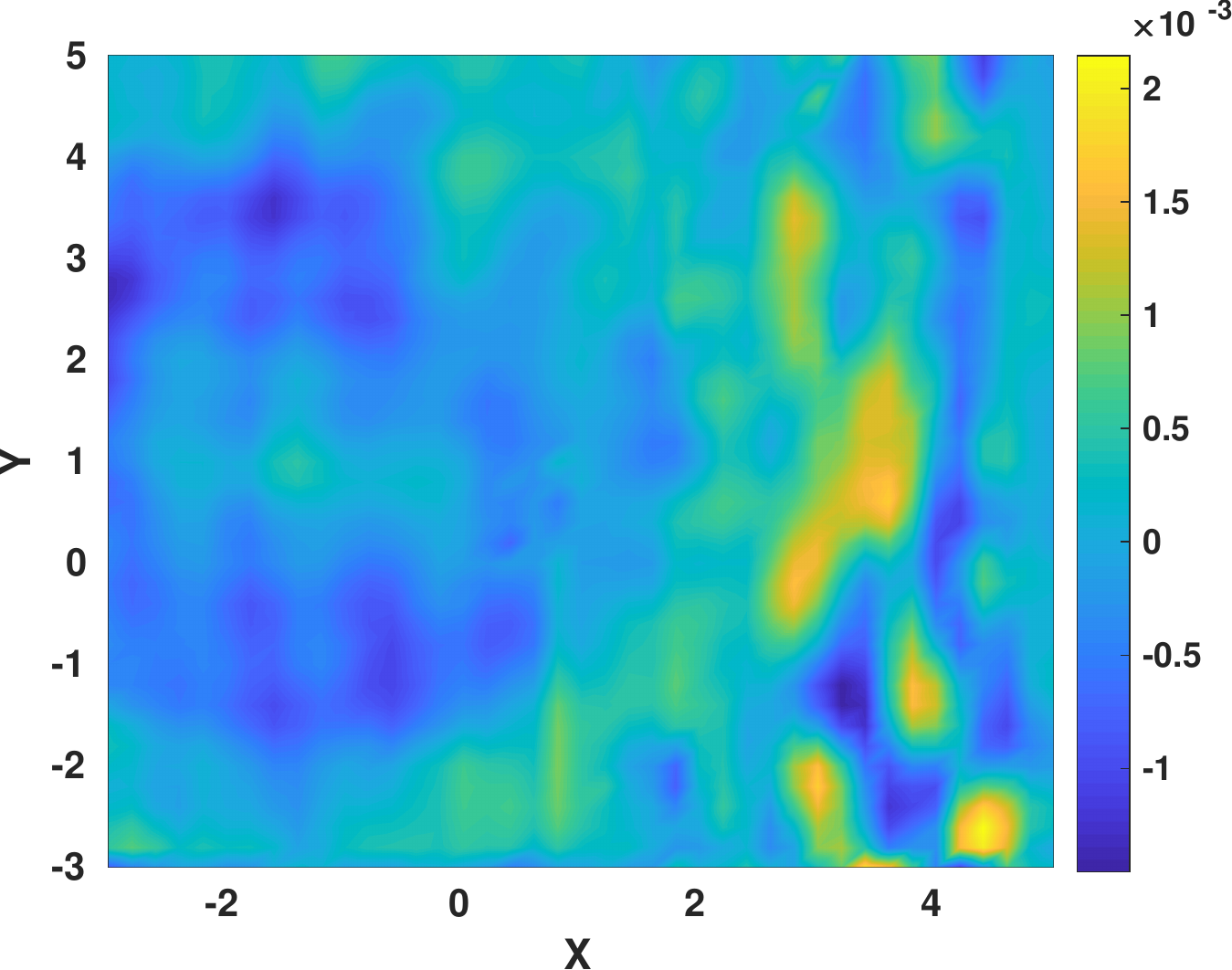}
				\includegraphics[scale=0.3]{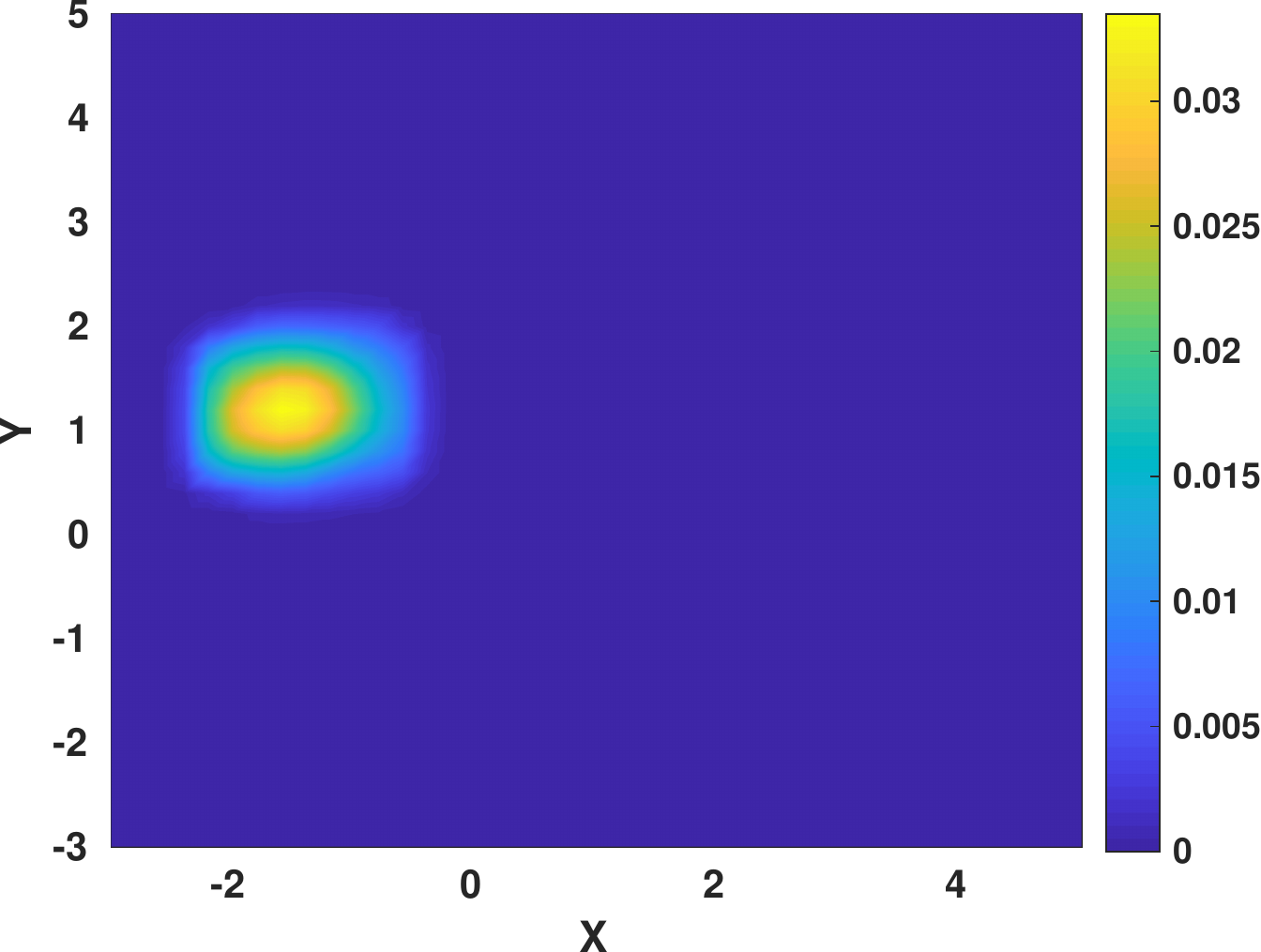}
				\par\end{centering}
		}
		\par\end{centering}
	\begin{centering}\vspace*{\fill}
		\subfloat[Left: Glass bottle (cf. \cite{Khoa2020}). Middle: Image of computed dielectric constant. Right: Image of computed conductivity\label{fig:bottle-1}]{ \begin{centering}
				\includegraphics[scale=0.2]{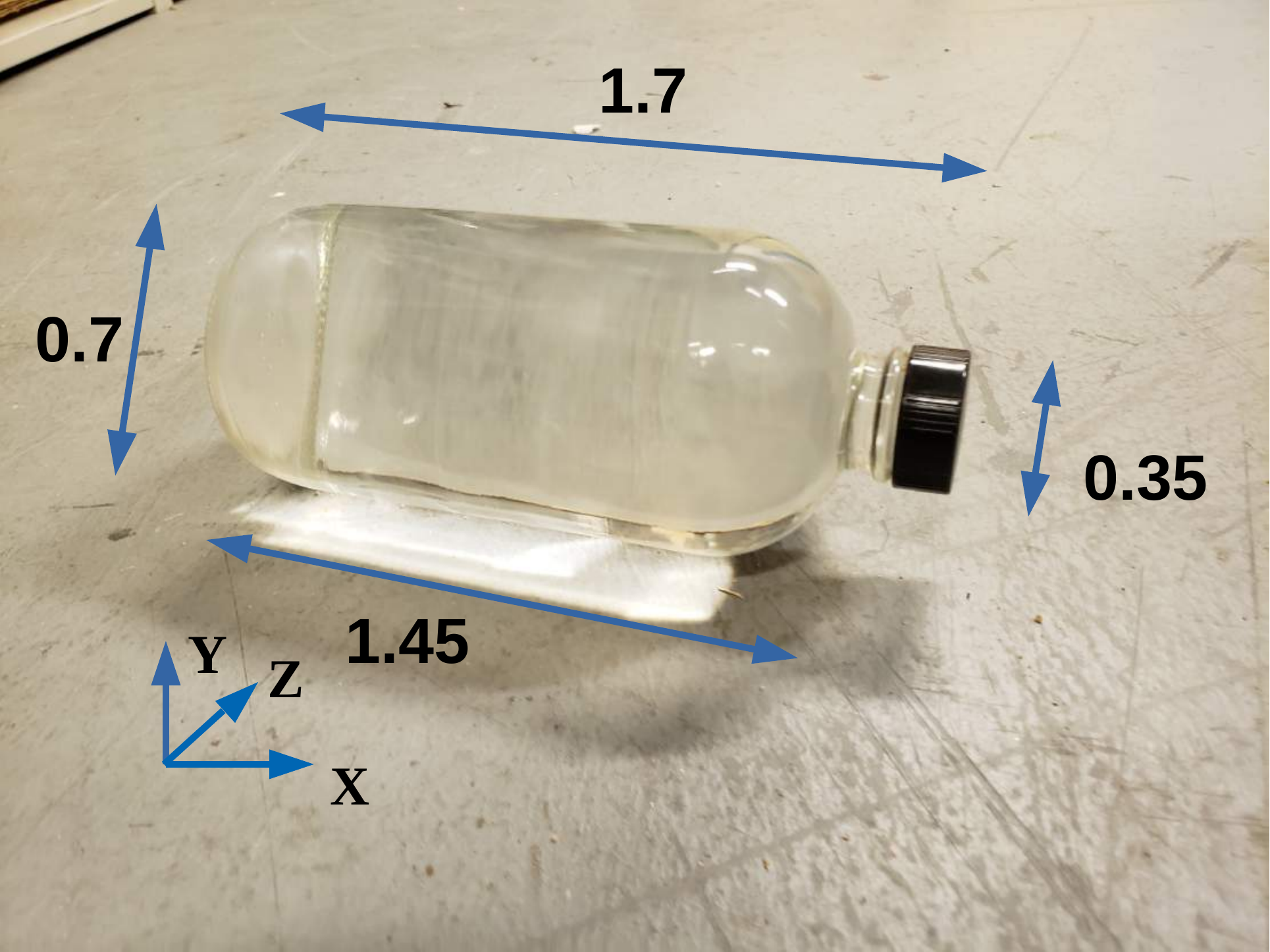}
				\includegraphics[scale=0.4]{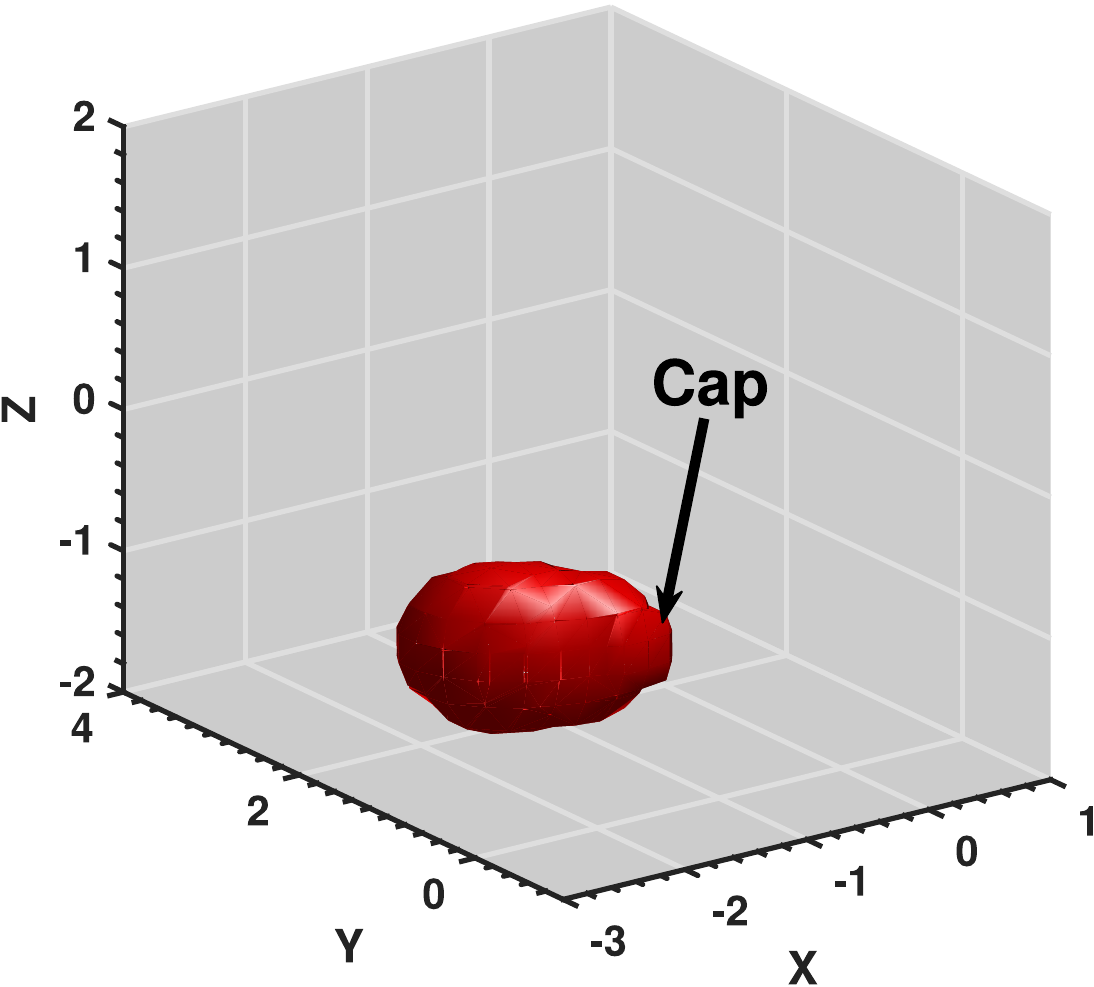}
				\includegraphics[scale=0.4]{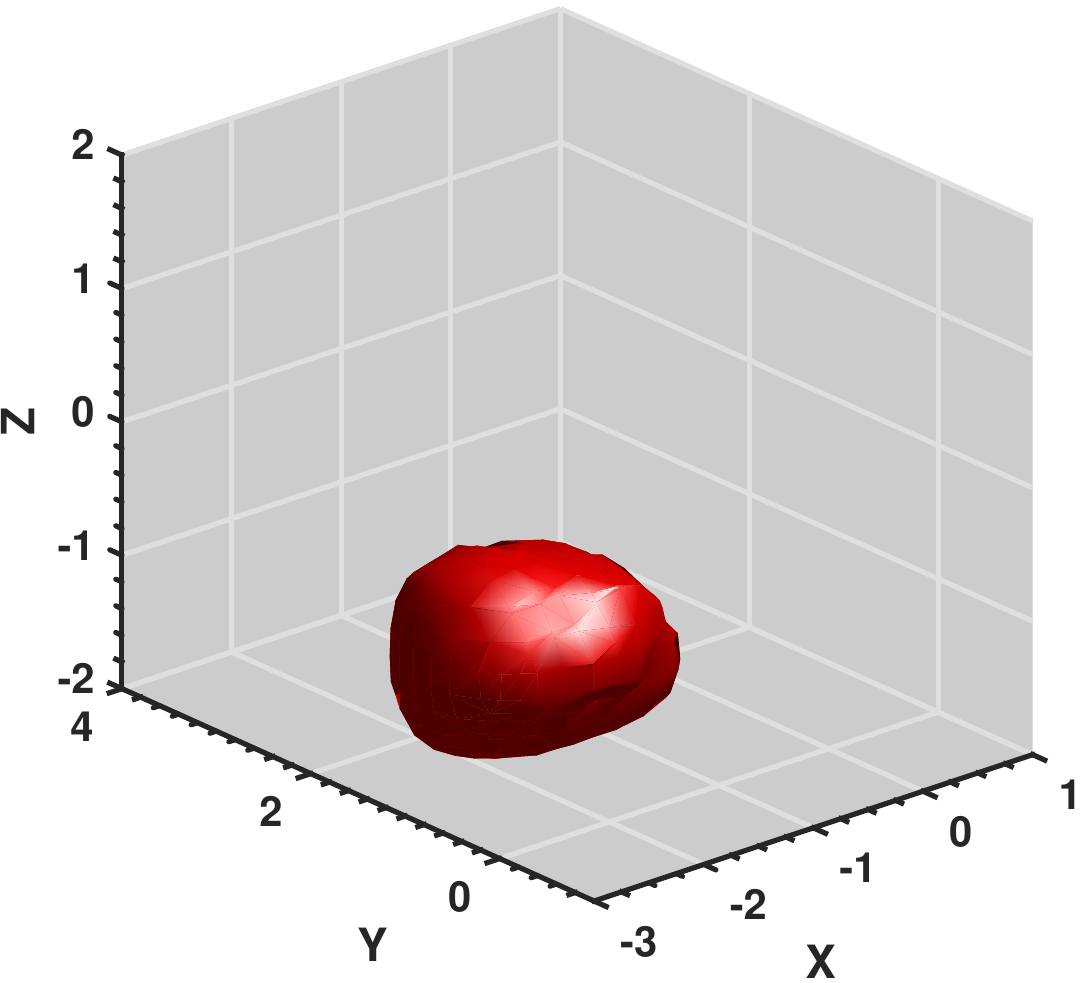}
				\par\end{centering}
		}
		\par\end{centering}\hspace*{\fill}
	\caption{A bottle of clear water; see Tables \ref{table:1}--\ref{table:2} for further details. Note that we can image even a tiny part of it: the cap of this bottle, at least when we compute the dielectric constant. (a) and (b)  indicate a serious data improvement due to the data propagation procedure.}\vspace{-5mm}
\end{figure}

\begin{figure}[H]\vspace{-7mm}
	\begin{centering}\hspace*{\fill}
		\subfloat[Real part of raw and propagated data at $\alpha=0.5$\label{fig:RawProp3}]{\begin{centering}
				\includegraphics[scale=0.3]{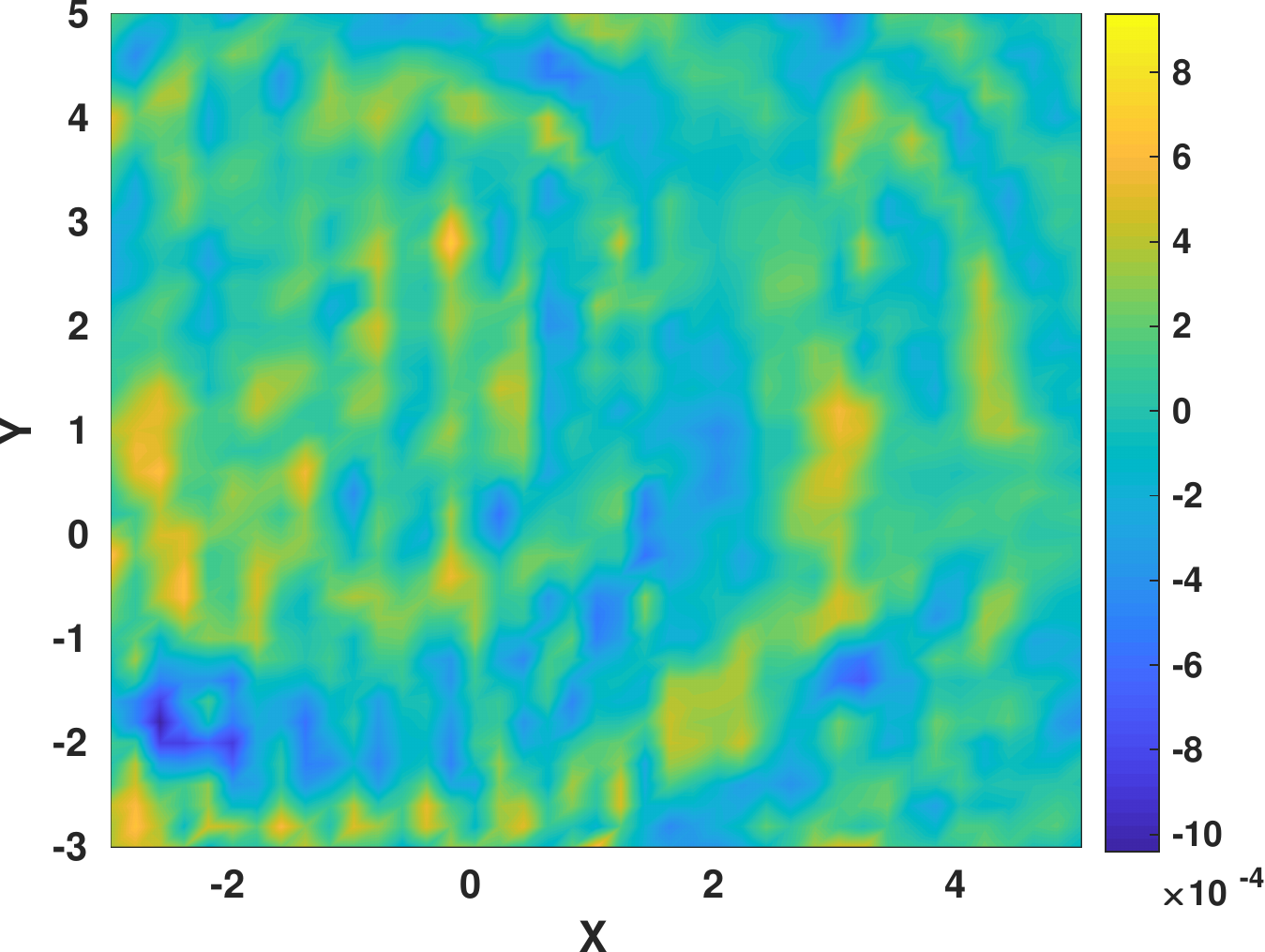}
				\includegraphics[scale=0.3]{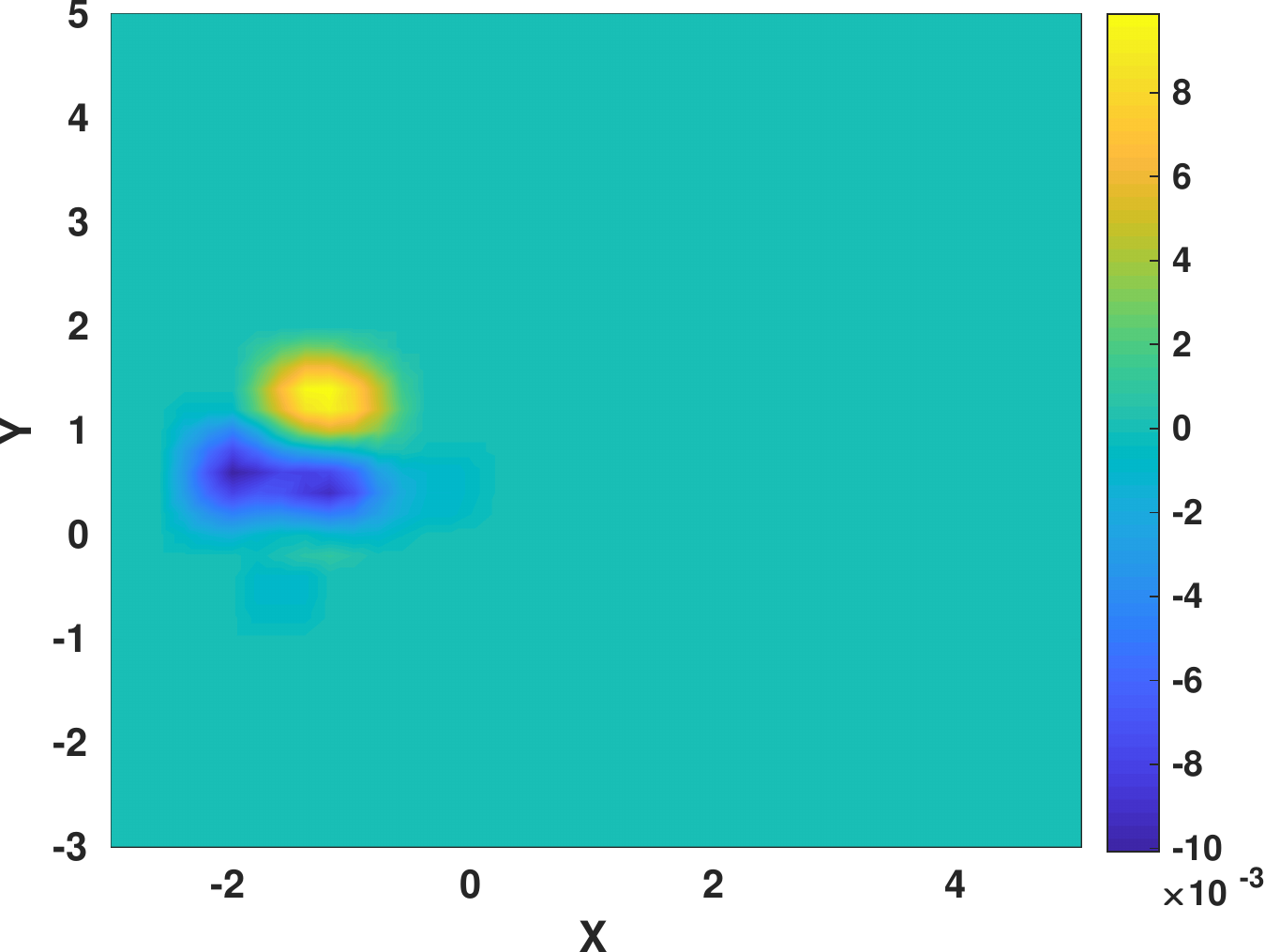}
				\par\end{centering}
		}\subfloat[Imaginary part of raw and propagated data at $\alpha=0.5$\label{fig:imgProp3}]{\begin{centering}
				\includegraphics[scale=0.3]{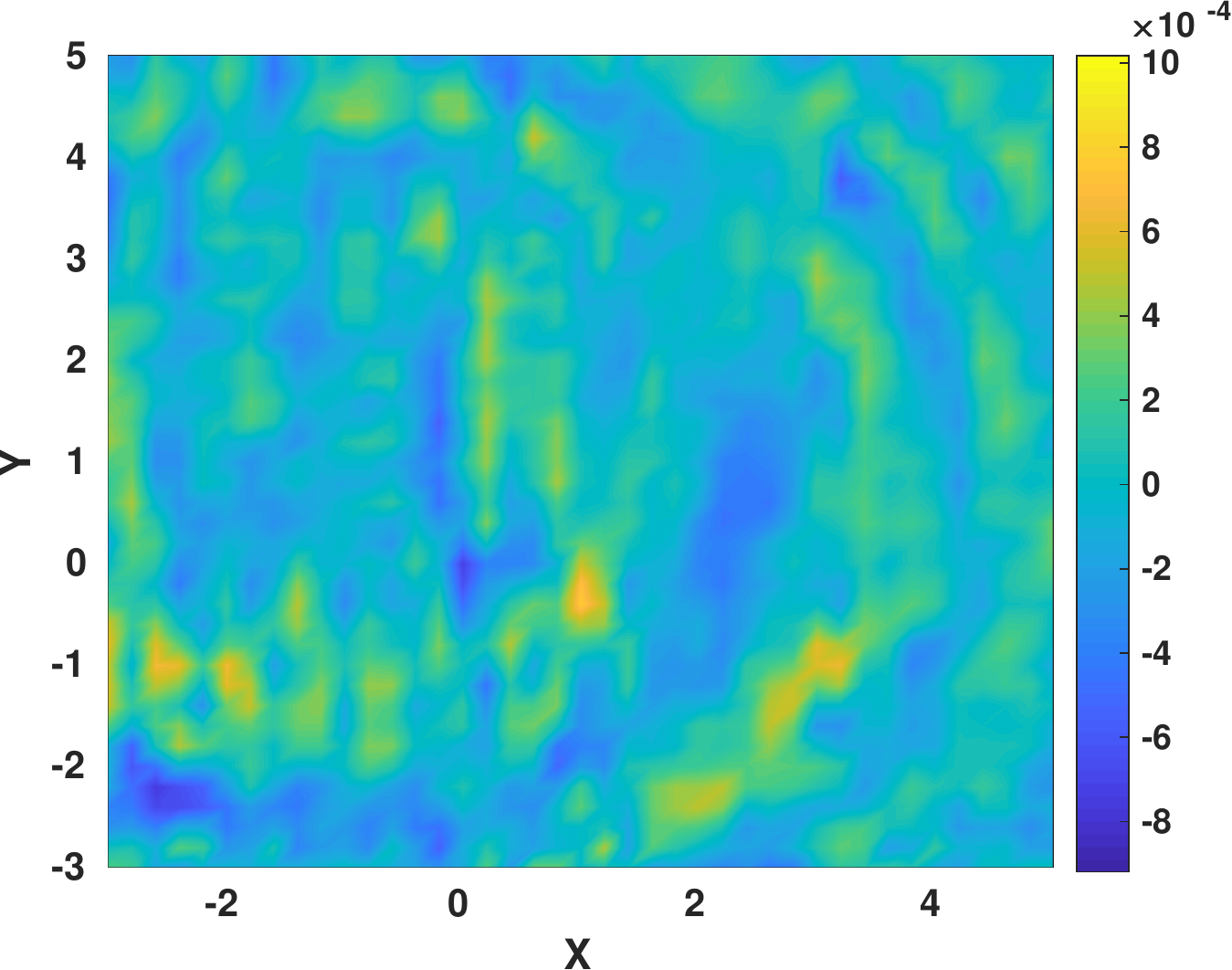}
				\includegraphics[scale=0.3]{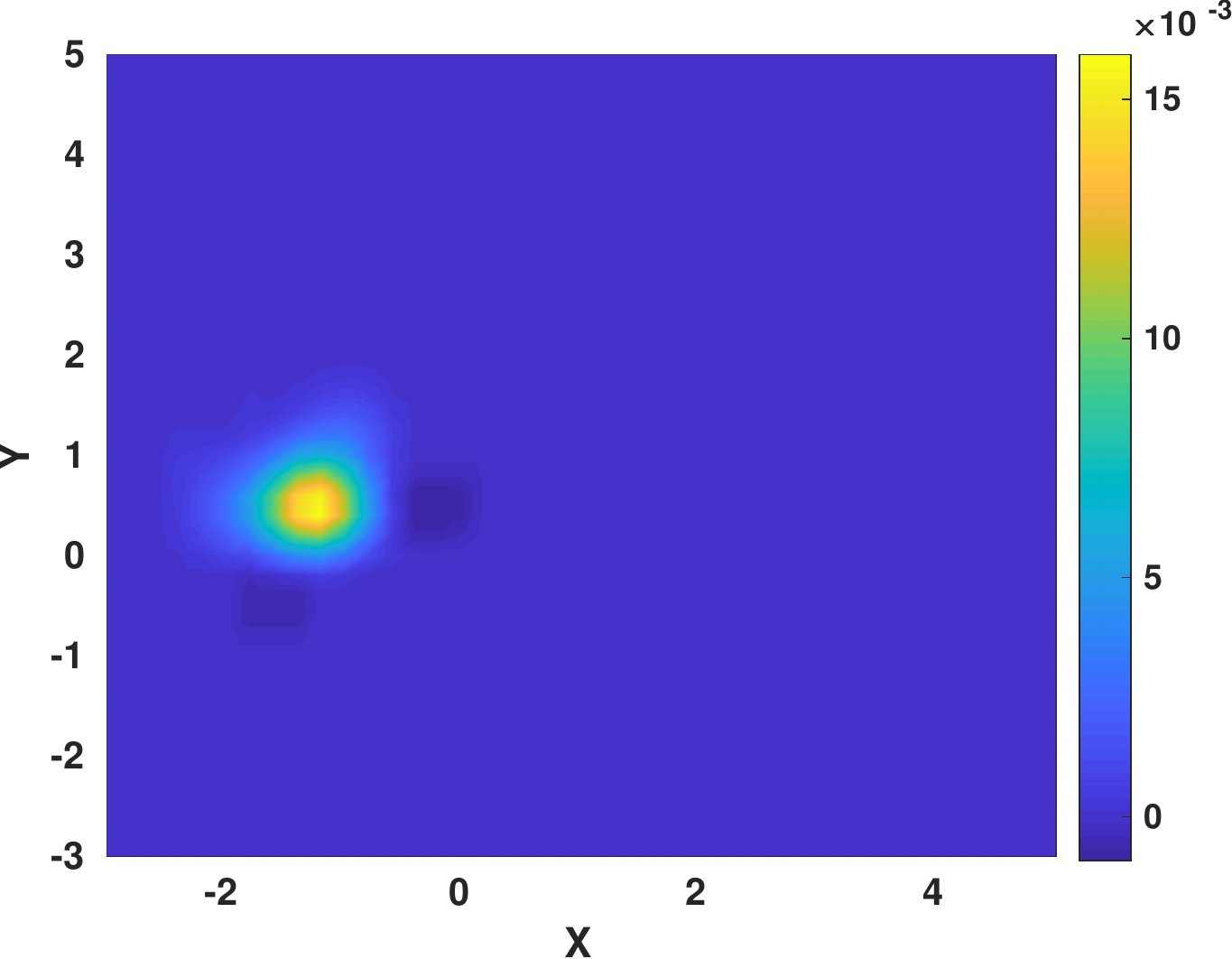}
				\par\end{centering}
		}
		\par\end{centering}
	\begin{centering}\vspace*{\fill}
		\subfloat[Left: U-shaped piece of dry wood (cf. \cite{Khoa2020}). Middle: Image of computed dielectric constant. Right: Image of computed conductivity\label{fig:Walnut}]{ \begin{centering}
				\includegraphics[scale=0.25]{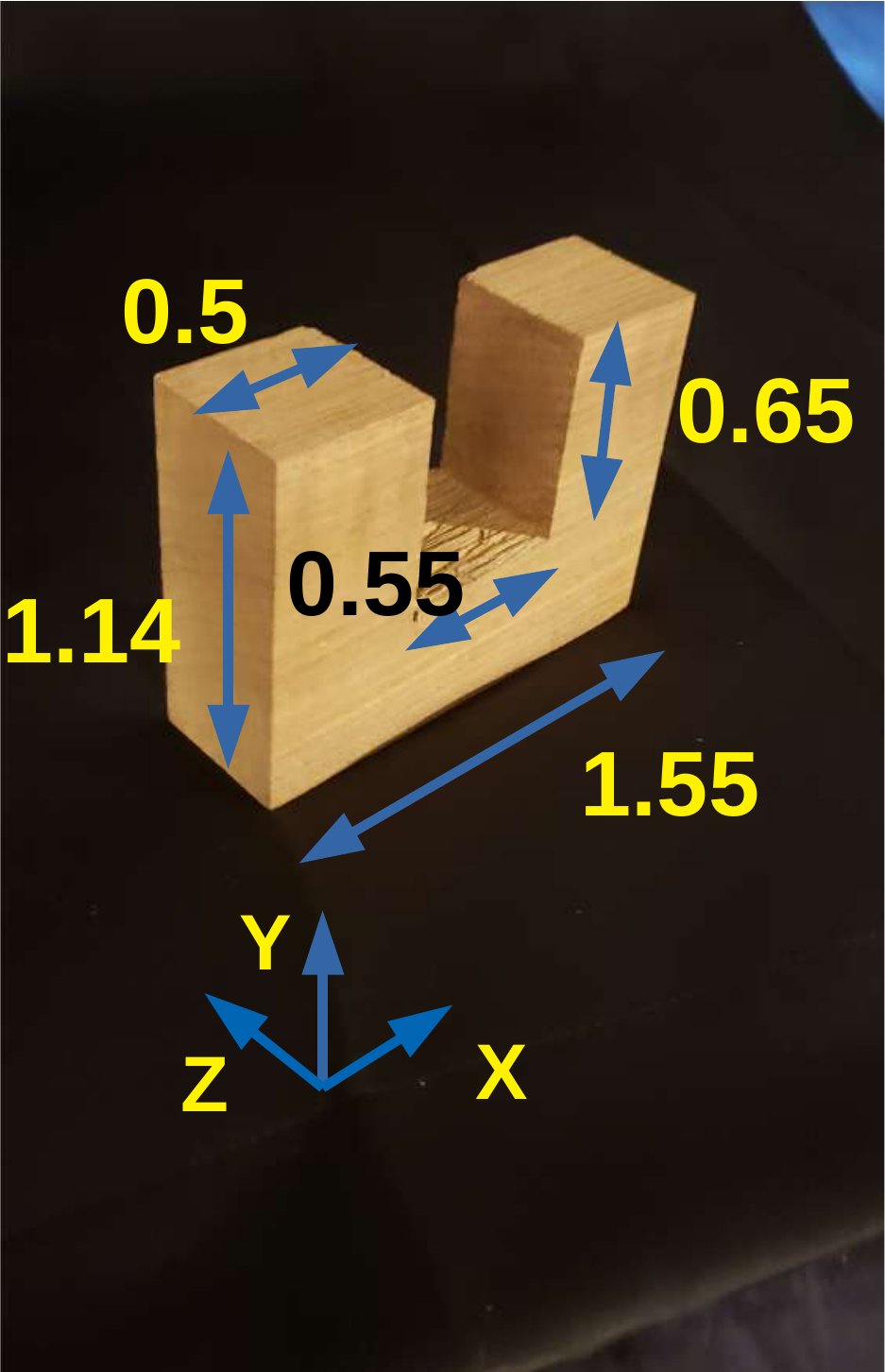}
				\includegraphics[scale=0.4]{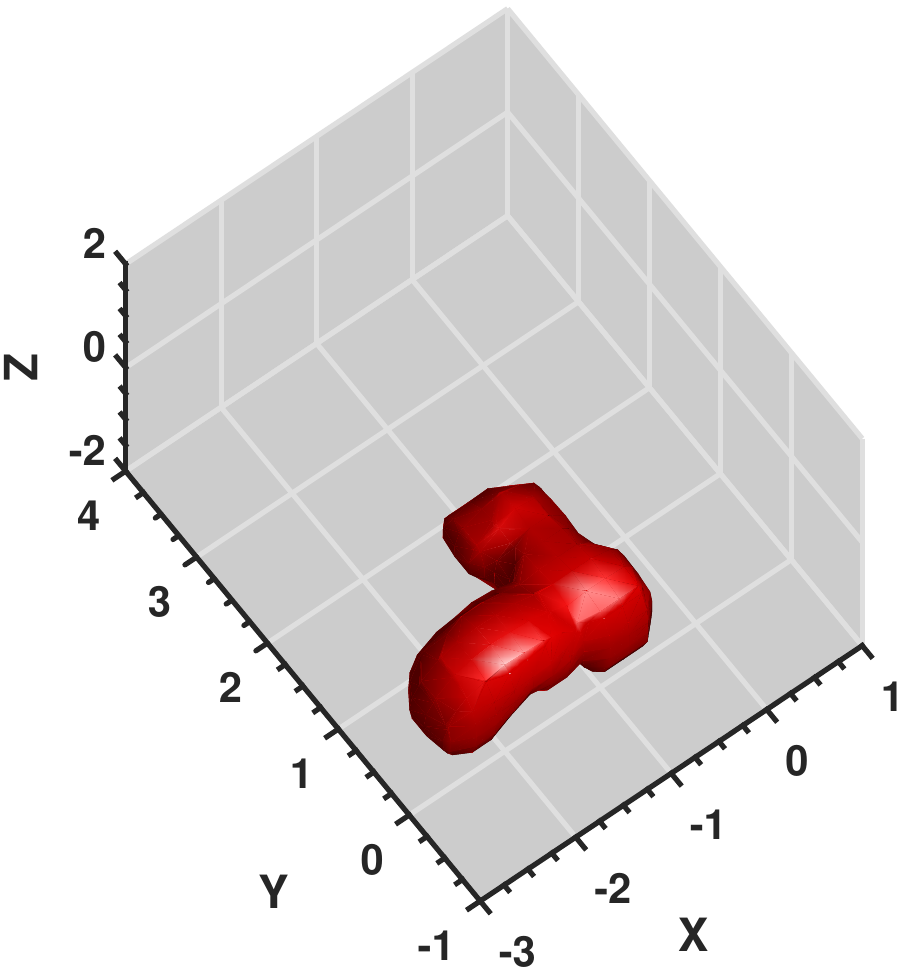}
				\includegraphics[scale=0.4]{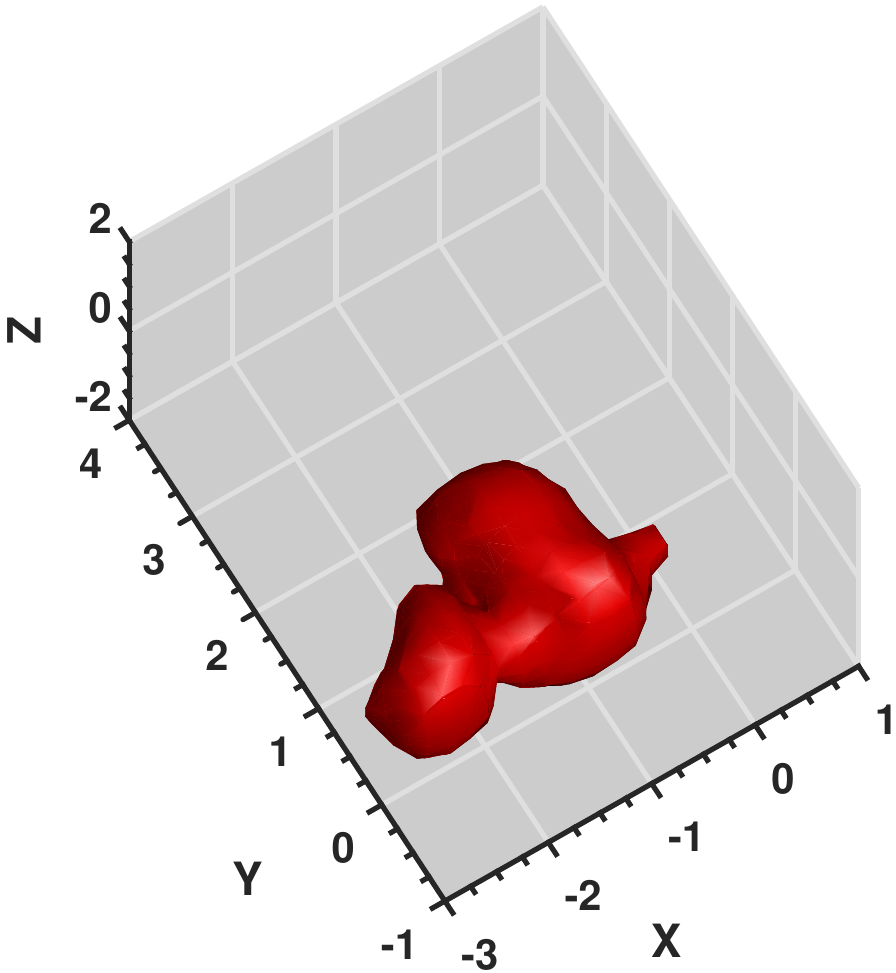}
				\par\end{centering}
		}
		\par\end{centering}\hspace*{\fill}
	\caption{U-shaped piece of dry wood; see Tables \ref{table:1}--\ref{table:2} for further details. Note that we can image even the void of this nonconvex target. It is well known that imaging of nonconvex targets with voids in them is a quite challenging goal. This is especially true for the most challenging case we consider: backscattering nonoverdetermined data.\label{fig:UU}}\vspace{-5mm}
\end{figure}

\begin{figure}[H]\vspace{-7mm}
	\begin{centering}\hspace*{\fill}
		\subfloat[Real part of raw and propagated data at $\alpha=0.2$\label{fig:RawProp4}]{\begin{centering}
				\includegraphics[scale=0.3]{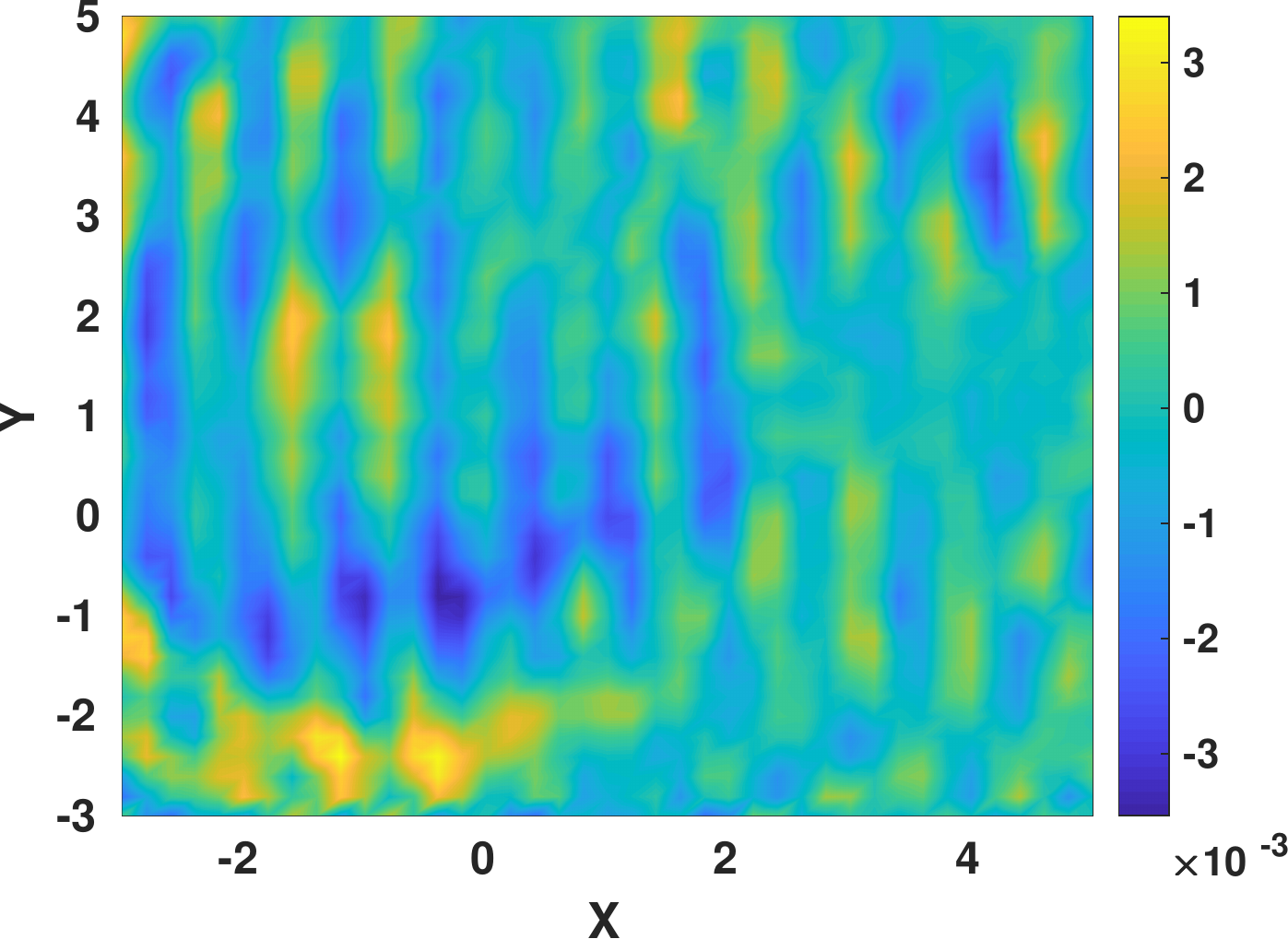}
				\includegraphics[scale=0.3]{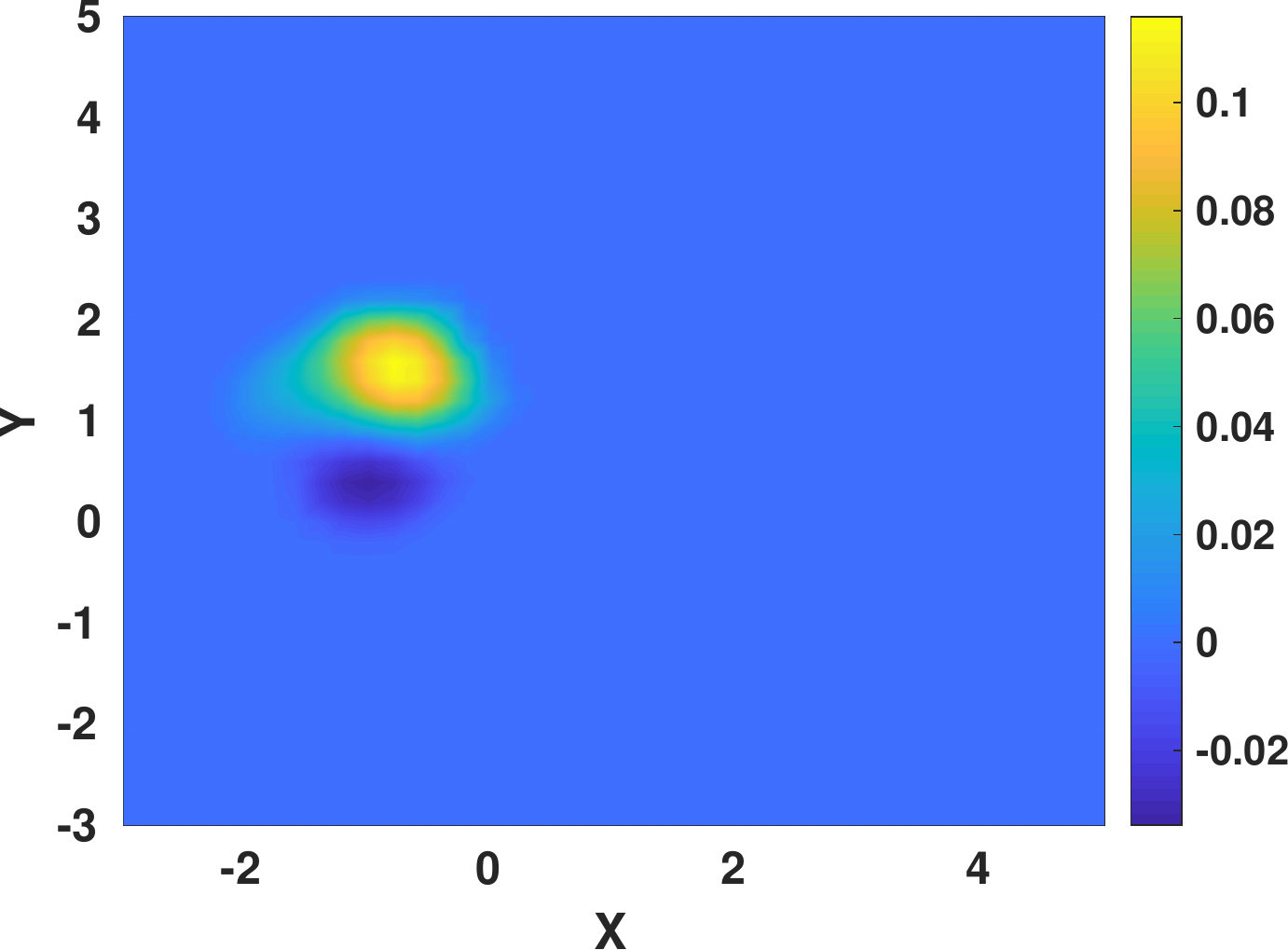}
				\par\end{centering}
		}\subfloat[Imaginary part of raw and propagated data at $\alpha=0.2$\label{fig:imgProp4}]{\begin{centering}
			\includegraphics[scale=0.3]{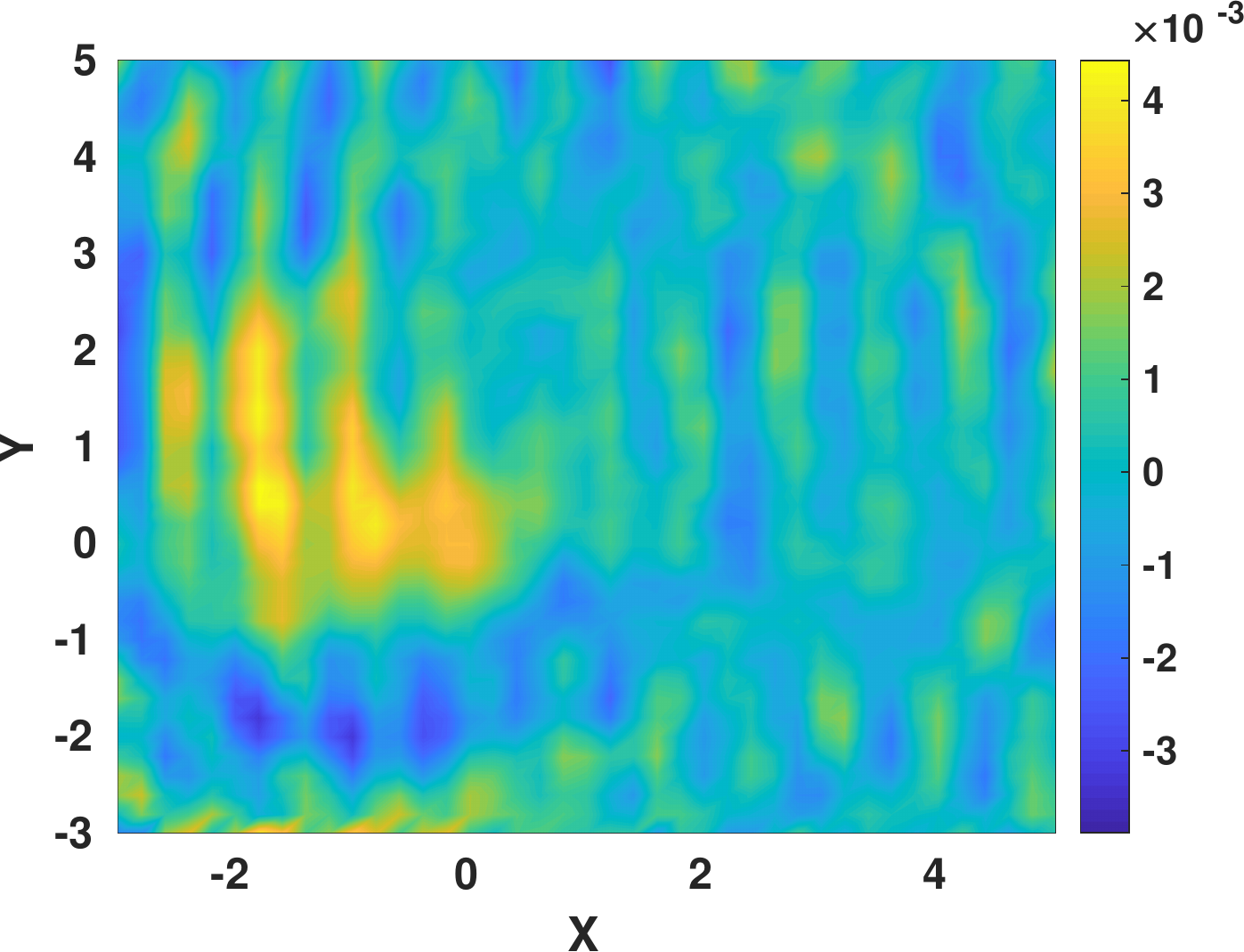}
			\includegraphics[scale=0.3]{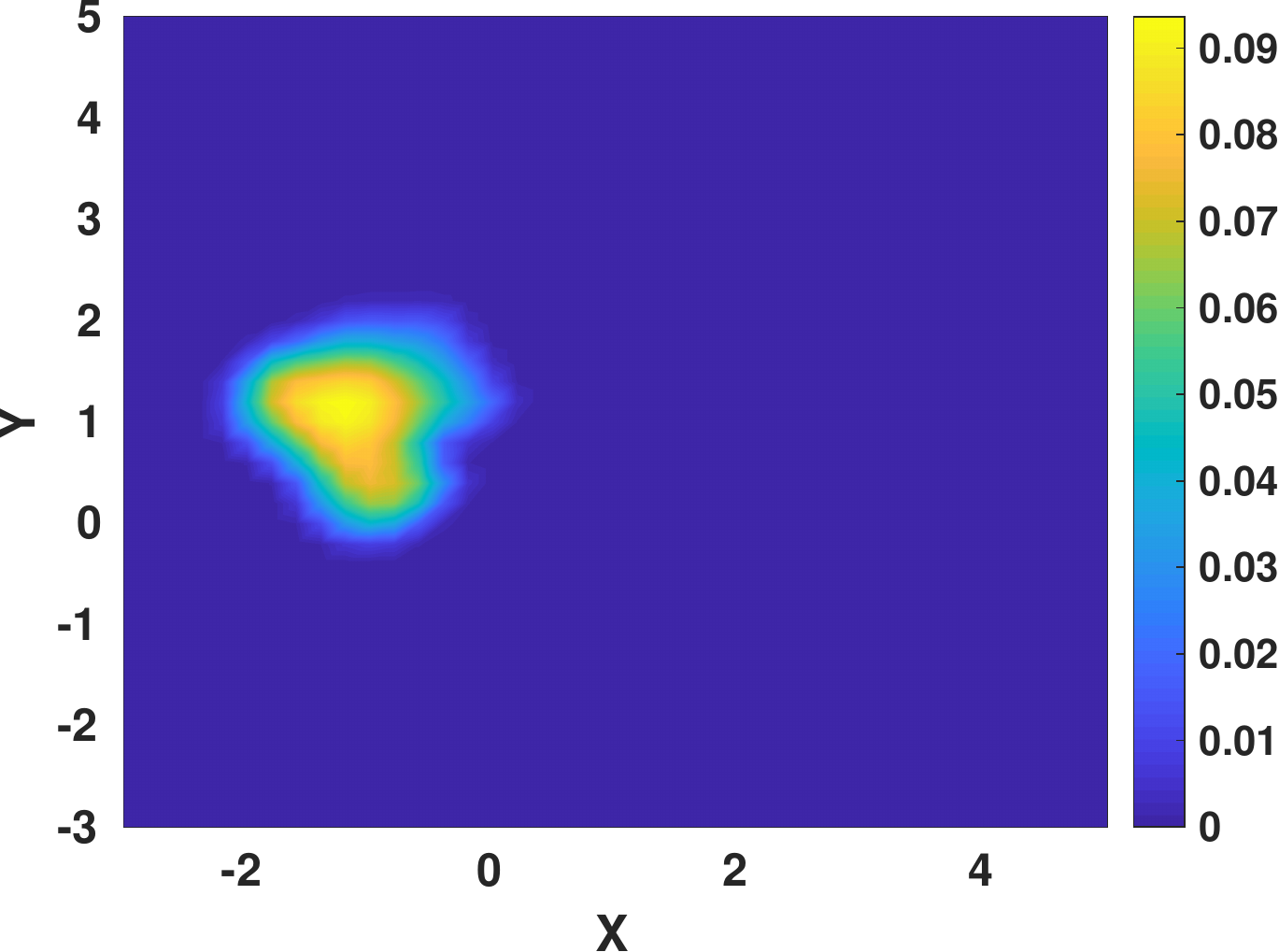}
			\par\end{centering}
		}
		\par\end{centering}
	\begin{centering}\vspace*{\fill}
		\subfloat[Left: Metallic letter ``A'' (cf. \cite{Khoa2020}). Middle: Image of computed dielectric constant. Right: Image of computed conductivity\label{fig:A}]{ \begin{centering}
				\includegraphics[scale=0.2]{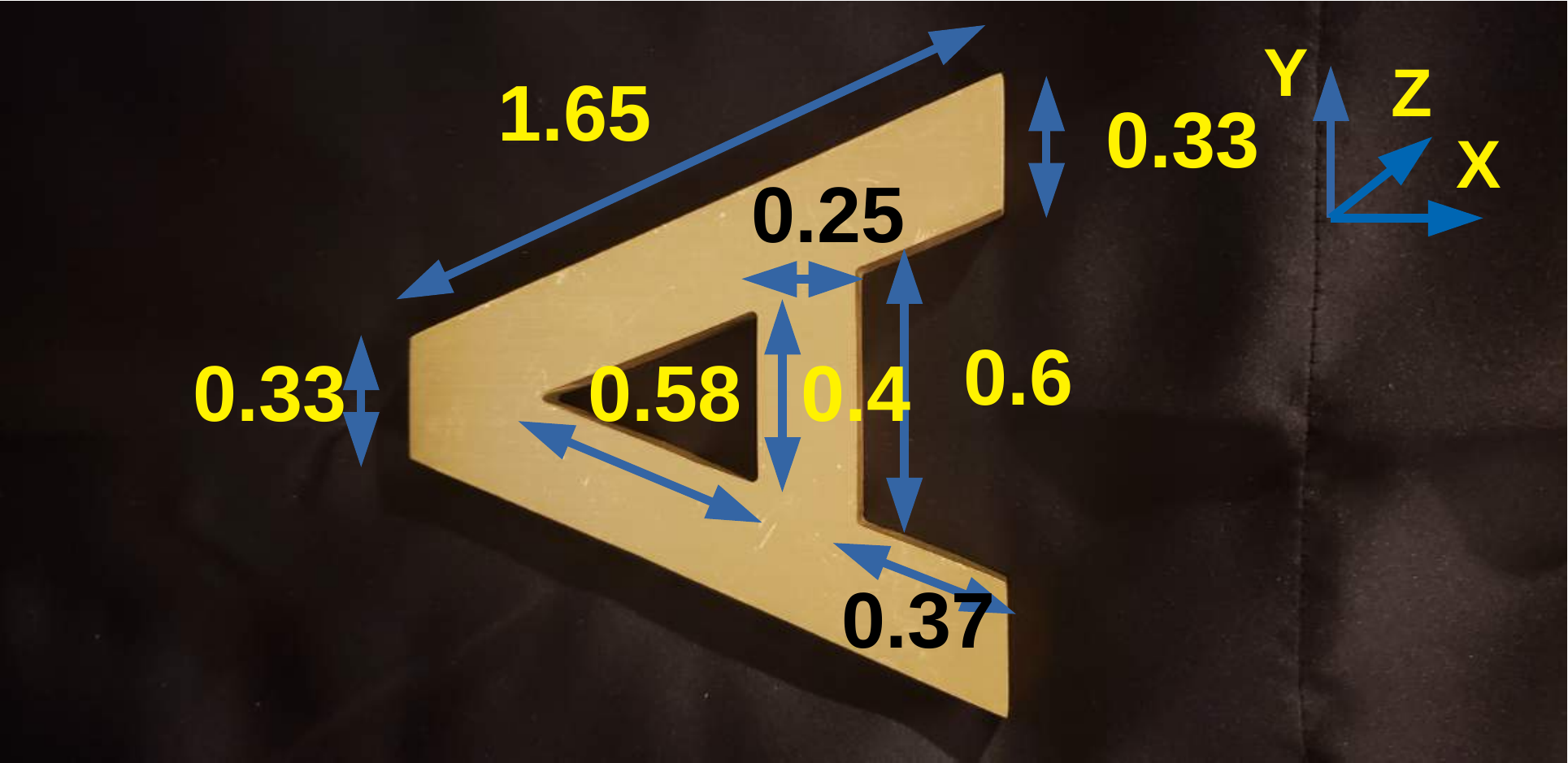}
				\includegraphics[scale=0.4]{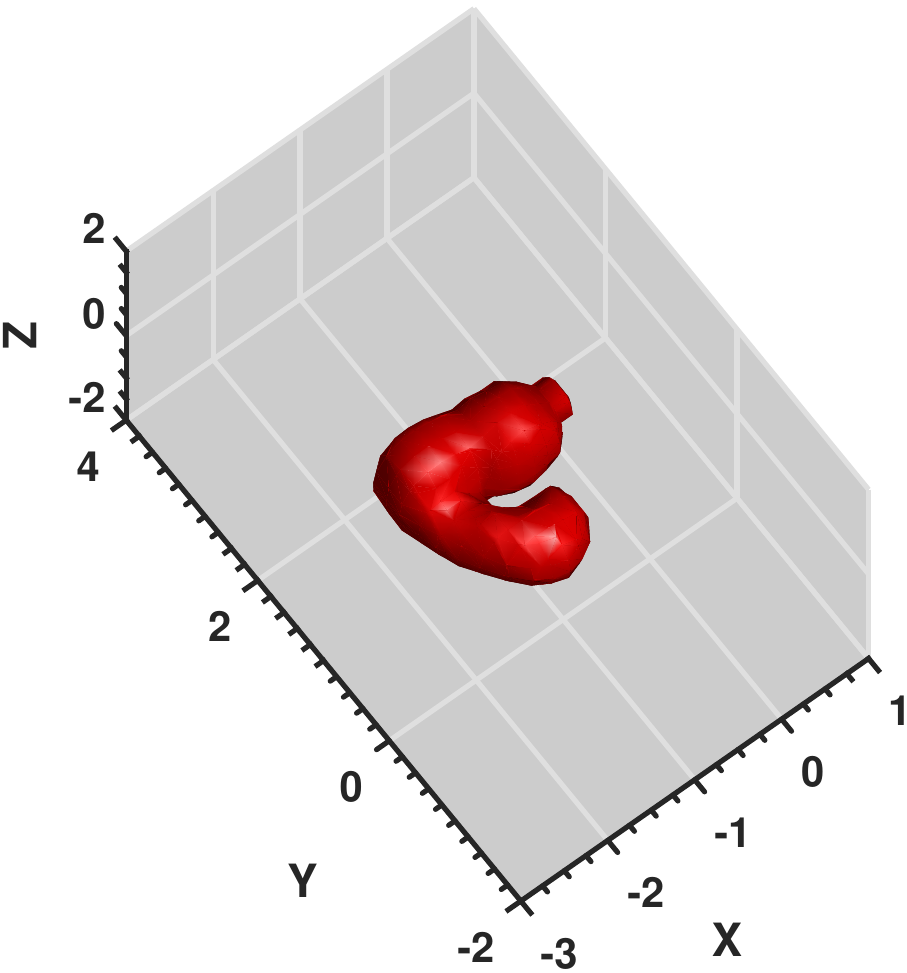}
				\includegraphics[scale=0.4]{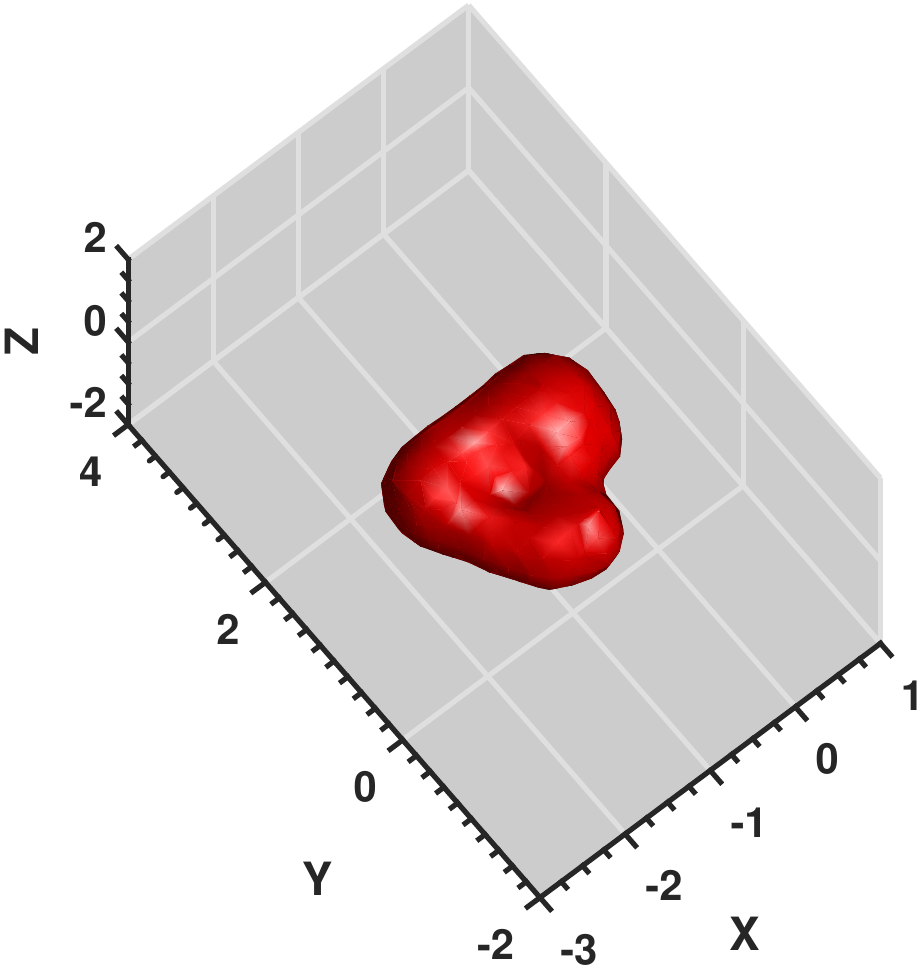}
				\par\end{centering}
		}
		\par\end{centering}\hspace*{\fill}
	\caption{A-shaped metallic target. The same comments as ones for Figure \ref{fig:UU} are applicable here.}\vspace{-5mm}
\end{figure}

\begin{figure}[H]\vspace{-7mm}
	\begin{centering}
		\subfloat[Real part of raw and propagated data at $\alpha=0.6$\label{fig:RawProp5}]{\begin{centering}
				\includegraphics[scale=0.3]{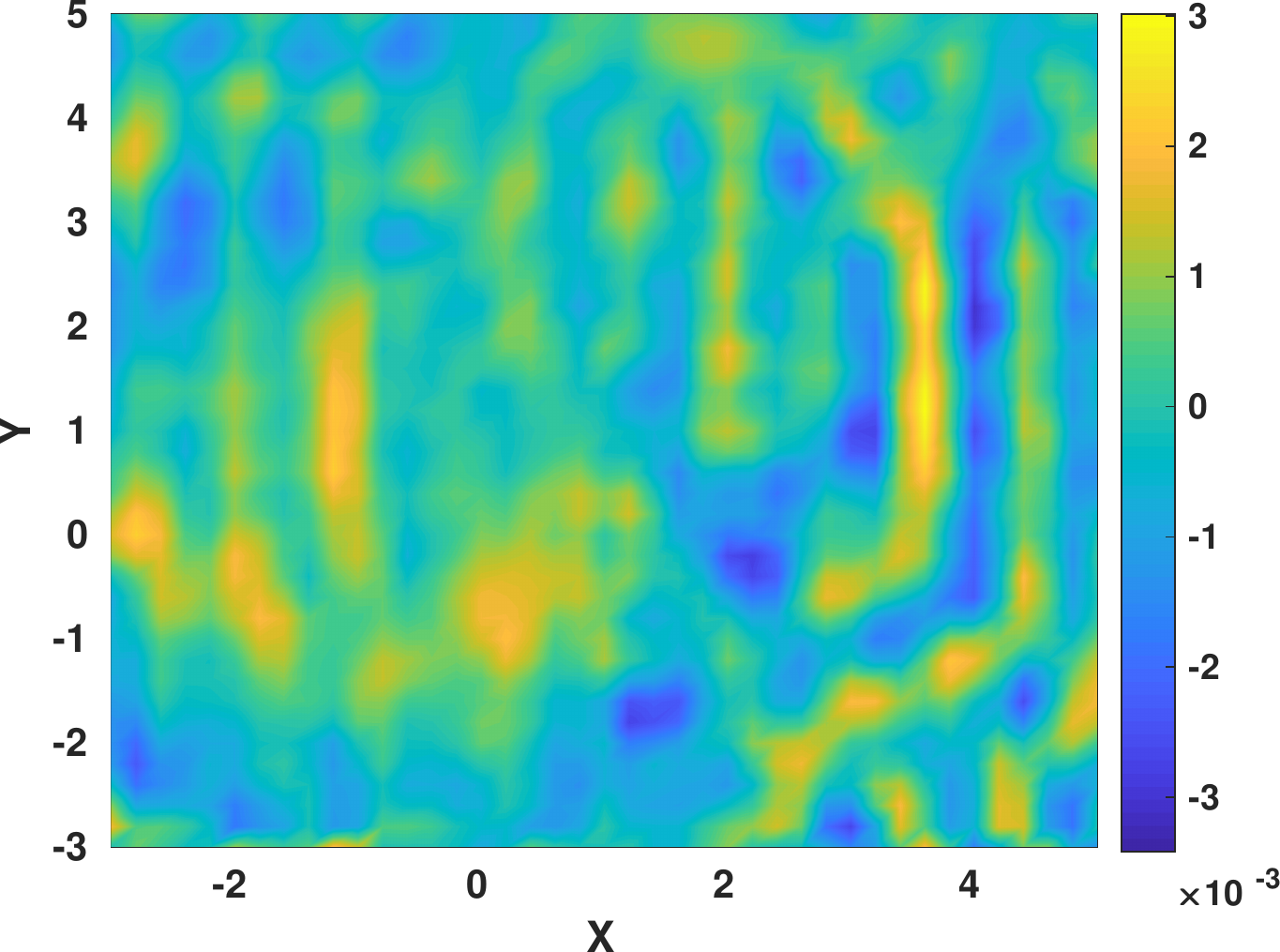}
				\includegraphics[scale=0.3]{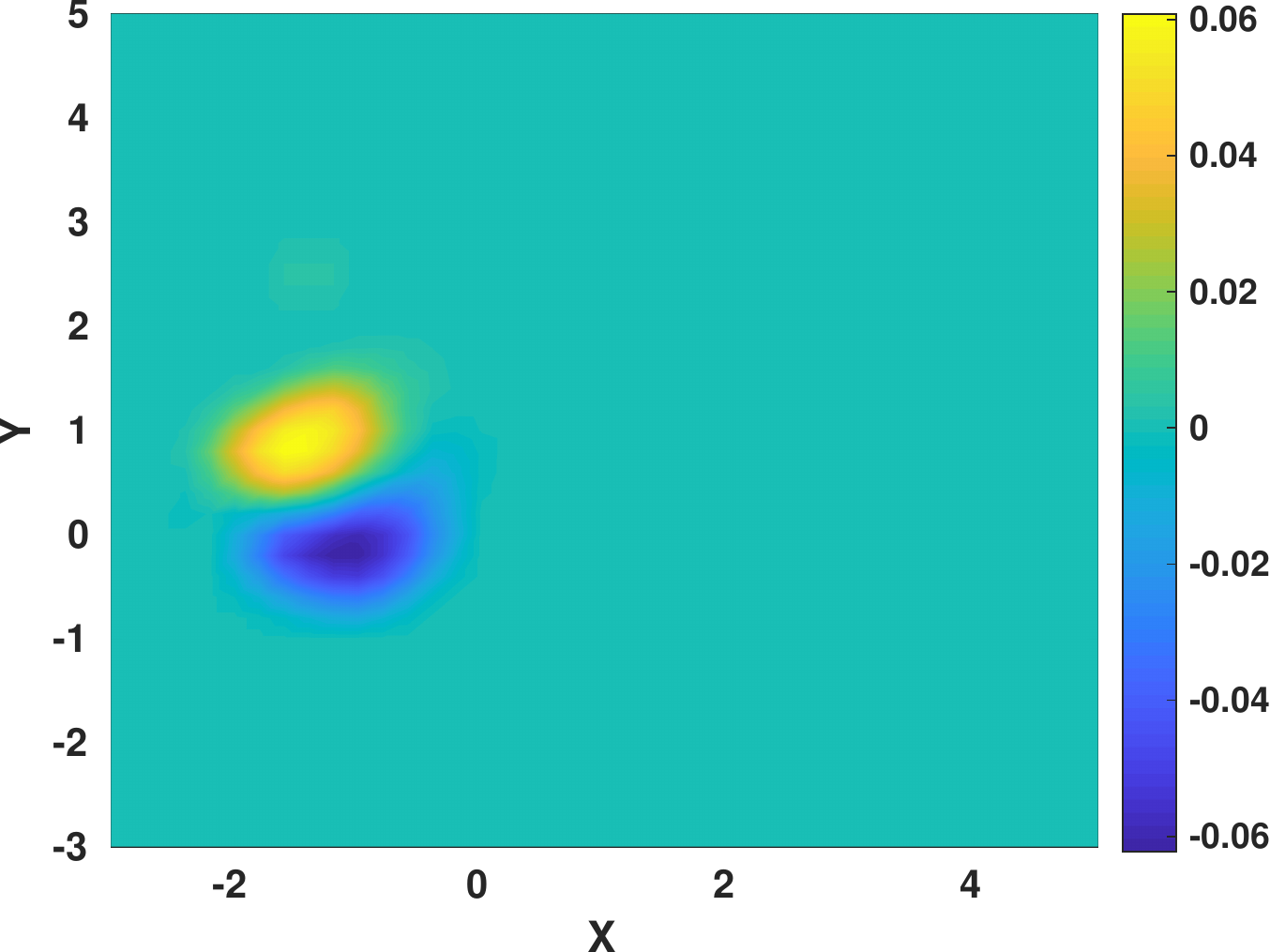}
				\par\end{centering}
		}\subfloat[Imaginary part of raw and propagated data at $\alpha=0.6$\label{fig:imgProp5}]{\begin{centering}
				\includegraphics[scale=0.3]{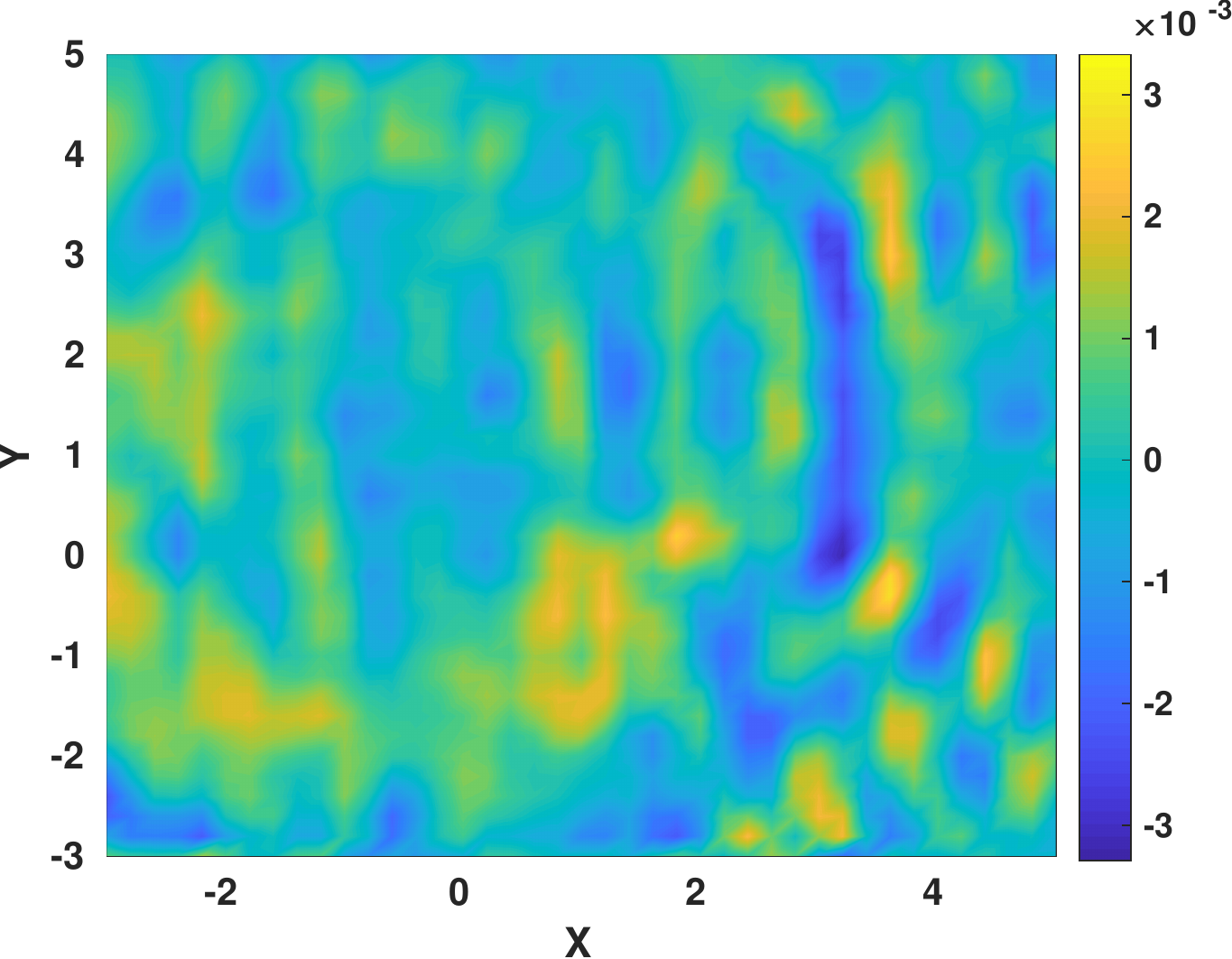}
				\includegraphics[scale=0.3]{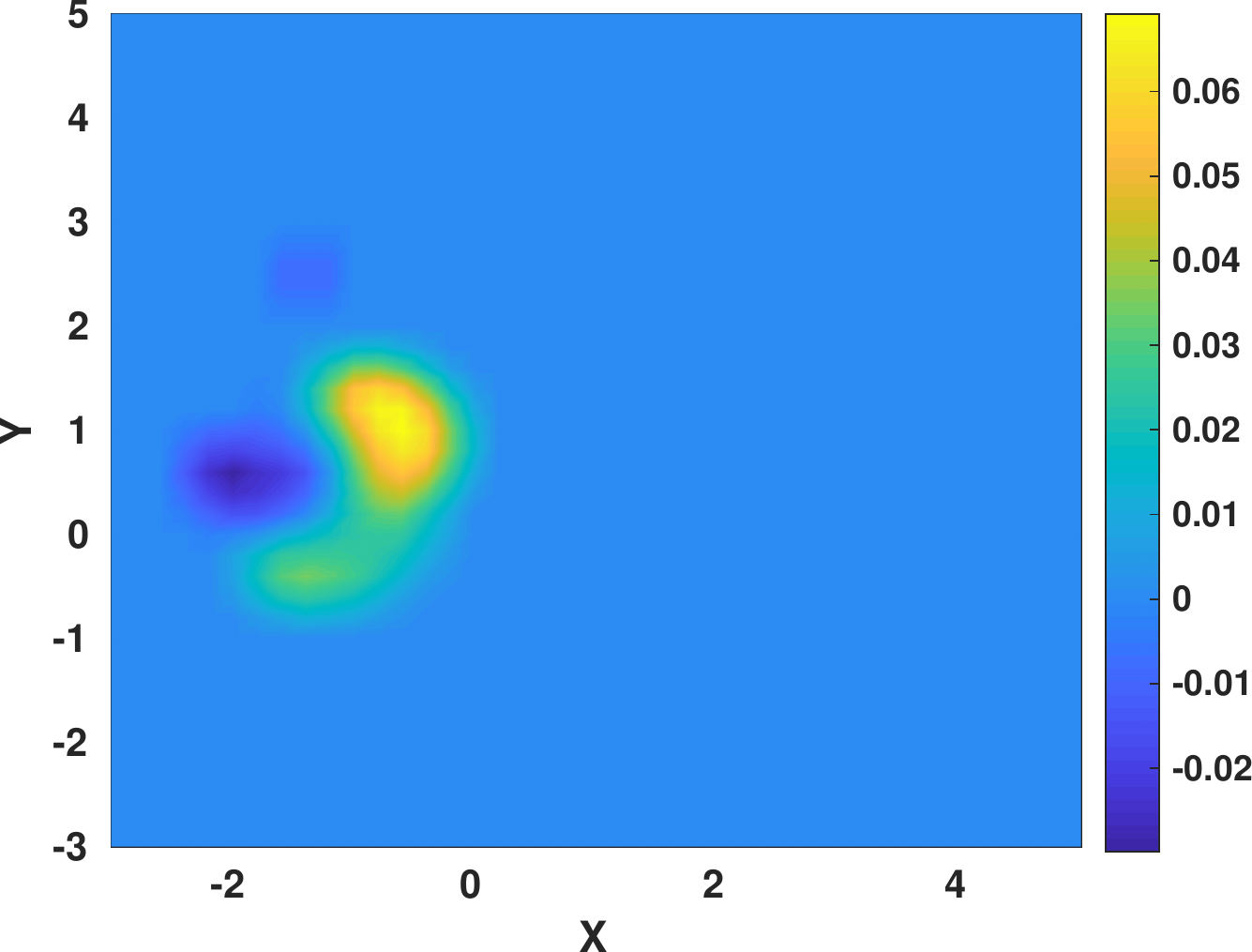}
				\par\end{centering}
		}
		\par\end{centering}
	\begin{centering}
		\subfloat[Left: Metallic letter ``O'' (cf. \cite{Khoa2020}). Middle: Image of computed dielectric constant. Right: Image of computed conductivity\label{fig:O}]{ \begin{centering}
				\includegraphics[scale=0.27]{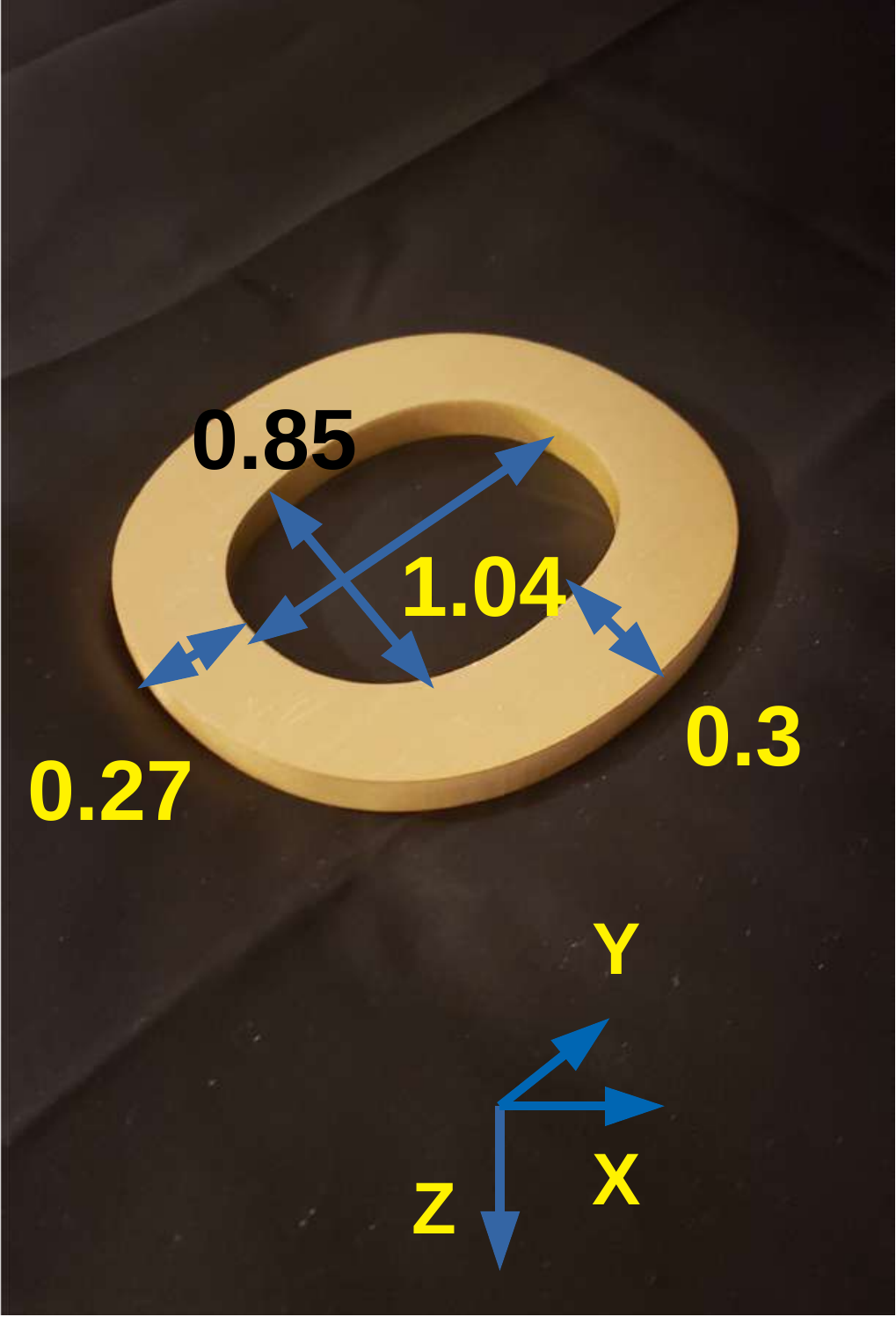}
				\includegraphics[scale=0.4]{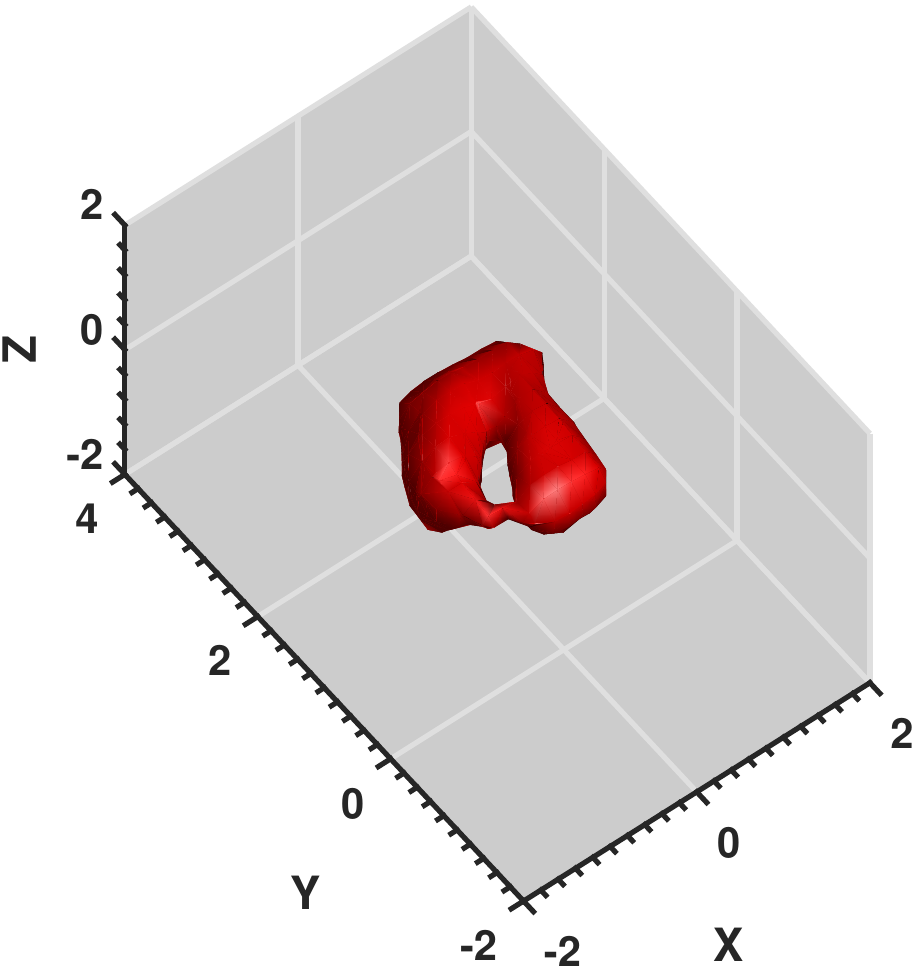}
				\includegraphics[scale=0.4]{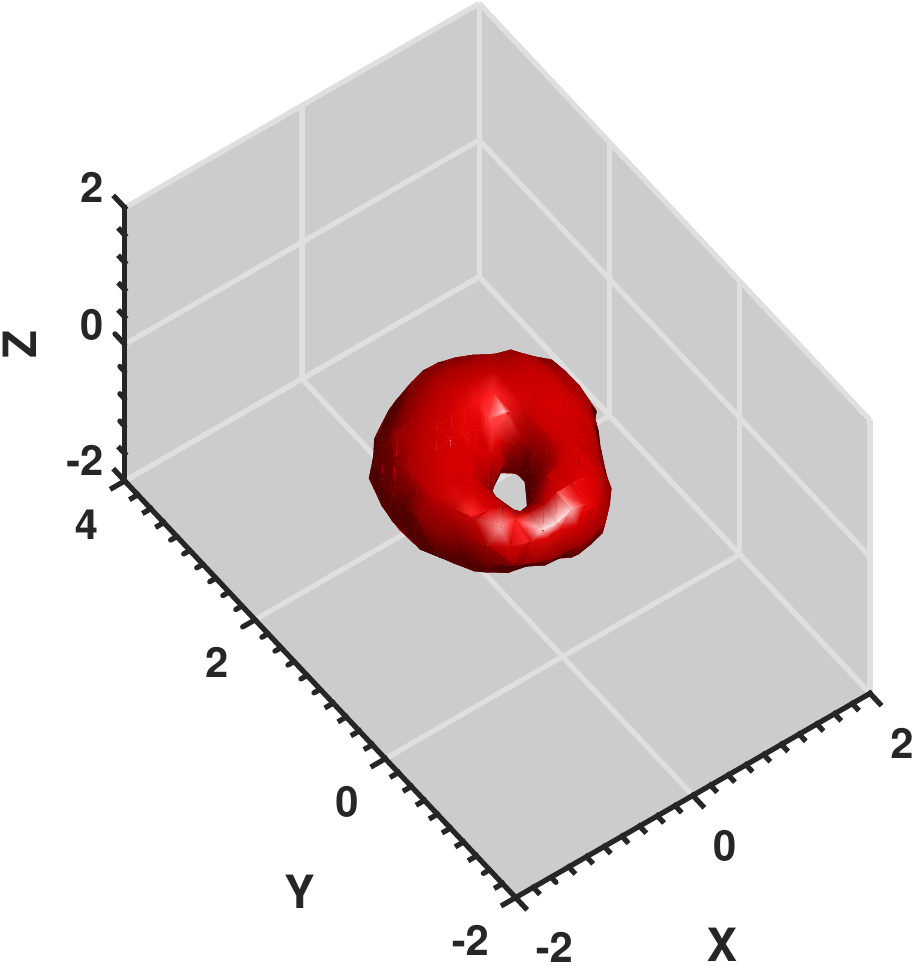}
				\par\end{centering}
		}
		\par\end{centering}\hspace*{\fill}
	\caption{O-shaped metallic target. The same comments as ones for Figure \ref{fig:UU} are applicable here.}\vspace{-5mm}
\end{figure}

\subsection{Comments}

In this work, we do not report the sizes of computed inclusions because this
was done in our work \cite{Khoa2020}. Instead, we apply the following
criterion to distinguish between conductive and non-conductive materials:
the imaged target is conductive if $\max \{\sigma _{p,q,s}\}>1$ and it is
non-conductive otherwise. However, cf. \cite[Section 2.5]{Balanis2012}, this
criterion cannot distinguish the intrinsic semiconductors and the
insulators, whose conductivities are all less than the unity. At this
moment, our Table \ref{table:2} shows that the algorithm can detect well the
conductors, which are usually metallic. Therefore, our algorithm might
potentially be helpful to decrease the number of false alarms.

We use the \texttt{isosurface} built-in function in MATLAB to depict 3D
images of computed inclusions. The associated \texttt{isovalue} is 10\% of
the maximal value of the inclusions.

We now briefly comment on shapes of imaged inclusions. Figure 3 shows that
we can image even a tiny part of that bottle: its cap, at least when we
image the dielectric constant. Images of Figures 4-6 show that we can image
non-convex targets and even voids inside of them. Obviously such features
are tough to image, given only limited backscattering and non-overdetermined
data.

\textbf{Acknowledgments.} The first author acknowledges Prof. Dr. Dinh-Liem Nguyen (Kansas, USA) for the fruitful discussions on the experimental setup.



\end{document}